\documentclass[reqno]{amsart}

\usepackage{amsfonts,amssymb,amsmath,mathtools,mathrsfs}
\usepackage{lmodern}
\usepackage[english]{babel}
\usepackage[utf8]{inputenc}
\usepackage[a4paper,margin=2.2cm]{geometry}
\usepackage{xcolor}
\usepackage{tikz}
\usepackage{cite}

\definecolor{RefBlue}{RGB}{0,0,180}
\definecolor{CiteGreen}{RGB}{0,0,180}
\definecolor{UrlBlue}{RGB}{0,0,180}
\usepackage[
  colorlinks=true,
  linkcolor=RefBlue,
  citecolor=CiteGreen,
  urlcolor=UrlBlue,
  filecolor=UrlBlue,
  bookmarksnumbered=true,
  bookmarksopen=true,
  linktoc=all,
  breaklinks=true,
  pdfborder={0 0 0},
  hypertexnames=false
]{hyperref}

\newtheorem{maintheorem}{Theorem}

\newtheorem{theorem}{Theorem}[section]
\newtheorem{proposition}[theorem]{Proposition}
\newtheorem{lemma}[theorem]{Lemma}
\newtheorem{corollary}[theorem]{Corollary}
\theoremstyle{definition}

\theoremstyle{remark}
\newtheorem{remark}[theorem]{Remark}

\numberwithin{equation}{section}
\DeclareMathOperator{\Div}{div}

\newcommand{\eps}{\varepsilon}
\newcommand{\Ystar}{Y^{*}}
\newcommand{\Qeps}{Q^{\eps}}
\newcommand{\Om}{\Omega}
\newcommand{\Omeps}{\Omega^{\eps}}
\newcommand{\dd}{\,d}
\newcommand{\R}{\mathbb R}
\newcommand{\weakto}{\rightharpoonup}
\newcommand{\tto}{\mathrel{\mathop{\longrightarrow}\limits^{2\mathrm{s}}}}
\newcommand{\twto}{\mathrel{\mathop{\rightharpoonup}\limits^{2}}}

\begin{document}

\title[Degenerate parabolic Signorini homogenization]{Homogenization of
 degenerate parabolic Signorini problems with vanishing storage and
 oscillating surface fluxes}

\author[Hermann Douanla]{Hermann Douanla$^{*}$}
\author[Rodrigue Kenne Bogning]{Rodrigue Kenne Bogning$^{\dagger}$}
\author[Erick Tetsadjio]{Erick Tetsadjio$^{\ddagger}$}
\thanks{\raggedright$^{*}$Department of Mathematics, University of Yaounde I, P.O. Box 812,
Yaounde, Cameroon.  E-mail:
\href{mailto:hermann.douanla@facsciences-uy1.cm}{\texttt{hermann.douanla@facsciences-uy1.cm}}.}
\thanks{\raggedright\hspace*{\normalparindent}$^{\dagger}$Department of Fundamental and Transversal Sciences, National
Advanced School of Public Works, P.O. Box 510, Yaounde, Cameroon.  E-mail:
\href{mailto:rodriguekenne@ymail.com}{\texttt{rodriguekenne@ymail.com}}.}
\thanks{\raggedright\hspace*{\normalparindent}$^{\ddagger}$Department of Mathematics, University of Yaounde I, P.O. Box 812,
Yaounde, Cameroon.  E-mail:
\href{mailto:tetsadjio@gmail.com}{\texttt{tetsadjio@gmail.com}}.}

\date{}
\subjclass[2020]{35B27, 35K65, 35K86, 35R35, 47J20, 76S05}
\keywords{Degenerate parabolic variational inequalities; Periodic homogenization;
Signorini boundary conditions; Weighted Sobolev spaces; Oscillating surface
fluxes; Two-scale convergence}

\begin{abstract}
We study the homogenization of a unilateral quasilinear parabolic problem in a
periodically perforated medium.  The model combines two degeneracies: the
coefficient of the time derivative is non-negative and may vanish, while the
directional diffusion weights are positive almost everywhere but need not be
bounded away from zero.  The exterior boundary carries a Signorini condition,
and the holes carry an unscaled zero-mean flux oscillating on the diffusive
time scale.  For any family of microscopic solutions we prove strong
space--time $L^2$ compactness through a weighted detection inequality, which
replaces the Aubin--Lions lemma when the storage vanishes, and we recover
spatial $H^1$ regularity of the limit by a closed-range argument for the
weighted two-scale gradients.  The limit solves a coupled two-scale problem in
Minty form, and microscopic existence holds under joint monotonicity.  Under a
strengthened condition, the limit decouples into a macroscopic Signorini
inequality whose uniqueness follows from an exponential regularization,
without time regularity of the limit, and the whole family converges.  The
surface datum enters both the effective flux and the exterior boundary
condition.
\end{abstract}

\maketitle

\section{Introduction}

Diffusion and saturated--unsaturated flow in periodically perforated porous
materials lead to parabolic problems whose storage or capacity coefficient may
vanish in regions where the stored quantity is negligible; the equation is then
elliptic--parabolic, as for the Richards equation \cite{AltLuckhaus1983}, and
only a positive cell average of the storage preserves an effective time
evolution.  The mobilities may also be direction-dependent and approach zero,
for instance near nearly impermeable inclusions, so that $H^1$ estimates are
not available and traces and compactness have to be recovered from an
integrability condition on the inverse weights.  Finally, the pores may
exchange mass with the matrix through a cyclic flux synchronized with the cell
diffusion time, while semipermeable walls, contact, threshold or one-sided
exchange laws on the exterior boundary are modelled by a Signorini condition
\cite{DuvautLions1972}.  This paper studies the periodic homogenization of the
model combining these features.

\subsection{The model}\label{sec:geometry}

Let $N\geq2$, $T>0$, $Y=(0,1)^N$ and $Z=(0,1)$.
Let $\Om\subset\mathbb R^N$ be a bounded connected Lipschitz domain whose
boundary is the union, up to a set of surface measure zero, of two disjoint
relatively open sets $\Gamma_1$ (Signorini) and $\Gamma_2$ (Dirichlet), with
$|\Gamma_2|_{N-1}>0$.  Let $T_0\Subset Y$ be a nonempty open reference hole with
Lipschitz boundary $S=\partial T_0$, and assume that
$\Ystar=Y\setminus\overline{T_0}$ is connected.  On $S$, let $n$ denote the
outward unit normal to $\Ystar$, pointing into $T_0$, and let $\nu$ denote
the outward unit normal to $\Om$.

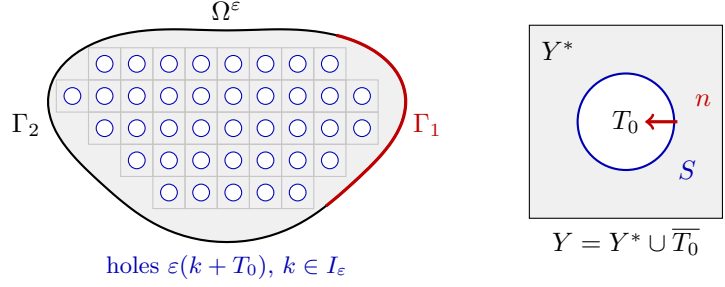
\begin{figure}[!htbp]
\centering
\begin{tikzpicture}[scale=0.85]
 \fill[gray!12] (5.403,1.900) -- (5.438,1.998) -- (5.463,2.099) -- (5.477,2.202) -- (5.478,2.305) -- (5.465,2.408) -- (5.439,2.511) -- (5.399,2.611) -- (5.345,2.708) -- (5.277,2.801) -- (5.195,2.889) -- (5.102,2.971) -- (4.997,3.046) -- (4.883,3.113) -- (4.760,3.173) -- (4.631,3.225) -- (4.496,3.269) -- (4.358,3.305) -- (4.217,3.333) -- (4.077,3.355) -- (3.937,3.370) -- (3.799,3.380) -- (3.663,3.385) -- (3.532,3.387) -- (3.403,3.386) -- (3.279,3.383) -- (3.158,3.379) -- (3.041,3.376) -- (2.926,3.373) -- (2.812,3.371) -- (2.700,3.370) -- (2.588,3.371) -- (2.474,3.373) -- (2.359,3.376) -- (2.242,3.379) -- (2.121,3.383) -- (1.997,3.386) -- (1.868,3.387) -- (1.737,3.385) -- (1.601,3.380) -- (1.463,3.370) -- (1.323,3.355) -- (1.183,3.333) -- (1.042,3.305) -- (0.904,3.269) -- (0.769,3.225) -- (0.640,3.173) -- (0.517,3.113) -- (0.403,3.046) -- (0.298,2.971) -- (0.205,2.889) -- (0.123,2.801) -- (0.055,2.708) -- (0.001,2.611) -- (-0.039,2.511) -- (-0.065,2.408) -- (-0.078,2.305) -- (-0.077,2.202) -- (-0.063,2.099) -- (-0.038,1.998) -- (-0.003,1.900) -- (0.041,1.804) -- (0.094,1.712) -- (0.152,1.623) -- (0.216,1.538) -- (0.283,1.456) -- (0.353,1.377) -- (0.425,1.301) -- (0.498,1.227) -- (0.572,1.156) -- (0.646,1.086) -- (0.720,1.018) -- (0.795,0.950) -- (0.870,0.883) -- (0.946,0.816) -- (1.024,0.750) -- (1.105,0.684) -- (1.188,0.619) -- (1.275,0.554) -- (1.367,0.491) -- (1.463,0.430) -- (1.565,0.371) -- (1.672,0.316) -- (1.785,0.265) -- (1.904,0.219) -- (2.028,0.178) -- (2.156,0.144) -- (2.288,0.116) -- (2.424,0.096) -- (2.561,0.084) -- (2.700,0.080) -- (2.839,0.084) -- (2.976,0.096) -- (3.112,0.116) -- (3.244,0.144) -- (3.372,0.178) -- (3.496,0.219) -- (3.615,0.265) -- (3.728,0.316) -- (3.835,0.371) -- (3.937,0.430) -- (4.033,0.491) -- (4.125,0.554) -- (4.212,0.619) -- (4.295,0.684) -- (4.376,0.750) -- (4.454,0.816) -- (4.530,0.883) -- (4.605,0.950) -- (4.680,1.018) -- (4.754,1.086) -- (4.828,1.156) -- (4.902,1.227) -- (4.975,1.301) -- (5.047,1.377) -- (5.117,1.456) -- (5.184,1.538) -- (5.248,1.623) -- (5.306,1.712) -- (5.359,1.804) -- cycle;
   \draw[gray!45,thin] (0.05,2.10) rectangle (0.55,2.60);\fill[white] (0.30,2.35) circle (0.13);\draw[blue!70!black] (0.30,2.35) circle (0.13);
   \draw[gray!45,thin] (0.55,1.60) rectangle (1.05,2.10);\fill[white] (0.80,1.85) circle (0.13);\draw[blue!70!black] (0.80,1.85) circle (0.13);
   \draw[gray!45,thin] (0.55,2.10) rectangle (1.05,2.60);\fill[white] (0.80,2.35) circle (0.13);\draw[blue!70!black] (0.80,2.35) circle (0.13);
   \draw[gray!45,thin] (0.55,2.60) rectangle (1.05,3.10);\fill[white] (0.80,2.85) circle (0.13);\draw[blue!70!black] (0.80,2.85) circle (0.13);
   \draw[gray!45,thin] (1.05,1.10) rectangle (1.55,1.60);\fill[white] (1.30,1.35) circle (0.13);\draw[blue!70!black] (1.30,1.35) circle (0.13);
   \draw[gray!45,thin] (1.05,1.60) rectangle (1.55,2.10);\fill[white] (1.30,1.85) circle (0.13);\draw[blue!70!black] (1.30,1.85) circle (0.13);
   \draw[gray!45,thin] (1.05,2.10) rectangle (1.55,2.60);\fill[white] (1.30,2.35) circle (0.13);\draw[blue!70!black] (1.30,2.35) circle (0.13);
   \draw[gray!45,thin] (1.05,2.60) rectangle (1.55,3.10);\fill[white] (1.30,2.85) circle (0.13);\draw[blue!70!black] (1.30,2.85) circle (0.13);
   \draw[gray!45,thin] (1.55,0.60) rectangle (2.05,1.10);\fill[white] (1.80,0.85) circle (0.13);\draw[blue!70!black] (1.80,0.85) circle (0.13);
   \draw[gray!45,thin] (1.55,1.10) rectangle (2.05,1.60);\fill[white] (1.80,1.35) circle (0.13);\draw[blue!70!black] (1.80,1.35) circle (0.13);
   \draw[gray!45,thin] (1.55,1.60) rectangle (2.05,2.10);\fill[white] (1.80,1.85) circle (0.13);\draw[blue!70!black] (1.80,1.85) circle (0.13);
   \draw[gray!45,thin] (1.55,2.10) rectangle (2.05,2.60);\fill[white] (1.80,2.35) circle (0.13);\draw[blue!70!black] (1.80,2.35) circle (0.13);
   \draw[gray!45,thin] (1.55,2.60) rectangle (2.05,3.10);\fill[white] (1.80,2.85) circle (0.13);\draw[blue!70!black] (1.80,2.85) circle (0.13);
   \draw[gray!45,thin] (2.05,0.60) rectangle (2.55,1.10);\fill[white] (2.30,0.85) circle (0.13);\draw[blue!70!black] (2.30,0.85) circle (0.13);
   \draw[gray!45,thin] (2.05,1.10) rectangle (2.55,1.60);\fill[white] (2.30,1.35) circle (0.13);\draw[blue!70!black] (2.30,1.35) circle (0.13);
   \draw[gray!45,thin] (2.05,1.60) rectangle (2.55,2.10);\fill[white] (2.30,1.85) circle (0.13);\draw[blue!70!black] (2.30,1.85) circle (0.13);
   \draw[gray!45,thin] (2.05,2.10) rectangle (2.55,2.60);\fill[white] (2.30,2.35) circle (0.13);\draw[blue!70!black] (2.30,2.35) circle (0.13);
   \draw[gray!45,thin] (2.05,2.60) rectangle (2.55,3.10);\fill[white] (2.30,2.85) circle (0.13);\draw[blue!70!black] (2.30,2.85) circle (0.13);
   \draw[gray!45,thin] (2.55,0.60) rectangle (3.05,1.10);\fill[white] (2.80,0.85) circle (0.13);\draw[blue!70!black] (2.80,0.85) circle (0.13);
   \draw[gray!45,thin] (2.55,1.10) rectangle (3.05,1.60);\fill[white] (2.80,1.35) circle (0.13);\draw[blue!70!black] (2.80,1.35) circle (0.13);
   \draw[gray!45,thin] (2.55,1.60) rectangle (3.05,2.10);\fill[white] (2.80,1.85) circle (0.13);\draw[blue!70!black] (2.80,1.85) circle (0.13);
   \draw[gray!45,thin] (2.55,2.10) rectangle (3.05,2.60);\fill[white] (2.80,2.35) circle (0.13);\draw[blue!70!black] (2.80,2.35) circle (0.13);
   \draw[gray!45,thin] (2.55,2.60) rectangle (3.05,3.10);\fill[white] (2.80,2.85) circle (0.13);\draw[blue!70!black] (2.80,2.85) circle (0.13);
   \draw[gray!45,thin] (3.05,0.60) rectangle (3.55,1.10);\fill[white] (3.30,0.85) circle (0.13);\draw[blue!70!black] (3.30,0.85) circle (0.13);
   \draw[gray!45,thin] (3.05,1.10) rectangle (3.55,1.60);\fill[white] (3.30,1.35) circle (0.13);\draw[blue!70!black] (3.30,1.35) circle (0.13);
   \draw[gray!45,thin] (3.05,1.60) rectangle (3.55,2.10);\fill[white] (3.30,1.85) circle (0.13);\draw[blue!70!black] (3.30,1.85) circle (0.13);
   \draw[gray!45,thin] (3.05,2.10) rectangle (3.55,2.60);\fill[white] (3.30,2.35) circle (0.13);\draw[blue!70!black] (3.30,2.35) circle (0.13);
   \draw[gray!45,thin] (3.05,2.60) rectangle (3.55,3.10);\fill[white] (3.30,2.85) circle (0.13);\draw[blue!70!black] (3.30,2.85) circle (0.13);
   \draw[gray!45,thin] (3.55,0.60) rectangle (4.05,1.10);\fill[white] (3.80,0.85) circle (0.13);\draw[blue!70!black] (3.80,0.85) circle (0.13);
   \draw[gray!45,thin] (3.55,1.10) rectangle (4.05,1.60);\fill[white] (3.80,1.35) circle (0.13);\draw[blue!70!black] (3.80,1.35) circle (0.13);
   \draw[gray!45,thin] (3.55,1.60) rectangle (4.05,2.10);\fill[white] (3.80,1.85) circle (0.13);\draw[blue!70!black] (3.80,1.85) circle (0.13);
   \draw[gray!45,thin] (3.55,2.10) rectangle (4.05,2.60);\fill[white] (3.80,2.35) circle (0.13);\draw[blue!70!black] (3.80,2.35) circle (0.13);
   \draw[gray!45,thin] (3.55,2.60) rectangle (4.05,3.10);\fill[white] (3.80,2.85) circle (0.13);\draw[blue!70!black] (3.80,2.85) circle (0.13);
   \draw[gray!45,thin] (4.05,1.10) rectangle (4.55,1.60);\fill[white] (4.30,1.35) circle (0.13);\draw[blue!70!black] (4.30,1.35) circle (0.13);
   \draw[gray!45,thin] (4.05,1.60) rectangle (4.55,2.10);\fill[white] (4.30,1.85) circle (0.13);\draw[blue!70!black] (4.30,1.85) circle (0.13);
   \draw[gray!45,thin] (4.05,2.10) rectangle (4.55,2.60);\fill[white] (4.30,2.35) circle (0.13);\draw[blue!70!black] (4.30,2.35) circle (0.13);
   \draw[gray!45,thin] (4.05,2.60) rectangle (4.55,3.10);\fill[white] (4.30,2.85) circle (0.13);\draw[blue!70!black] (4.30,2.85) circle (0.13);
   \draw[gray!45,thin] (4.55,1.60) rectangle (5.05,2.10);\fill[white] (4.80,1.85) circle (0.13);\draw[blue!70!black] (4.80,1.85) circle (0.13);
   \draw[gray!45,thin] (4.55,2.10) rectangle (5.05,2.60);\fill[white] (4.80,2.35) circle (0.13);\draw[blue!70!black] (4.80,2.35) circle (0.13);
 \draw[thick] (5.403,1.900) -- (5.438,1.998) -- (5.463,2.099) -- (5.477,2.202) -- (5.478,2.305) -- (5.465,2.408) -- (5.439,2.511) -- (5.399,2.611) -- (5.345,2.708) -- (5.277,2.801) -- (5.195,2.889) -- (5.102,2.971) -- (4.997,3.046) -- (4.883,3.113) -- (4.760,3.173) -- (4.631,3.225) -- (4.496,3.269) -- (4.358,3.305) -- (4.217,3.333) -- (4.077,3.355) -- (3.937,3.370) -- (3.799,3.380) -- (3.663,3.385) -- (3.532,3.387) -- (3.403,3.386) -- (3.279,3.383) -- (3.158,3.379) -- (3.041,3.376) -- (2.926,3.373) -- (2.812,3.371) -- (2.700,3.370) -- (2.588,3.371) -- (2.474,3.373) -- (2.359,3.376) -- (2.242,3.379) -- (2.121,3.383) -- (1.997,3.386) -- (1.868,3.387) -- (1.737,3.385) -- (1.601,3.380) -- (1.463,3.370) -- (1.323,3.355) -- (1.183,3.333) -- (1.042,3.305) -- (0.904,3.269) -- (0.769,3.225) -- (0.640,3.173) -- (0.517,3.113) -- (0.403,3.046) -- (0.298,2.971) -- (0.205,2.889) -- (0.123,2.801) -- (0.055,2.708) -- (0.001,2.611) -- (-0.039,2.511) -- (-0.065,2.408) -- (-0.078,2.305) -- (-0.077,2.202) -- (-0.063,2.099) -- (-0.038,1.998) -- (-0.003,1.900) -- (0.041,1.804) -- (0.094,1.712) -- (0.152,1.623) -- (0.216,1.538) -- (0.283,1.456) -- (0.353,1.377) -- (0.425,1.301) -- (0.498,1.227) -- (0.572,1.156) -- (0.646,1.086) -- (0.720,1.018) -- (0.795,0.950) -- (0.870,0.883) -- (0.946,0.816) -- (1.024,0.750) -- (1.105,0.684) -- (1.188,0.619) -- (1.275,0.554) -- (1.367,0.491) -- (1.463,0.430) -- (1.565,0.371) -- (1.672,0.316) -- (1.785,0.265) -- (1.904,0.219) -- (2.028,0.178) -- (2.156,0.144) -- (2.288,0.116) -- (2.424,0.096) -- (2.561,0.084) -- (2.700,0.080) -- (2.839,0.084) -- (2.976,0.096) -- (3.112,0.116) -- (3.244,0.144) -- (3.372,0.178) -- (3.496,0.219) -- (3.615,0.265) -- (3.728,0.316) -- (3.835,0.371) -- (3.937,0.430) -- (4.033,0.491) -- (4.125,0.554) -- (4.212,0.619) -- (4.295,0.684) -- (4.376,0.750) -- (4.454,0.816) -- (4.530,0.883) -- (4.605,0.950) -- (4.680,1.018) -- (4.754,1.086) -- (4.828,1.156) -- (4.902,1.227) -- (4.975,1.301) -- (5.047,1.377) -- (5.117,1.456) -- (5.184,1.538) -- (5.248,1.623) -- (5.306,1.712) -- (5.359,1.804) -- cycle;
 \draw[very thick,red!75!black] (4.240,0.640) -- (4.295,0.684) -- (4.349,0.728) -- (4.402,0.772) -- (4.454,0.816) -- (4.505,0.861) -- (4.555,0.905) -- (4.605,0.950) -- (4.655,0.995) -- (4.704,1.040) -- (4.754,1.086) -- (4.803,1.133) -- (4.852,1.180) -- (4.902,1.227) -- (4.950,1.276) -- (4.999,1.326) -- (5.047,1.377) -- (5.094,1.429) -- (5.140,1.483) -- (5.184,1.538) -- (5.227,1.594) -- (5.268,1.652) -- (5.306,1.712) -- (5.342,1.773) -- (5.374,1.836) -- (5.403,1.900) -- (5.428,1.965) -- (5.448,2.032) -- (5.463,2.099) -- (5.474,2.168) -- (5.478,2.236) -- (5.478,2.305) -- (5.471,2.374) -- (5.458,2.443) -- (5.439,2.511) -- (5.414,2.578) -- (5.382,2.644) -- (5.345,2.708) -- (5.301,2.771) -- (5.251,2.831) -- (5.195,2.889) -- (5.134,2.944) -- (5.068,2.996) -- (4.997,3.046) -- (4.922,3.091) -- (4.843,3.134) -- (4.760,3.173) -- (4.674,3.209) -- (4.586,3.240) -- (4.496,3.269) -- (4.404,3.294);
 \node[right,red!75!black] at (5.45,1.90) {$\Gamma_1$};
 \node[left] at (-0.05,1.90) {$\Gamma_2$};
 \node[above] at (2.70,3.39) {$\Omeps$};
 \node[below,blue!70!black] at (2.70,0.03) {\small holes $\eps(k+T_0)$, $k\in I_\eps$};
 \begin{scope}[xshift=7.4cm,yshift=0.45cm]
  \fill[gray!12] (0,0) rectangle (3,3);
  \draw (0,0) rectangle (3,3);
  \fill[white] (1.5,1.5) circle (0.75);
  \draw[blue!70!black,thick] (1.5,1.5) circle (0.75);
  \node at (1.5,1.5) {$T_0$};
  \draw[->,very thick,red!75!black] (2.30,1.50) -- (1.8,1.50);
  \node[above right,red!75!black] at (2.43,1.60) {$n$};
  \node at (0.45,2.6) {$\Ystar$};
  \node[blue!70!black] at (2.45,0.75) {$S$};
  \node[below] at (1.5,0) {$Y=\Ystar\cup\overline{T_0}$};
 \end{scope}
\end{tikzpicture}
\caption{The perforated domain $\Omeps$, with the Signorini part $\Gamma_1$
and the Dirichlet part $\Gamma_2$ of the exterior boundary and an unperforated
boundary strip of width $O(\eps)$, and the reference cell; $n$ is the outward unit normal to $\Ystar$ on $S$, pointing into $T_0$.}
\label{fig:geometry}
\end{figure}

For $0<\eps<1$ define
\begin{align}
 I_\eps&=\{k\in\mathbb Z^N:\eps(k+\overline Y)\subset\Om\},
 &T^\eps&=\bigcup_{k\in I_\eps}\eps(k+T_0),\nonumber\\
 \Omeps&=\Om\setminus\overline{T^\eps},
 &S_\eps&=\partial T^\eps,
 &G_\eps&=\bigcup_{k\in I_\eps}\eps(k+\Ystar).
\end{align}
Thus only complete cells are perforated, $\Omeps$ is connected, and a
boundary strip of width $O(\eps)$ is left untouched
(Figure~\ref{fig:geometry}).  In particular,
\begin{equation}
 |\Omeps\setminus G_\eps|+
 |\Om\setminus\textstyle\bigcup_{k\in I_\eps}\eps(k+Y)|
 \leq C\eps.                                                   \label{eq:strip}
\end{equation}

The outward unit normal to $\Omeps$ is $\nu$ on $\partial\Om$ and is
denoted by $\nu_\eps$ on $S_\eps$; the latter is the rescaled normal $n$. We write $Q=\Om\times(0,T)$ and $\Qeps=\Omeps\times(0,T)$ and for periodic coefficients $b$, $A_0$ and
$A=(A_1,\ldots,A_N)$, whose precise assumptions are collected in
Section~\ref{sec:hypotheses}, we set
\[
 b^\eps(x)=b(x/\eps),\qquad
 A_0^\eps(x,t,r)=A_0(x/\eps,t/\eps^2,r),
\]
and
\[
 A^\eps(x,t,r,\xi)
 =A(x/\eps,t/\eps^2,r,\xi).
\]
The unknown $u_\eps$ is governed by
\begin{equation}
\left\{
\begin{aligned}
 \partial_t(b^\eps u_\eps)-\Div A^\eps(x,t,u_\eps,Du_\eps)
     +A_0^\eps(x,t,u_\eps)&=f &&\text{in }\Qeps,\\
 \nu\cdot A^\eps(x,t,u_\eps,Du_\eps)&\leq0,\quad
 u_\eps\leq\varphi
     &&\text{on }\Gamma _1\times(0,T),\\
 (u_\eps-\varphi)\nu\cdot A^\eps(x,t,u_\eps,Du_\eps)&=0
     &&\text{on }\Gamma _1\times(0,T),\\
 u_\eps&=0&&\text{on }\Gamma _2\times(0,T),\\
 \nu_\eps\cdot A^\eps(x,t,u_\eps,Du_\eps)
     &=g^\eps
     &&\text{on }S_\eps\times(0,T),\\
 b^\eps u_\eps(0)&=0&&\text{in }\Omeps.
\end{aligned}
\right.                                                     \label{eq:strong-micro}
\end{equation}
The Signorini relations express the unilateral character of the exterior
boundary interaction.  Here $g^\eps$ is only a shorthand for the weak
cellwise surface functional defined in Section~\ref{sec:micro}; no pointwise
representative of $g$ on $S$ is assumed.

\subsection{Related work}

We group the literature according to the four difficulties combined in
\eqref{eq:strong-micro}.

For elliptic--parabolic equations, uniqueness is delicate even without
constraints \cite{Carrillo1999,Otto1996}.  Abstract existence for degenerate
evolution inequalities and inclusions $(Bu)'+Au\ni f$ goes back to
\cite{DiBenedettoShowalter1982,Kuttler1984}, and was extended to set-valued
pseudomonotone operators in \cite{KuttlerShillor1999}.  We pose the
microscopic problem in the integrated weak form of Lions
\cite[Chapter~2]{Lions1969} (see also \cite[Section~4]{MigorskiDudek2025}),
and use the energy identities of \cite[Lemmas~2 and~3]{Kuttler1984}.

Degenerate diffusion has been studied for doubly degenerate parabolic
equations \cite{FisherLiu2002,Liu2000} and for a doubly degenerate Signorini
inequality \cite{LiuLiu2013}, a closely related microscopic model.  The
homogenization of degenerate quasilinear parabolic equations was considered in
\cite{Huang2006,NandakumaranRajesh2001,ZhangHuang2005}, and that of
pseudomonotone operators with $p$-growth, whose degeneracy is not due to
spatial weights, in \cite{Woukeng2010}.  Clark and Packer
\cite{ClarkPacker1999} obtained a two-scale limit for nonlinear degenerate
evolution equations using monotone operators; here the additional exterior
Signorini constraint, perforations with unscaled surface fluxes and weighted
diffusion require a different compactness and limit-identification argument.

For perforations we use two-scale convergence \cite{Allaire1992,Nguetseng1989}
in the space--time form of \cite{Holmbom1997}, with the diffusive scaling
$t/\eps^2$ \cite{BensoussanLionsPapanicolaou1978}.  Uniform extension operators
go back to \cite{AcerbiChiadoDalMasoPercivale1992,CioranescuSaintJeanPaulin1979}
(see also \cite{HopkerBohm2014}), and periodic unfolding offers an alternative
framework \cite{CioranescuGrisoDamlamian2012}.  Unscaled nonhomogeneous Neumann
data on the holes were treated in \cite{CioranescuDonato1988}, and multiscale
parabolic problems in perforated or fractured media, with uniformly elliptic
diffusion and positive storage, in \cite{DouanlaTetsadjio2017,DouanlaWoukeng2015}.
Lobkova \cite{Lobkova2019} treats resonant spatial and temporal oscillations
in periodically perforated domains, but for a linear parabolic equation with
homogeneous Neumann conditions on the holes, rather than the present
quasilinear weighted Signorini problem with unscaled surface forcing.

Unilateral conditions on the holes were homogenized in the elliptic case in
\cite{Pastukhova2001,Yosifian2001}, with semilinear terms and correctors in
\cite{Avila2024,Avila2025}, and with a strange term at critical hole size in
\cite{ConcaMuratTimofte2003}.  Parabolic Signorini problems were treated in
perforated domains \cite{Beliaev2004,Sandrakov2005} and in thick junctions
\cite{MelnykNakvasiuk2012}.  Dynamic unilateral conditions on critical
perforations produce nonlocal effective operators in the evolutionary regime
\cite{DiazShaposhnikovaPodolskiy2024,PodolskiyShaposhnikova2022}, and the
time-periodic regime is studied in \cite{DiazPodolskiyShaposhnikova2026}.
Two-scale evolutionary inequalities were homogenized in the rate-independent
setting \cite{MielkeTimofte2007} and, for degenerate pseudoparabolic
inequalities without microscopic uniqueness, in \cite{Ptashnyk2019}.

In the cited homogenization results for perforated Signorini problems
\cite{Avila2024,Avila2025,Beliaev2004,ConcaMuratTimofte2003,Pastukhova2001,Sandrakov2005,Yosifian2001},
the unilateral constraint typically acts on the holes.  Our constraint instead
acts on a fixed part of the exterior boundary, while the holes carry an
unscaled oscillating flux; vanishing storage and non-uniformly positive
diffusion weights are treated simultaneously.

\subsection{Standing hypotheses}\label{sec:hypotheses}

We use the following periodic-space convention throughout.
If $E$ is one of the periodicity cells $Y$, $Z$ or $Y\times Z$, then the
subscript $\#$ denotes periodicity with respect to $E$ (componentwise when
$E=Y\times Z$).  Thus, for $1\leq p\leq\infty$ and a Banach space $X$,
$L^p_\#(E;X)$ is the space of $E$-periodic, locally $L^p$ functions with
values in $X$, identified with their restrictions to $E$; we write
$L^p_\#(E)$ in the scalar case.  The same convention is used for
$C^\infty_\#$, $C_\#$, $H^1_\#$ and $\mathcal D'_\#$.  When the periodicity
cell is clear from the displayed variables, it is not repeated in the
notation.

\medskip
\noindent
All periodic functions are identified with their periodic extensions.  The
diffusion weights satisfy, for
every $i\in\{1,\ldots,N\}$,
\begin{equation}
 w_i\in L^\infty_\#(Y),\qquad w_i>0\quad\text{a.e. in }Y,
 \qquad w_i^{-s}\in L^1(Y),                                  \tag{$H_1$}\label{ass:weights}
\end{equation}
for one number $s>N/2$.  We set
$\omega=(w_1,\ldots,w_N)$ and
\begin{equation}
 q=\frac{2s}{s+1}\in(1,2);\qquad\text{then}\qquad
 s>\frac N2\iff q>\frac{2N}{N+2}.                            \label{eq:q}
\end{equation}
The values of $w_i$ inside $T_0$, like those of $b$ and of the coefficients
below, are only an auxiliary periodic extension needed on the unperforated
boundary strip; only their restrictions to $\Ystar$ enter the homogenized
problem.  The condition $s>N/2$ ensures the compact embedding
$W^{1,q}\hookrightarrow\hookrightarrow L^2$.  The weights are positive a.e. but need not
be bounded away from zero, as in the inverse-integrability framework of
\cite[Hypothesis~1.1]{Bonafede1993}.

\medskip
\noindent
The mass coefficient satisfies
\begin{equation}
 b\in L^\infty_\#(Y),\qquad b\geq0,\qquad
 \overline b:=\int_{\Ystar}b(y)\dd y>0.                       \tag{$H_2$}\label{ass:mass}
\end{equation}
Thus $b$ may vanish on a set of positive measure; only its integral over the
material cell is required to be strictly positive.

\medskip
\noindent
The zero-order coefficient $A_0=A_0(y,\tau,r)$ is periodic and measurable in
$(y,\tau)$ and continuous in $r$ a.e.  For $i=1,\ldots,N$,
$A_i=A_i(y,\tau,r,\xi)$ is periodic and measurable in $(y,\tau)$ and
continuous in $(r,\xi)$ a.e.  Write $A=(A_1,\ldots,A_N)$.  There exist
constants $c_0,c_1>0$ and functions
$k_1\in L^2_\#(Y\times Z)$ and
$k_2\in L^1_\#(Y\times Z)$, with $k_1,k_2\geq0$ a.e., such that, for a.e.
$(y,\tau)\in Y\times Z$, all $r\in\R$ and all $\xi,\zeta\in\R^N$,
\begin{align}
 |A_0(y,\tau,r)|
 &\leq c_1\left(k_1+|r|\right),                              \tag{$H_3$}\label{ass:A1a}\\
 |A_i(y,\tau,r,\xi)|
 &\leq c_1\sqrt{w_i}\left(k_1+|r|+
                   \sum_{j=1}^N\sqrt{w_j}|\xi_j|\right),
                                                               \tag{$H_4$}\label{ass:A1b}\\
 A_0(y,\tau,r)r+A(y,\tau,r,\xi)\cdot\xi
 &\geq c_0\sum_{i=1}^Nw_i|\xi_i|^2-c_1k_1|r|-k_2,
                                                               \tag{$H_5$}\label{ass:A2}\\
 [A(y,\tau,r,\xi)-A(y,\tau,r,\zeta)]
       \cdot(\xi-\zeta)&\geq0.                              \tag{$H_6$}\label{ass:A3}
\end{align}
Here \eqref{ass:A1a} and \eqref{ass:A1b} are the growth conditions,
\eqref{ass:A2} is coercivity, and \eqref{ass:A3} is monotonicity in the
gradient variable.  Coercivity is required in the weighted gradient only: since
$|\Gamma_2|_{N-1}>0$, the uniform Poincar\'e inequality of
Lemma~\ref{lem:uniform-poincare} controls the state.  In particular, if
$A(y,\tau,r,\xi)\cdot\xi\geq c_0\sum_{i=1}^Nw_i|\xi_i|^2-k_2$, then \eqref{ass:A2}
holds for $A_0\equiv0$ and, more generally, for every $A_0$ nondecreasing in
$r$ and satisfying \eqref{ass:A1a}, because then
$A_0(y,\tau,r)r\geq A_0(y,\tau,0)r\geq-c_1k_1|r|$.

\medskip
\noindent
The source and the obstacle satisfy
\begin{equation}
\begin{gathered}
 f\in L^2(Q),\quad \Phi\in H^1(0,T;H^1(\Om)),\quad
 \gamma_0\Phi=0\text{ on }\Gamma_2,\\
 \varphi=\gamma_0\Phi|_{\Gamma_1},\quad
 \varphi(0)\geq0\quad\text{a.e. on }\Gamma_1.
\end{gathered}                                                  \tag{$H_7$}\label{ass:data}
\end{equation}
The microscopic obstacle is the restriction of this single lifting.

\medskip
\noindent
The oscillating surface datum satisfies the compatibility condition
\begin{equation}
 g\in L^2_\#\bigl(Z;(W^{1-1/q,q}(S))'\bigr),\qquad
 \langle g(\tau),1\rangle_S=0\quad\text{a.e. }\tau.
                                                               \tag{$H_8$}\label{ass:surface}
\end{equation}
The zero-mean condition is the compatibility condition for the unscaled
surface flux; without it the surface contribution is of order $\eps^{-1}$.
Its first moment is
\begin{equation}
 (m_g)_i=\int_Z\langle g(\tau),y_i|_S\rangle_S\dd\tau,
 \qquad i=1,\ldots,N.
\end{equation}

\medskip
\noindent
For the microscopic existence result we additionally assume the joint
state--gradient monotonicity:  for a.e.
$(y,\tau)\in Y\times Z$ and all $r,\bar r\in\R$, $\xi,\zeta\in\R^N$,
\begin{equation}
 [A_0(y,\tau,r)-A_0(y,\tau,\bar r)](r-\bar r)+[A(y,\tau,r,\xi)-A(y,\tau,\bar r,\zeta)]\cdot(\xi-\zeta)\geq0.
 \tag{$H_9$}\label{ass:A4}
\end{equation}
It implies \eqref{ass:A3} and is used only for microscopic solvability.

\medskip
\noindent
For microscopic uniqueness, the single-valued cell law and the decoupled
macroscopic problem, we assume that there exists $c_u>0$ such
that, for a.e. $(y,\tau)$ and all $r,\bar r\in\mathbb R$,
$\eta,\zeta\in\mathbb R^N$,
\begin{align}
 &[A_0(y,\tau,r)-A_0(y,\tau,\bar r)](r-\bar r)
 +[A(y,\tau,r,\eta)-A(y,\tau,\bar r,\zeta)]\cdot(\eta-\zeta)\nonumber\\
 &\qquad\geq c_u\sum_{i=1}^Nw_i(y)|\eta_i-\zeta_i|^2.
 \tag{$H_{10}$}\label{ass:Adu}
\end{align}
Thus \eqref{ass:Adu} $\Rightarrow$ \eqref{ass:A4} $\Rightarrow$
\eqref{ass:A3}, and for $r=\bar r$ \eqref{ass:Adu} is a uniform strict
monotonicity in the weighted gradient.

\subsection{Main results and new ingredients}

Under the geometric assumptions of Section~\ref{sec:geometry} and the standing
hypotheses above, coercivity is required in the weighted gradient only.  Microscopic solutions are understood in the integrated weak form of
Lions; they exist under the joint monotonicity condition \((H_9)\)
(Theorem~\ref{thm:existence-eps}), are unique under \((H_{10})\) although $b$
may vanish (Proposition~\ref{prop:micro-uniqueness}), and need not be unique
under \((H_9)\) alone (Remark~\ref{rem:non-uniqueness}).

\begin{maintheorem}\label{thm:A}
Assume \((H_1)\)--\((H_8)\), and let $(u_\eps)$ be any family of microscopic
solutions.  It satisfies estimates uniform in $\eps$
(Theorem~\ref{thm:uniform}).  Along a subsequence, the extended states
converge strongly in $L^2(Q)$ to some $u_0\in L^2(0,T;H^1(\Om))$, the weighted
gradients converge two-scale to $\chi_{\Ystar}\sqrt{w_i}(D_{x_i}u_0+D_{y_i}u_1)$
(Theorem~\ref{thm:weighted-compactness}), and $(u_0,u_1)$ solves a coupled
two-scale problem in Minty comparison form (Theorem~\ref{thm:main}), with
$u_1$ solving a fast-time cell equation at the frozen limit state
(Proposition~\ref{prop:frozen-cell}).
\end{maintheorem}

\begin{maintheorem}\label{thm:B}
Assume \((H_1)\)--\((H_8)\) and \((H_{10})\).  Then $u_1=\pi(u_0,D_xu_0)$ is given by a
single-valued cell law, $u_0\in H^1(0,T;H^{-1}(\Om))$ with $u_0(0)=0$ in
$H^{-1}(\Om)$, and $u_0$ is the unique function in
$L^2(0,T;H^1_{\Gamma_2}(\Om))$ with $\gamma_0u_0\leq\varphi$ on
$\Gamma_1\times(0,T)$ such that
\begin{align*}
 &\int_0^T\langle\overline b\,\partial_tv,v-u_0\rangle\dd t
 +\int_Q A_0^{\rm hom}(u_0)(v-u_0)\dd x\dd t\\
 &\quad+\int_Q\bigl(A^{\rm hom}(u_0,D_xu_0)-m_g\bigr)
 \cdot D_x(v-u_0)\dd x\dd t
 \geq\theta\int_Q f(v-u_0)\dd x\dd t
 \qquad\forall\,v\in\mathcal V_{\rm ad}.
\end{align*}
Here $\mathcal V_{\rm ad}$ is the class from
Theorem~\ref{thm:decoupled}: its elements satisfy
$v\in L^2(0,T;H^1_{\Gamma_2}(\Om))$,
$\partial_tv\in L^2(0,T;(H^1_{\Gamma_2}(\Om))')$, $v(0)=0$, and
$\gamma_0v\leq\varphi$ on $\Gamma_1\times(0,T)$.  The effective
coefficients are defined in Section~\ref{sec:decoupling}, $m_g$ is the first
surface-flux moment, $\theta=|\Ystar|$, and $\overline b=\int_{\Ystar}b$.
The whole family converges (Corollary~\ref{cor:macro}).
\end{maintheorem}

Two analytical mechanisms address the loss of standard compactness.
First, because $b$ may vanish, direct Aubin--Lions compactness is unavailable.
Its nonlinear or degenerate variants \cite{Maitre2003,Moussa2016} require time
control of the sequence itself or of a fixed nonlinear transform of it,
whereas the available estimate controls $\partial_t(\rho_\eps u_\eps)$,
with an $\eps$-dependent, non-injective multiplier $\rho_\eps$.  A detection
inequality that holds for any weakly-star convergent family of densities with
positive limit (Lemma~\ref{lem:detection}), combined with time-translation
estimates and Simon's criterion \cite{Simon1986}, gives strong compactness.  Second, the
uniform estimates control the extended states only in
$L^2(0,T;W^{1,q}(\Om))$, $q<2$; the limit nevertheless satisfies
$D_xu_0\in L^2(Q)^N$, because the two-scale limit of the weighted gradients
lies in the closed range of the affine cell map
$(\xi,z)\mapsto(\sqrt{w_i}(\xi_i+D_iz))_i$
(Theorem~\ref{thm:weighted-compactness}).

Since $\partial_tu_0$ is only known in $L^2(0,T;H^{-1}(\Om))$, the macroscopic
inequality is identified in a Lions-type integrated comparison formulation,
with comparisons vanishing at $t=0$, and an exponential
regularization preserving the fixed constraint set gives uniqueness at both
scales (Lemma~\ref{lem:abstract-uniqueness}).  The surface datum enters through
the cell problem, hence through $A^{\rm hom}$, and through the boundary term
$m_g\cdot\nu$ on $\Gamma_1$; for a linear, symmetric diffusion independent of
the fast time these combine into $d_g\cdot\nu$, with
$(d_g)_i=\int_Z\langle g,(y_i+\kappa_i)|_S\rangle_S$
(Remark~\ref{rem:strong-form}).

\subsection{Organization of the paper}

Section~\ref{sec:framework} develops the functional setting and two-scale tools.
Section~\ref{sec:micro} treats
the microscopic problem, Section~\ref{sec:surface} the surface lifting and the
uniform estimates, Section~\ref{sec:compactness} compactness and the surface
limit, and Section~\ref{sec:limit} the coupled limit problem, the cell problem,
and the macroscopic problem with uniqueness and full convergence.

\section{Functional framework and two-scale tools}\label{sec:framework}

\subsection{Weighted functional setting}\label{sec:setting}

Throughout, $H^{-1}(\Om):=(H_0^1(\Om))'$; it must not be confused with
$(H^1_{\Gamma_2}(\Om))'$, the dual of the trace space used in
Section~\ref{sec:limit}.  The weighted functional setting uses \((H_1)\)--\((H_2)\); the nonlinear
two-scale recovery results also invoke the growth assumptions
\((H_3)\)--\((H_4)\).  The exponent $q$ is defined in \eqref{eq:q}.  For a bounded open set $D\subset\Om$, such as $\Omeps$ (perforated cells
together with the unperforated boundary strip) or a union of material cells, we
denote by $H^1_\omega(D)$ the completion of smooth functions on $\overline D$
for the norm below, in which the weights $w_i^\eps$ depend on $\eps$; this
dependence is omitted in the notation $H^1_\omega$:
\begin{equation}
 \|v\|_{H^1_\omega(D)}^2=\|v\|_{L^2(D)}^2+
 \sum_{i=1}^N\|\sqrt{w_i^\eps}D_i v\|_{L^2(D)}^2,
 \qquad w_i^\eps(x)=w_i(x/\eps).                         \label{eq:energynorm}
\end{equation}
By Lemma~\ref{lem:weighted-control} below, a smooth Cauchy sequence $(v_n)$
for this norm converges in $W^{1,q}(D)$, and the $L^2$ limit of
$\sqrt{w_i^\eps}D_iv_n$ is $\sqrt{w_i^\eps}D_iv$.  Hence every element of
$H^1_\omega(D)$ has a distributional gradient in $L^q(D)^N$ and the weighted
completion is unambiguous; cf. \cite[Sections~2--3]{Zhikov1998}.

\begin{lemma}\label{lem:weighted-control}
There is a constant independent of $\eps$ such that
\begin{equation}
 \|D_i v\|_{L^q(\Omeps)}
 \leq C\|\sqrt{w_i^\eps}D_i v\|_{L^2(\Omeps)},
 \qquad i=1,\ldots,N.                                    \label{eq:weighted-holder}
\end{equation}
Consequently $H^1_\omega(\Omeps)$ is continuously embedded in
$W^{1,q}(\Omeps)$, its elements have the usual $W^{1,q}$ traces, and
$W^{1,q}(D)\hookrightarrow\hookrightarrow L^2(D)$ on each fixed bounded
Lipschitz domain $D$.
\end{lemma}

\begin{proof}
Since $q/(2-q)=s$, H\"older's inequality with exponents $2/q$ and $2/(2-q)$
gives
\[
 \|D_iv\|_{L^q(\Omeps)}
 \leq\Bigl(\int_{\Omeps}(w_i^\eps)^{-s}\dd x\Bigr)^{(2-q)/(2q)}
       \|\sqrt{w_i^\eps}D_iv\|_{L^2(\Omeps)},
\]
and, covering $\Om$ by at most $C_\Om\eps^{-N}$ cells,
$\int_{\Omeps}(w_i^\eps)^{-s}\leq C_\Om\int_Yw_i^{-s}$.  This proves
\eqref{eq:weighted-holder} and the embedding.  We denote by $\gamma_0$ the
trace on $\partial\Omeps$ and by $\gamma_0^y$ the trace
$W^{1,q}(\Ystar)\to W^{1-1/q,q}(\partial\Ystar)$, restricted to $S$ when
paired with $g$.  Finally, with the Sobolev exponent $q^*=Nq/(N-q)$, $q>2N/(N+2)$ is equivalent to $q^*>2$ when
$q<N$, and Rellich's theorem gives the compact embedding.
\end{proof}

Let
\begin{equation}
 V_\eps=\{v\in H^1_\omega(\Omeps):\gamma_0 v=0
 \text{ on }\Gamma _2\},
\end{equation}
with norm \eqref{eq:energynorm}.  It is a Hilbert space.  The exterior trace
is understood through Lemma~\ref{lem:weighted-control}.

\begin{lemma}\label{lem:extension}
There are linear operators $P_\eps:V_\eps\to W^{1,q}(\Om)$ such that
$P_\eps v=v$ on $\Omeps$, $\gamma_0 P_\eps v=0$ on $\Gamma _2$, and
\begin{equation}
 \|DP_\eps v\|_{L^q(\Om)}\leq C\|Dv\|_{L^q(\Omeps)},\qquad
 \|P_\eps v\|_{W^{1,q}(\Om)}\leq C\|v\|_{V_\eps},        \label{eq:extension}
\end{equation}
where $C$ does not depend on $\eps$.
\end{lemma}

\begin{proof}
The construction follows \cite{AcerbiChiadoDalMasoPercivale1992,CioranescuSaintJeanPaulin1979}.
Let $E_0:W^{1,q}(\Ystar)\to W^{1,q}(\mathbb R^N)$ be a Sobolev extension
operator and $Ez=z_{\Ystar}+E_0(z-z_{\Ystar})|_Y$, where $z_{\Ystar}$ is the
mean of $z$ over $\Ystar$.  By the Poincar\'e inequality on the connected cell,
$\|D_yEz\|_{L^q(Y)}\leq C\|D_yz\|_{L^q(\Ystar)}$ and
$\|Ez\|_{L^q(Y)}\leq C\|z\|_{W^{1,q}(\Ystar)}$.  For $k\in I_\eps$ define
$P_\eps v(\eps(k+y))=E[v(\eps(k+\cdot))](y)$, and $P_\eps v=v$ elsewhere.
Since $T_0\Subset Y$, $P_\eps v=v$ near every cell face and near $\partial\Om$,
so $P_\eps v\in W^{1,q}(\Om)$ and $\gamma_0P_\eps v=0$ on $\Gamma_2$.  Scaling
the cell estimates, the first of which is scale invariant, and summing over
the disjoint cells gives $\|DP_\eps v\|_{L^q(\Om)}\leq C\|Dv\|_{L^q(\Omeps)}$
and $\|P_\eps v\|_{W^{1,q}(\Om)}\leq C\|v\|_{W^{1,q}(\Omeps)}$;
\eqref{eq:extension} then follows from Lemma~\ref{lem:weighted-control}.
\end{proof}

\begin{lemma}\label{lem:uniform-poincare}
There is $C_P>0$, independent of $\eps$, such that
\begin{equation}
 \|v\|_{L^2(\Omeps)}^2\leq C_P\sum_{i=1}^N\|\sqrt{w_i^\eps}D_iv\|_{L^2(\Omeps)}^2
 \qquad\forall\,v\in V_\eps .                               \label{eq:uniform-poincare}
\end{equation}
Consequently
$\|v\|_{V_\eps}^2\leq(1+C_P)\sum_{i=1}^N\|\sqrt{w_i^\eps}D_iv\|_{L^2(\Omeps)}^2$.
\end{lemma}

\begin{proof}
Since $\Om$ is connected and $|\Gamma_2|_{N-1}>0$, the Poincar\'e inequality
$\|z\|_{W^{1,q}(\Om)}\leq C\|Dz\|_{L^q(\Om)}$ holds for every
$z\in W^{1,q}(\Om)$ with $\gamma_0z=0$ on $\Gamma_2$, and
$W^{1,q}(\Om)\hookrightarrow L^2(\Om)$ because $q>2N/(N+2)$.  Applying this to
$z=P_\eps v$ and using \eqref{eq:extension} and \eqref{eq:weighted-holder},
\[
 \|v\|_{L^2(\Omeps)}\leq\|P_\eps v\|_{L^2(\Om)}
 \leq C\|DP_\eps v\|_{L^q(\Om)}\leq C\|Dv\|_{L^q(\Omeps)}
 \leq C\Bigl(\sum_{i=1}^N\|\sqrt{w_i^\eps}D_iv\|_{L^2(\Omeps)}^2\Bigr)^{1/2}.
\]
\end{proof}

On the reference material cell define
\[
 \|v\|_{\omega,\Ystar}^2=\|v\|_{L^2(\Ystar)}^2
 +\sum_{i=1}^N\int_{\Ystar}w_i(y)|D_i v(y)|^2\dd y.
\]
Let $X_\omega(\Ystar)$ be the completion of
$C^\infty(\overline{\Ystar})$ for this norm, and let
$H^1_{\#,\omega}(\Ystar)$ be the completion, for the same norm, of the
restrictions to $\Ystar$ of functions in $C^\infty_\#(Y)$.  Constants are
factored out in $X_\omega(\Ystar)/\R$.  Set
\begin{equation}
 V_\omega=\left\{v\in H^1_{\#,\omega}(\Ystar):
       \int_{\Ystar}b(y)v(y)\dd y=0\right\}.
\end{equation}

\begin{lemma}\label{lem:cell-poincare}
The weighted gradient seminorm is a norm on the preceding quotient spaces.
Moreover, for some $C>0$,
\begin{align}
 \|v\|_{L^2(\Ystar)}^2
 &\leq C\sum_{i=1}^N\int_{\Ystar}w_i|D_iv|^2
 &&\forall\,v\in V_\omega,                                \label{eq:cell-poincare}\\
 |\xi|^2+\|v\|_{L^2(\Ystar)}^2+
 \sum_{i=1}^N\int_{\Ystar}w_i|D_iv|^2
 &\leq C\sum_{i=1}^N\int_{\Ystar}w_i|\xi_i+D_iv|^2
 &&\forall\,\xi\in\R^N,\ \forall\,v\in V_\omega.        \label{eq:cell-affine}
\end{align}
\end{lemma}

\begin{proof}
Both inequalities follow by the standard contradiction argument: by the cell
version of \eqref{eq:weighted-holder}, a normalized sequence violating them is
bounded in $W^{1,q}(\Ystar)\hookrightarrow\hookrightarrow L^2(\Ystar)$, and its
limit has constant gradient on the connected set $\Ystar$.  For
\eqref{eq:cell-poincare} the gradient vanishes, so the limit is a constant
with zero $b$-mean, hence zero; the same argument with zero ordinary mean
applies to $X_\omega(\Ystar)/\mathbb R$.  For \eqref{eq:cell-affine} the limit
is $v=c-\xi\cdot y$; its traces on the opposite faces $\{y_i=0\}$ and
$\{y_i=1\}$ coincide by periodicity and differ by $-\xi_i$, so $\xi=0$ and
then $v=0$.  Since
$\sum_{i=1}^N\int_{\Ystar}w_i|D_iv_n|^2\leq2\sum_{i=1}^N\int_{\Ystar}w_i|\xi_i^n+D_iv_n|^2
+2|\xi^n|^2\max_{1\leq i\leq N}\|w_i\|_{L^1(Y)}$, this contradicts the normalization.
\end{proof}

Let $H_b$ be the completion of $V_\omega/N_b$, where
$N_b=\{v\in V_\omega:\int_{\Ystar}b|v|^2=0\}$, for the norm
$\|[v]\|_{H_b}^2=\int_{\Ystar}b|v|^2$.  Denote by
$J:V_\omega\to H_b$ the natural continuous dense map and set $B=J^*J$.
Thus $\langle Bv,\zeta\rangle=\int_{\Ystar}bv\zeta$.  The map $J$ need not
be injective; retaining its kernel is essential when $b$ vanishes.
Set
\begin{equation}
 \mathcal W_\omega=\{z\in L^2_\#(Z;V_\omega):
      \partial_\tau(bz)\in L^2_\#(Z;V_\omega')\}.
\end{equation}

\begin{lemma}\label{lem:time-space}
For every $z\in\mathcal W_\omega$, the class $Jz$ has a representative in
$C_\#([0,1];H_b)$ and
\begin{equation}
 \int_Z\langle\partial_\tau(bz),\zeta\rangle\dd\tau
 =-\int_Z\langle\partial_\tau(b\zeta),z\rangle\dd\tau,
 \qquad
 \int_Z\langle\partial_\tau(bz),z\rangle\dd\tau=0
 \qquad\forall\,z,\zeta\in\mathcal W_\omega.     \label{eq:fast-ibp}
\end{equation}
\end{lemma}

\begin{proof}
Since $b$ is independent of $\tau$, periodic mollification in $\tau$
approximates $z$ by smooth $z_n$ in the graph norm of $\mathcal W_\omega$.
For smooth functions,
$\frac{\dd}{\dd\tau}(Jz_n,J\zeta_n)_{H_b}
=\langle\partial_\tau(bz_n),\zeta_n\rangle+\langle\partial_\tau(b\zeta_n),z_n\rangle$,
which yields the uniform bound
$\sup_{\tau\in[0,1]}\|Jz_n\|_{H_b}^2\leq C\|z_n\|_{\mathcal W_\omega}^2$ and, after
integration over $Z$, \eqref{eq:fast-ibp}.  Passing to the limit gives the
claims; this is the classical integration-by-parts formula in the degenerate
setting $B=J^*J$ of \cite[Chapter~III, Section~6]{Showalter1997}.
\end{proof}

\begin{lemma}\label{lem:periodic-degenerate-derivative}
Let
\[
 \mathcal H_\omega=L^2_\#(Z;V_\omega),
 \qquad
 \mathcal Lz=\partial_\tau(bz),
 \qquad D(\mathcal L)=\mathcal W_\omega.
\]
Then $\mathcal L:D(\mathcal L)\subset\mathcal H_\omega\to
\mathcal H_\omega'$ is maximal monotone.  More precisely, if
$R_\omega:V_\omega\to V_\omega'$ denotes the Riesz isomorphism and
$\mathbf R_\omega$ its pointwise realization on $\mathcal H_\omega$, then
\begin{equation}
 \operatorname{Ran}(\mathcal L+\mathbf R_\omega)
 =\mathcal H_\omega'.                                      \label{eq:periodic-derivative-range}
\end{equation}
\end{lemma}

\begin{proof}
Monotonicity follows immediately from \eqref{eq:fast-ibp}:
\[
 \langle\mathcal Lz-\mathcal L\zeta,z-\zeta\rangle=0
 \qquad\forall\,z,\zeta\in\mathcal W_\omega.
\]
We prove \eqref{eq:periodic-derivative-range}.  Set
\[
 \mathcal M=R_\omega^{-1}B\in\mathcal L(V_\omega).
\]
Because $B=J^*J$, $\mathcal M$ is bounded, self-adjoint and non-negative in the Hilbert
space $V_\omega$ endowed with its Riesz inner product.  Let
$\mathfrak F\in\mathcal H_\omega'$ and put
$\mathfrak g=\mathbf R_\omega^{-1}\mathfrak F\in\mathcal H_\omega$.  Work first in the
complexification of $V_\omega$ and write the Fourier series
\[
 \mathfrak g(\tau)=\sum_{n\in\mathbb Z}\mathfrak g_ne^{2\pi\mathrm{i}n\tau}.
\]
For every $n\in\mathbb Z$, the operator
$I+2\pi\mathrm{i}n\mathcal M$ is bijective and
\begin{equation}
 \|(I+2\pi\mathrm{i}n\mathcal M)^{-1}\|_{\mathcal L(V_\omega)}\leq1. \label{eq:Fourier-resolvent}
\end{equation}
Indeed,
\[
 \Re\bigl((I+2\pi\mathrm{i}n\mathcal M)v,v\bigr)_{V_\omega}=\|v\|_{V_\omega}^2,
\]
which gives injectivity, closed range and the estimate in
\eqref{eq:Fourier-resolvent}; applying the same argument to the adjoint
$I-2\pi\mathrm{i}n\mathcal M$ shows that the range is dense, hence all of $V_\omega$.
Define
\[
 z_n=(I+2\pi\mathrm{i}n\mathcal M)^{-1}\mathfrak g_n,
 \qquad
 z(\tau)=\sum_{n\in\mathbb Z}z_ne^{2\pi\mathrm{i}n\tau}.
\]
Parseval's identity and \eqref{eq:Fourier-resolvent} give
\[
 \sum_{n\in\mathbb Z}\|z_n\|_{V_\omega}^2
 \leq\sum_{n\in\mathbb Z}\|\mathfrak g_n\|_{V_\omega}^2.
\]
Moreover
\[
 2\pi\mathrm{i}n\mathcal Mz_n=\mathfrak g_n-z_n,
\]
so
\[
 \sum_{n\in\mathbb Z}\|2\pi\mathrm{i}n\mathcal Mz_n\|_{V_\omega}^2
 \leq4\sum_{n\in\mathbb Z}\|\mathfrak g_n\|_{V_\omega}^2.
\]
Thus $z\in\mathcal H_\omega$ and
$\partial_\tau(\mathcal Mz)\in\mathcal H_\omega$.  Since $B=R_\omega \mathcal M$,
\[
 \partial_\tau(bz)=\partial_\tau(Bz)
 =R_\omega\partial_\tau(\mathcal Mz)\in\mathcal H_\omega',
\]
so $z\in\mathcal W_\omega$, and the Fourier identities give
\[
 \mathcal Lz+\mathbf R_\omega z=\mathfrak F.
\]
The real and imaginary parts decouple, so real data give a real solution.
Equation \eqref{eq:periodic-derivative-range} follows, and Minty's criterion
proves maximal monotonicity.
\end{proof}

\begin{lemma}
\label{lem:graph-core}
For every $z\in L^2(Q;\mathcal W_\omega)$ there are
\begin{equation}
 z^n(x,t,y,\tau)=\sum_{\ell=1}^{L_n}
 \alpha_{\ell n}(x,t)\zeta_{\ell n}(y,\tau),
\end{equation}
where $\alpha_{\ell n}\in C_c^\infty(Q)$ and each
$\zeta_{\ell n}$ is the restriction to $\Ystar\times Z$ of a function in
$C_\#^\infty(Y\times Z)$ satisfying
\[
 \int_{\Ystar}b(y)\zeta_{\ell n}(y,\tau)\dd y=0
 \quad\text{for every }\tau,
\]
such that
\begin{align}
 z^n&\longrightarrow z
 &&\text{in }L^2(Q\times Z;V_\omega),\nonumber\\
 \partial_\tau(bz^n)&\longrightarrow\partial_\tau(bz)
 &&\text{in }L^2(Q\times Z;V_\omega').
\end{align}
\end{lemma}

\begin{proof}
Periodic convolution in $\tau$ commutes with multiplication by $b$, which
reduces the claim to the case $z\in L^2(Q;H^1_\#(Z;V_\omega))$.  The set
$\mathcal D_\omega$ of the functions $\psi-\overline b^{-1}\int_{\Ystar}b\psi$,
$\psi\in C_\#^\infty(Y)|_{\Ystar}$, is dense in $V_\omega$ by definition of
$H^1_{\#,\omega}(\Ystar)$, so
$C_c^\infty(Q)\otimes C_\#^\infty(Z)\otimes\mathcal D_\omega$ is dense in
$L^2(Q;H^1_\#(Z;V_\omega))$; a diagonal argument concludes.
\end{proof}

\begin{remark}[Sufficient conditions for $(H_{10})$ and an example]\label{rem:suff-H10}
If $[A_0(r)-A_0(\bar r)](r-\bar r)\geq\mu_0|r-\bar r|^2$,
$[A(r,\eta)-A(r,\zeta)]\cdot(\eta-\zeta)\geq\alpha_0\sum_{i=1}^Nw_i|\eta_i-\zeta_i|^2$
and $\sum_{i=1}^Nw_i^{-1}|A_i(r,\zeta)-A_i(\bar r,\zeta)|^2\leq L_A^2|r-\bar r|^2$,
splitting the flux difference at $r$ and Young's inequality show that
\eqref{ass:Adu} holds with $c_u=\alpha_0/2$ whenever $\mu_0\geq L_A^2/(2\alpha_0)$;
in particular when $A$ does not depend on $r$ and $A_0$ is nondecreasing.

For example, let $y_0$ be an interior point of $\Ystar$ and put
$w_i(y)=|y-y_0|^\alpha$ near $y_0$, with $\alpha>0$, extended to a periodic function that is
bounded and bounded away from zero away from $y_0$.  Then $w_i^{-s}\in L^1(Y)$
if and only if $\alpha s<N$, so \eqref{ass:weights} holds for some $s>N/2$ if
and only if $0<\alpha<2$; the diffusion degenerates at every point
$\eps(k+y_0)$.  Any $b=\chi_E$ with $E\subset\Ystar$ of positive measure
satisfies \eqref{ass:mass} and vanishes on $\Ystar\setminus E$.  Let $a$ be measurable and periodic with $0<a_*\leq a\leq a^*$, let
$\lambda_0>0$ and $\beta\in\R^N\setminus\{0\}$, and put
\[
 A_0(y,\tau,r)=\mu r,\qquad
 A_i(y,\tau,r,\xi)=w_i(y)a(y,\tau)\bigl(\xi_i+\lambda_0\tanh\xi_i\bigr)
 +\sqrt{w_i(y)}\,\beta_ir,
 \qquad \mu\geq\frac{|\beta|^2}{2a_*}.
\]
Since $\tanh$ is nondecreasing and bounded, \((H_3)\)--\((H_6)\) hold, with
$c_0=a_*/2$, and the sufficient conditions above, with $\alpha_0=a_*$,
$L_A=|\beta|$ and $\mu_0=\mu$, give \((H_{10})\), hence \((H_9)\), with
$c_u=a_*/2$.  The flux is nonlinear in the
gradient and depends on the state, in addition to the two degeneracies.
\end{remark}

\subsection{Space--time two-scale convergence}

We use the space--time version \cite[Section~2]{Holmbom1997} of the
two-scale convergence of Nguetseng \cite{Nguetseng1989} and
Allaire~\cite{Allaire1992}.  For $1<p<\infty$, a bounded sequence $(z_\eps)$ in
$L^p(Q)$ converges weakly two-scale to $z\in L^p(Q\times Y\times Z)$, written
$z_\eps\twto z$, if
$\int_Qz_\eps\psi^\eps\dd x\dd t\to\int_Q\int_Y\int_Zz\psi\dd\tau\dd y\dd x\dd t$
for every $\psi\in L^{p'}(Q;C_\#(Y\times Z))$, where
$\psi^\eps(x,t)=\psi(x,t,x/\eps,t/\eps^2)$; it converges strongly two-scale,
written $z_\eps\tto z$, if moreover
$\|z_\eps\|_{L^p(Q)}\to\|z\|_{L^p(Q\times Y\times Z)}$
\cite[Definition~2.4]{Holmbom1997}.  Every bounded sequence in $L^p(Q)$ has a
weakly two-scale convergent subsequence, and the product of a weakly and a
strongly two-scale convergent sequence in $L^2(Q)$ converges to the integral of
the product of the limits \cite[Theorem~1.8]{Allaire1992},
\cite[Section~2]{Holmbom1997}.

\begin{proposition}
\label{prop:classical-gradient-decomposition}
Let $(z_\eps)$ be bounded in $L^p(0,T;W^{1,p}(\Om))$, where $1<p<\infty$.  Then,
after extraction, there are
\[
 z_0\in L^p(0,T;L^p_\#(Z;W^{1,p}(\Om))),\qquad
 z_1\in L^p(Q;L^p_\#(Z;W^{1,p}_\#(Y)/\mathbb R))
\]
such that
$z_\eps\twto z_0(x,t,\tau)$ and $D_iz_\eps\twto D_{x_i}z_0+D_{y_i}z_1$; if
moreover $z_\eps\to\bar z_0$ strongly in $L^p(Q)$, then $z_0=\bar z_0(x,t)$.
\end{proposition}

\noindent This is \cite[Proposition~1.14 and Corollary~1.15]{Allaire1992} for
$p=2$, with the fast time as a parameter \cite[Section~2]{Holmbom1997}; the
duality argument is the same for $1<p<\infty$.

\subsection{Nonlinear two-scale recovery}

For a.e. $y\in Y$ and $\zeta\in\mathbb R^N$, set
\begin{equation}
 \xi_j(y,\zeta)=\frac{\zeta_j}{\sqrt{w_j(y)}},\qquad
 j=1,\ldots,N,                                                \label{eq:normalized-xi}
\end{equation}
and define the normalized fluxes by
\begin{equation}
 \widehat A_0(y,\tau,r,\zeta)=A_0(y,\tau,r),\qquad
 \widehat A_i(y,\tau,r,\zeta)
 =\frac{A_i(y,\tau,r,\xi(y,\zeta))}{\sqrt{w_i(y)}},
 \quad i=1,\ldots,N.                                        \label{eq:normalizedA}
\end{equation}
This is well defined a.e. because $w_j>0$ a.e. for every
$j=1,\ldots,N$ by \eqref{ass:weights}.  The normalized coefficients have a
common linear-growth bound in $(r,\zeta)$.

\begin{proposition}
\label{prop:corrected31}
Use the normalized coefficients in \eqref{eq:normalized-xi}--\eqref{eq:normalizedA}
as functions of $(r,\zeta)$.  Let
$(r_\eps,\zeta_\eps)$ converge strongly two-scale in $L^2(Q)$ to
$(r,\zeta)$.  Then, for $i=1,\ldots,N$,
\begin{equation}
 \widehat A_i(x/\eps,t/\eps^2,r_\eps,\zeta_\eps)
 \tto \widehat A_i(y,\tau,r,\zeta)\quad\text{in }L^2(Q).  \label{eq:nemytskii-strong}
\end{equation}
If $r_\eps\tto r$ strongly in $L^2(Q)$, then also
\begin{equation}
 A_0(x/\eps,t/\eps^2,r_\eps)
 \tto A_0(y,\tau,r)\quad\text{in }L^2(Q).                 \label{eq:A0-strong}
\end{equation}
Consequently, the product of any strongly two-scale convergent output above
with a weakly two-scale convergent $L^2$ sequence converges to the product of
the corresponding limits.
\end{proposition}

\begin{proof}
Set $\mathcal C=Y\times Z$, $R_\eps=(r_\eps,\zeta_\eps)$, $R=(r,\zeta)$,
and recall that
$|\widehat A_i(y,\tau,z)|\leq C(k_1(y,\tau)+|z|)$ by \eqref{ass:A1b}.

We first choose $R^m\in C_c^\infty(Q;C_\#^\infty(\mathcal C))^{N+1}$
with $R^m\to R$ in $L^2(Q\times\mathcal C)$.  As in the proof of
\cite[Theorem~1.8]{Allaire1992}, strong two-scale convergence gives
\begin{equation}
 \lim_{\eps\to0}\|R_\eps-(R^m)^\eps\|_{L^2(Q)}^2
 =\|R-R^m\|_{L^2(Q\times\mathcal C)}^2,
 \label{eq:direct-recovery}
\end{equation}
which implies that $\{|R_\eps|^2\}$ is uniformly integrable; by periodicity,
so is $\{|k_1^\eps|^2\}$.

Next, for fixed $m$, $H^m=\widehat A_i(\cdot,R^m)$ belongs to
$L^2_\#(\mathcal C;C(\overline Q))$, with
$\sup_{\overline Q}|H^m|\leq C(k_1+\|R^m\|_\infty)$.  By the mean-value property
of $L^1_\#(\mathcal C;C(\overline Q))$-functions
\cite[Section~5]{Allaire1992}, \cite[Section~2]{Holmbom1997}, applied to $H^m\psi$
and to $|H^m|^2$, $(H^m)^\eps\tto H^m$.

It remains to show
$\lim_{m\to\infty}\limsup_{\eps\downarrow0}\|\widehat A_i^\eps(R_\eps)-\widehat A_i^\eps((R^m)^\eps)\|_{L^2(Q)}=0$.  On
the set where $|R_\eps|>M$ or $|(R^m)^\eps|>M$, whose measure is $O(M^{-2})$,
the growth bound and \eqref{eq:direct-recovery} make the contribution small as $M\to\infty$,
uniformly in $\eps$ and large $m$.  On its complement, the Scorza-Dragoni
theorem \cite{ScorzaDragoni1948} provides a compact $\mathcal C_\gamma\subset\mathcal C$
with $|\mathcal C\setminus \mathcal C_\gamma|<\gamma$ on which $\widehat A_i$ is uniformly
continuous on $\mathcal C_\gamma\times\overline B_M$; the set where
$(x/\eps,t/\eps^2)\notin \mathcal C_\gamma$ has measure at most $|Q|\gamma+o(1)$, and the
set where $|R_\eps-(R^m)^\eps|\geq a$ has measure at most
$a^{-2}\|R-R^m\|^2+o(1)$ by \eqref{eq:direct-recovery}.  Uniform integrability
of $|k_1^\eps|^2$ concludes.  The same argument gives
$\widehat A_i(\cdot,R^m)\to \widehat A_i(\cdot,R)$ in $L^2(Q\times\mathcal C)$.  Combining with
the convergence $(H^m)^\eps\tto H^m$ yields \eqref{eq:nemytskii-strong}; \eqref{eq:A0-strong} is identical,
and the product statement is \cite[Theorem~1.8]{Allaire1992}.
\end{proof}

\begin{corollary}
\label{cor:minty-recovery}
Let $v_0\in L^2(0,T;H^1(\Om))$ and let
\begin{equation*}
 v_1(x,t,y,\tau)=\sum_{\ell=1}^L
       \alpha_\ell(x,t)z_\ell(y,\tau),
\end{equation*}
where $\alpha_\ell\in C_c^\infty(Q)$ and each $z_\ell$ is the restriction to
$\Ystar\times Z$ of a smooth $Y\times Z$-periodic function.  Using these
periodic extensions, define on $Q$
\begin{equation*}
 v_\eps=v_0+\eps v_1(x,t,x/\eps,t/\eps^2).
\end{equation*}
Then, for every $i=1,\ldots,N$,
\begin{align}
 \sqrt{w_i^\eps}D_i v_\eps
 &\tto \sqrt{w_i(y)}(D_{x_i}v_0+D_{y_i}v_1)
       &&\text{strongly in }L^2(Q),                    \label{eq:minty-recovery-full}\\
 \chi_{\Omeps}\sqrt{w_i^\eps}D_i v_\eps
 &\tto \chi_{\Ystar}\sqrt{w_i(y)}
       (D_{x_i}v_0+D_{y_i}v_1)
       &&\text{strongly in }L^2(Q).                    \label{eq:minty-recovery}
\end{align}
Moreover $v_\eps-v_0\to0$ strongly in $L^2(Q)$, and if
$r_\eps\to r$ strongly in $L^2(Q)$, then, after extending $r_\eps$ to $Q$,
\begin{equation*}
 \chi_{\Omeps}
 \widehat A_i(x/\eps,t/\eps^2,r_\eps,
       (\sqrt{w_j^\eps}D_jv_\eps)_{j=1}^N)
 \tto
 \chi_{\Ystar}\widehat A_i(y,\tau,r,
       (\sqrt{w_j}(D_{x_j}v_0+D_{y_j}v_1))_{j=1}^N)
\end{equation*}
strongly in $L^2(Q)$, for $i=1,\ldots,N$.
\end{corollary}

\begin{proof}
Choose $v_0^m\in C^\infty(\overline Q)$ such that
$v_0^m\to v_0$ in $L^2(0,T;H^1(\Om))$.  For fixed $m$ and
$v_\eps^m=v_0^m+\eps v_1^\eps$, we have
$D_iv_\eps^m=D_{x_i}v_0^m+(D_{y_i}v_1)^\eps+\eps(D_{x_i}v_1)^\eps$, whose first
two terms are admissible and whose last tends to zero; approximating
$\sqrt{w_i}$ by $h_k\in C^\infty_\#(Y)$ in $L^2(Y)$ gives an error bounded by
$C_m\|h_k-\sqrt{w_i}\|_{L^2(Y)}$ uniformly in $\eps$, where $C_m$ depends on
$v_0^m$ and $v_1$.  Letting $\eps\to0$ and then $k\to\infty$ gives, for each
fixed $m$, the strong two-scale recovery for $v_\eps^m$.  Since $w_i\in L^\infty(Y)$,
\begin{align*}
 \|\sqrt{w_i^\eps}D_i(v_0-v_0^m)\|_{L^2(Q)}
 &\leq \|w_i\|_{L^\infty(Y)}^{1/2}
       \|D_i(v_0-v_0^m)\|_{L^2(Q)}.
\end{align*}
Thus, letting finally $m\to\infty$ proves the
full-domain recovery \eqref{eq:minty-recovery-full}.  Multiplication by the
strongly two-scale convergent characteristic function
$\chi_{\Omeps}\tto\chi_{\Ystar}$ then gives the material-domain
recovery \eqref{eq:minty-recovery}.  The term
$\eps D_xv_1^\eps$ tends to zero strongly in $L^2(Q)$ by periodic averaging,
and $\eps v_1^\eps\to0$ strongly in $L^2(Q)$, which proves the state recovery.
The nonlinear conclusion follows from Proposition~\ref{prop:corrected31}
and the strong two-scale convergence
$\chi_{\Omeps}\tto\chi_{\Ystar}$.
\end{proof}

\section{The microscopic variational problem and solvability}\label{sec:micro}

We interpret \eqref{eq:strong-micro} in the weighted spaces of
Section~\ref{sec:setting} and prove solvability for each fixed $\eps>0$.  We use
\[
 A_i^\eps(x,t,r,\xi)=A_i(x/\eps,t/\eps^2,r,\xi),\qquad
 A_0^\eps(x,t,r)=A_0(x/\eps,t/\eps^2,r),
\]
and $b^\eps(x)=b(x/\eps)$; we abbreviate
$A^\eps[u]=A^\eps(x,t,u,Du)$ and
$A_0^\eps[u]=A_0^\eps(x,t,u)$.

For $v\in V_\eps$ define the oscillating surface functional by
\begin{equation}
 \langle\ell_\eps(t),v\rangle
 =\eps^{N-1}\sum_{k\in I_\eps}
 \left\langle g(t/\eps^2),
  \gamma_0^y[v(\eps(k+\cdot))]\right\rangle_S.                \label{eq:surface-functional}
\end{equation}
For smooth $g$ this coincides with
$\int_{S_\eps}g(x/\eps,t/\eps^2)v\,\dd\sigma_x$.

The surface equality in \eqref{eq:strong-micro} is interpreted by the
cellwise duality \eqref{eq:surface-functional}.  We now define
\begin{align}
 K_\eps(t)&=\{v\in V_\eps:\gamma_0 v\leq\varphi(t)
                   \text{ on }\Gamma _1\},               \\
 \langle B_\eps u,v\rangle&=\int_{\Omeps}b^\eps uv\dd x,\nonumber\\
 \langle\mathcal A_\eps(t)u,v\rangle
 &=\int_{\Omeps}\left[A_0^\eps(x,t,u)v+
              A^\eps(x,t,u,Du)\cdot Dv\right]\dd x.
\end{align}
Put
\[
 X_\eps=L^2(0,T;V_\eps),\qquad \Psi(t)=\Phi(t)-\Phi(0),
\]
and
\[
 C_{\eps,0}=\{z\in V_\eps:\gamma_0z\leq\varphi(0)
                    \text{ on }\Gamma_1\}.
\]
The compatibility $\varphi(0)\geq0$ gives $0\in C_{\eps,0}$, and
\[
 u(t)\in K_\eps(t)\quad\Longleftrightarrow\quad
 z(t):=u(t)-\Psi(t)\in C_{\eps,0}.
\]
For $0<\sigma\leq T$ let
\[
 \mathcal K_\eps^1(\sigma)=\{v\in H^1(0,\sigma;V_\eps):
       v(t)\in K_\eps(t)\text{ for a.e. }t\in(0,\sigma)\}.
\]
The endpoint values of such comparisons are meaningful in $V_\eps$.
We use the following integrated notion of solution, in which the time
derivative is transferred to the regular comparison and the initial condition
enters through the quadratic mass term (cf. \cite[Chapter~2]{Lions1969},
\cite[Section~4]{MigorskiDudek2025}).

A function $u_\eps\in X_\eps$ is a \emph{microscopic variational solution} if
$u_\eps(t)\in K_\eps(t)$ for a.e. $t$, if
\begin{equation}
 m_\eps:=\chi_{\Omeps}b^\eps u_\eps
 \in H^1(0,T;H^{-1}(\Om)),\qquad m_\eps(0)=0,             \label{eq:initial-mass}
\end{equation}
and if
\begin{align}
 &\int_0^\sigma\langle B_\eps\partial_tv(t),v(t)-u_\eps(t)\rangle\dd t
 +\int_0^\sigma\langle\mathcal A_\eps(t)u_\eps(t),v(t)-u_\eps(t)\rangle\dd t\nonumber\\
 &\qquad
 +\frac12\langle B_\eps v(0),v(0)\rangle
 \geq\int_0^\sigma\int_{\Omeps}f(t)(v-u_\eps)\dd x\dd t
 +\int_0^\sigma\langle\ell_\eps(t),v(t)-u_\eps(t)\rangle\dd t\nonumber\\
 &\hspace{50mm}\forall\,\sigma\in(0,T],\ \forall\,v\in\mathcal K_\eps^1(\sigma).
                                                               \label{eq:microVI}
\end{align}
In addition, the unconstrained directions satisfy the distributional balance
\begin{align}
 &\int_0^T\langle m_\eps'(t),\zeta(t)\rangle\dd t
 +\int_{\Qeps}\bigl(A_0^\eps[u_\eps]\zeta
                    +A^\eps[u_\eps]\cdot D\zeta\bigr)\dd x\dd t\nonumber\\
 &\qquad=\int_{\Qeps}f\zeta\dd x\dd t
 +\int_0^T\langle\ell_\eps(t),\zeta(t)\rangle\dd t
 \quad\forall\,\zeta\in L^2(0,T;H_0^1(\Om)).          \label{eq:micro-interior-balance}
\end{align}
The term $\frac12\langle B_\eps v(0),v(0)\rangle$ encodes the zero initial
mass without requiring a trace of $u_\eps$ itself.

\begin{remark}
If $B_\eps u_\eps\in H^1(0,T;V_\eps')$, $B_\eps u_\eps(0)=0$ and the usual
pointwise-in-time inequality holds, the energy identity for $B_\eps$
\cite[Lemmas~2 and~3]{Kuttler1984} gives \eqref{eq:microVI} and
\eqref{eq:micro-interior-balance}; the integrated formulation is thus weaker
and does not require a bound on the time derivative in $V_\eps'$, which is not
available in general for a Signorini constraint.
\end{remark}

\begin{lemma}
\label{lem:positive-mass-VI}
Let $V$ be a real separable Hilbert space, $C\subset V$ non-empty, closed and
convex with $0\in C$, and let $M\in\mathcal L(V,V')$ be symmetric and
coercive:
\begin{equation}
 m_M\|v\|_V^2\leq\langle Mv,v\rangle\leq M_M\|v\|_V^2
 \qquad\forall\,v\in V
\end{equation}
for some $m_M,M_M>0$.  Let $\mathcal N:(0,T)\times V\to V'$ satisfy:
\begin{enumerate}
\item $t\mapsto \mathcal N(t,v)$ is measurable for every $v\in V$;
\item $\mathcal N(t,\cdot)$ is monotone and hemicontinuous for a.e. $t$;
\item for some $a\in L^2(0,T)$ and $c_{\mathcal N}>0$,
\[
 \|\mathcal N(t,v)\|_{V'}\leq a(t)+c_{\mathcal N}\|v\|_V;
\]
\item for some $\alpha_{\mathcal N}>0$ and $\beta_{\mathcal N}\in L^1(0,T)$,
\[
 \langle \mathcal N(t,v),v\rangle\geq
 \alpha_{\mathcal N}\|v\|_V^2-\beta_{\mathcal N}(t).
\]
\end{enumerate}
Then, for every $\mathfrak G\in L^2(0,T;V')$, there exists
\[
 z\in H^1(0,T;V),\qquad z(0)=0,\qquad z(t)\in C\quad\text{a.e.},
\]
such that
\begin{equation}
 \langle Mz'(t)+\mathcal N(t,z(t))-\mathfrak G(t),w-z(t)\rangle\geq0
 \quad\forall\,w\in C,\ \text{a.e. }t.                       \label{eq:positive-mass-VI}
\end{equation}
\end{lemma}

\begin{proof}
Endow $V$ with the scalar product $(v,w)_M=\langle Mv,w\rangle$, whose Riesz
map is $M$, and let $X_M=L^2(0,T;V)$ with the scalar product
$(u,v)_{X_M}=\int_0^T(u,v)_M\dd t$; its norm is equivalent to that of
$L^2(0,T;V)$ by the bounds on $M$.  Let $\mathcal K_M$ be the set of
$C$-valued elements of $X_M$ and
\[
 \mathcal Lz=z'+N_{\mathcal K_M}(z),\qquad
 D(\mathcal L)=\{z\in H^1(0,T;V)\cap\mathcal K_M:z(0)=0\}.
\]
We first show that $\mathcal L$ is maximal monotone.  By
\cite[Chapter~III, Theorem~3.6]{Brezis1973}, applied in
the Hilbert space $(V,(\cdot,\cdot)_M)$ to the proper, convex, lower
semicontinuous function $I_C+\frac12\|\cdot\|_M^2$, whose subdifferential is
$N_C+I$, the problem $z'+N_C(z)+z\ni h$, $z(0)=0$, has a solution in $D(\mathcal L)$ for every $h\in X_M$; together with
the monotonicity of $\mathcal L$, which follows from $(z',z)_{X_M}=\frac12\|z(T)\|_M^2$
for $z(0)=0$ and from the monotonicity of $N_{\mathcal K_M}$, this shows that $\mathcal L$ is maximal monotone in $X_M$.

Next, let $\mathfrak N z=M^{-1}\mathcal N(\cdot,z)$.  By
(1)--(3), $\mathfrak N:X_M\to X_M$ is monotone, hemicontinuous and bounded,
$\|\mathfrak Nz\|_{X_M}\leq C(\|a\|_{L^2(0,T)}+\|z\|_{X_M})$; hence it is
pseudomonotone.

It remains to check the coercivity of the sum.  For $z\in D(\mathcal L)$ and
$\xi\in\mathcal Lz$,
$\xi=z'+\eta$ with $\eta\in N_{\mathcal K_M}(z)$; since $0\in\mathcal K_M$,
$(\eta,z)_{X_M}\geq0$, and $(z',z)_{X_M}=\frac12\|z(T)\|_M^2\geq0$.  Hence, by
(4) and $\|v\|_V^2\geq M_M^{-1}\|v\|_M^2$,
\[
 (\xi+\mathfrak Nz-M^{-1}\mathfrak G,z)_{X_M}
 \geq\frac{\alpha_{\mathcal N}}{M_M}\|z\|_{X_M}^2-\|\beta_{\mathcal N}\|_{L^1(0,T)}
 -\|M^{-1}\mathfrak G\|_{X_M}\|z\|_{X_M}.
\]
Since $0\in D(\mathcal L)$, this inequality gives condition (2.1) of
\cite[Theorem~2.2]{Le2011} with $u_0=0$ and $R$ large, and that range theorem
for bounded pseudomonotone perturbations of maximal monotone operators in a
reflexive space yields $z\in D(\mathcal L)$ with
$M^{-1}\mathfrak G\in\mathcal Lz+\mathfrak Nz$.
Localizing the normal cone on a countable dense subset of $C$ yields
\eqref{eq:positive-mass-VI}.
\end{proof}

\begin{lemma}\label{lem:exp-regularization}
Let $V$ be a real separable Hilbert space, $C\subset V$ closed and convex with
$0\in C$, $0<\sigma\leq T$, and $z\in L^2(0,\sigma;V)$ with $z(t)\in C$ for a.e. $t$.
For $h>0$ put
\begin{equation}
 z_h(t)=\frac1h\int_0^t e^{-(t-r)/h}z(r)\dd r .
\end{equation}
Then $z_h\in H^1(0,\sigma;V)$, $z_h(0)=0$, $z_h(t)\in C$ for every $t$,
$hz_h'+z_h=z$, and $z_h\to z$ in $L^2(0,\sigma;V)$ as $h\to0$.  Consequently,
for every symmetric non-negative $B\in\mathcal L(V,V')$,
\begin{equation}
 \int_0^\sigma\langle Bz_h',z_h-z\rangle\dd t
 =-h\int_0^\sigma\langle Bz_h',z_h'\rangle\dd t\leq0.          \label{eq:exp-reg-B}
\end{equation}
\end{lemma}

\begin{proof}
Differentiation gives $hz_h'=z-z_h\in L^2(0,\sigma;V)$, and $z_h(0)=0$.  For
$t>0$, $z_h(t)=(1-e^{-t/h})\bar z(t)$, where $\bar z(t)$ is the mean of $z$ on
$(0,t)$ for the probability density proportional to $e^{-(t-r)/h}$; since $C$
is closed and convex, $\bar z(t)\in C$, and $z_h(t)\in C$ as a convex
combination of $\bar z(t)$ and $0$.  The convergence is the approximate-identity
property of the kernel $h^{-1}e^{-r/h}\mathbf 1_{r>0}$, and \eqref{eq:exp-reg-B}
follows from $z_h-z=-hz_h'$.
\end{proof}

This averaging preserves any fixed closed convex constraint set containing $0$;
it is the regularization device of \cite[Chapter~2]{Lions1969} and
\cite{Brezis1972} for weak solutions of evolution inequalities.  It also gives
uniqueness for integrated inequalities with a non-negative, possibly
non-injective, mass.

\begin{lemma}\label{lem:abstract-uniqueness}
Let $V$ and $C$ be as in Lemma~\ref{lem:exp-regularization}, let
$B\in\mathcal L(V,V')$ be symmetric and non-negative, let
$\Psi\in H^1(0,T;V)$ with $\Psi(0)=0$, and let
$\mathfrak G,\mathcal A^1,\mathcal A^2\in L^2(0,T;V')$.  Assume that
$u^1,u^2\in L^2(0,T;V)$ satisfy $u^j(t)-\Psi(t)\in C$ for a.e. $t$ and
\begin{equation}
 \int_0^T\langle Bv'+\mathcal A^j-\mathfrak G,v-u^j\rangle\dd t\geq0,
 \quad\forall\,v\in\mathcal D_C,\quad j=1,2,                  \label{eq:abstract-VI}
\end{equation}
where $\mathcal D_C=\{v\in H^1(0,T;V):v(0)=0,\ 
 v(t)-\Psi(t)\in C\text{ for a.e. }t\}$.  Then
\begin{equation}
 \int_0^T\langle\mathcal A^1-\mathcal A^2,u^1-u^2\rangle\dd t\leq0.
                                                           \label{eq:abstract-uniqueness}
\end{equation}
\end{lemma}

\begin{proof}
Put $z^j=u^j-\Psi$, $\bar z=\frac12(z^1+z^2)$, let $\bar z_h$ be the
regularization of $\bar z$ on $(0,T)$ given by
Lemma~\ref{lem:exp-regularization}, and $v=\Psi+\bar z_h$.  Then
$v\in H^1(0,T;V)$, $v(0)=0$ and $v-\Psi\in C$, so $v$ is admissible in
\eqref{eq:abstract-VI} for $j=1,2$.  Since $2v-u^1-u^2=2(\bar z_h-\bar z)$,
adding the two inequalities and using \eqref{eq:exp-reg-B} gives
\[
 \sum_{j=1}^2\int_0^T\langle\mathcal A^j-\mathfrak G,v-u^j\rangle\dd t
 \geq-2\int_0^T\langle B\Psi',\bar z_h-\bar z\rangle\dd t
 +2h\int_0^T\langle B\bar z_h',\bar z_h'\rangle\dd t
 \geq-2\int_0^T\langle B\Psi',\bar z_h-\bar z\rangle\dd t .
\]
As $h\to0$, $v\to\frac12(u^1+u^2)$ in $L^2(0,T;V)$ and the right-hand side
tends to zero.  Since
$\frac12(u^1+u^2)-u^1=\frac12(u^2-u^1)=-\bigl(\frac12(u^1+u^2)-u^2\bigr)$, the
terms containing $\mathfrak G$ cancel and we obtain
$-\frac12\int_0^T\langle\mathcal A^1-\mathcal A^2,u^1-u^2\rangle\dd t\geq0$.
\end{proof}

\begin{theorem}\label{thm:existence-eps}
Under \((H_1)\)--\((H_9)\), for every fixed $\eps>0$ there exists at least
one microscopic variational solution in the sense of
\eqref{eq:initial-mass}--\eqref{eq:micro-interior-balance}.
\end{theorem}

\begin{proof}
By Lemma~\ref{lem:weighted-control}, $K_\eps(t)$ and $C_{\eps,0}$ are closed
and convex, $0\in C_{\eps,0}$, and, $I_\eps$ being finite, the cell trace
theorem gives $F_\eps:=f+\ell_\eps\in L^2(0,T;V_\eps')$.  By
\eqref{ass:A1a}--\eqref{ass:A2} and \eqref{ass:A4}, $\mathcal A_\eps(t)$ is
monotone, hemicontinuous, of linear growth, and coercive up to a term in
$L^1(0,T)$; $B_\eps$ is symmetric and non-negative.

We now regularize the mass.  Let $\mathcal R_\eps$ be the Riesz
map of $V_\eps$, $B_{\eps,\delta}=B_\eps+\delta\mathcal R_\eps$ ($0<\delta\leq1$),
$\widetilde{\mathcal A}_\eps(t,z)=\mathcal A_\eps(t)(z+\Psi(t))$ and
$\mathfrak G_{\eps,\delta}=F_\eps-B_{\eps,\delta}\Psi'$.  Since $\Psi\in H^1(0,T;V_\eps)$,
\eqref{ass:A1a}--\eqref{ass:A2}, Young's inequality and
Lemma~\ref{lem:uniform-poincare} give
$\langle\widetilde{\mathcal A}_\eps(t,z),z\rangle\geq\frac{c_0}{4(1+C_P)}\|z\|_{V_\eps}^2-\mathfrak r_\eps(t)$
with $\mathfrak r_\eps\in L^1(0,T)$.  Lemma~\ref{lem:positive-mass-VI}, with
$M=B_{\eps,\delta}$, $C=C_{\eps,0}$, $\mathcal N=\widetilde{\mathcal A}_\eps$ and
$\mathfrak G=\mathfrak G_{\eps,\delta}$, gives $z_\delta\in H^1(0,T;V_\eps)$, $z_\delta(0)=0$,
$z_\delta(t)\in C_{\eps,0}$ a.e., with
\begin{equation}
 \langle B_{\eps,\delta}z_\delta'+\widetilde{\mathcal A}_\eps(t,z_\delta)
 -\mathfrak G_{\eps,\delta},w-z_\delta\rangle\geq0
 \qquad\forall\,w\in C_{\eps,0},\ \text{a.e. }t.                     \label{eq:regularized-fixed-VI}
\end{equation}
Testing with $w=0$ gives a bound on $z_\delta$ in $X_\eps$ independent of
$\delta$.  Hence, along a subsequence, $z_\delta\rightharpoonup z$ in $X_\eps$,
$\widetilde{\mathcal A}_\eps(\cdot,z_\delta)\rightharpoonup \mathfrak h$ in $X_\eps'$, and
$z(t)\in C_{\eps,0}$ a.e.

To identify $\mathfrak h$, fix $0<\sigma\leq T$; testing
\eqref{eq:regularized-fixed-VI} with $w=v(t)$, using the product rule for
$v-z_\delta$ and discarding the non-negative terminal term gives
\begin{align}
 &\int_0^\sigma\langle B_{\eps,\delta}v'+\widetilde{\mathcal A}_\eps(t,z_\delta)
 -\mathfrak G_{\eps,\delta},v-z_\delta\rangle\dd t
 +\frac12\langle B_{\eps,\delta}v(0),v(0)\rangle\geq0\nonumber\\
 &\hspace{30mm}\forall\,v\in H^1(0,\sigma;V_\eps)\text{ with }v(t)\in C_{\eps,0}\text{ on }(0,\sigma).
                                                         \label{eq:regularized-integrated}
\end{align}
Let $z_h$ be the regularization of $z$ on $(0,\sigma)$ given by
Lemma~\ref{lem:exp-regularization} with $V=V_\eps$ and $C=C_{\eps,0}$;
in particular $z_h(0)=0$.  Taking $v=z_h$ in
\eqref{eq:regularized-integrated}, using the uniform $X_\eps$ bound and
$B_{\eps,\delta}\to B_\eps$ in $\mathcal L(V_\eps,V_\eps')$, and then letting
$\delta\downarrow0$, gives
\[
 \limsup_{\delta\downarrow0}\int_0^\sigma\langle\widetilde{\mathcal A}_\eps z_\delta,
       z_\delta-z_h\rangle\dd t
 \leq\int_0^\sigma\langle B_\eps z_h'-\mathfrak G_\eps,z_h-z\rangle\dd t,
 \quad\mathfrak G_\eps=F_\eps-B_\eps\Psi'.
\]
Adding the weak limit of
$\int_0^\sigma\langle\widetilde{\mathcal A}_\eps z_\delta,z_h-z\rangle\dd t$
yields
\begin{align*}
 \limsup_{\delta\downarrow0}\int_0^\sigma
 \langle\widetilde{\mathcal A}_\eps z_\delta,z_\delta-z\rangle\dd t
 &\leq\int_0^\sigma\langle B_\eps z_h',z_h-z\rangle\dd t\\
 &\quad+\int_0^\sigma\langle\mathfrak h-\mathfrak G_\eps,z_h-z\rangle\dd t.
\end{align*}
The first integral is non-positive by \eqref{eq:exp-reg-B};
the second tends to zero as $h\downarrow0$.  Hence
\begin{equation}
 \limsup_{\delta\downarrow0}\int_0^\sigma
 \langle\widetilde{\mathcal A}_\eps z_\delta,z_\delta-z\rangle\dd t\leq0.
 \label{eq:mass-minty-defect}
\end{equation}
For every $w\in L^2(0,\sigma;V_\eps)$, monotonicity implies
$\int_0^\sigma\langle\widetilde{\mathcal A}_\eps z_\delta-
\widetilde{\mathcal A}_\eps w,z_\delta-w\rangle\dd t\geq0$.
Using \eqref{eq:mass-minty-defect} and weak convergence gives
$\int_0^\sigma\langle\mathfrak h-
\widetilde{\mathcal A}_\eps w,z-w\rangle\dd t\geq0$.
Taking $w=z\pm\lambda\zeta$ and letting $\lambda\downarrow0$ by
hemicontinuity shows $\mathfrak h=\widetilde{\mathcal A}_\eps(\cdot,z)$.
Monotonicity with $w=z$ then gives the matching lower bound zero in
\eqref{eq:mass-minty-defect}, so the defect converges to zero.
Passing to the limit in \eqref{eq:regularized-integrated} and setting
$u_\eps=z+\Psi$, $v=\widehat v-\Psi$ with $\widehat v\in\mathcal K_\eps^1(\sigma)$,
gives \eqref{eq:microVI}.

It remains to recover the interior balance and the initial mass.  With
$u_\delta=z_\delta+\Psi$, the tests $u_\delta\pm r\zeta|_{\Omeps}$ give
\begin{equation}
 \langle(B_{\eps,\delta}u_\delta)',\zeta\rangle
 +\langle\mathcal A_\eps(t)u_\delta,\zeta\rangle
 =\langle F_\eps,\zeta\rangle
 \qquad\forall\,\zeta\in H_0^1(\Om),\ \text{a.e. }t.     \label{eq:regularized-interior-balance}
\end{equation}
The growth bounds and the uniform $L^2(0,T;V_\eps)$ estimate make the
right-hand side of \eqref{eq:regularized-interior-balance} bounded in
$L^2(0,T;H^{-1}(\Om))$, independently of $\delta$.  Thus
$M_\delta=B_{\eps,\delta}u_\delta|_{H_0^1(\Om)}$ is bounded in
$H^1(0,T;H^{-1}(\Om))$, $M_\delta(0)=0$, and, for fixed $\eps$,
\[
 \|M_\delta-\chi_{\Omeps}b^\eps u_\delta\|_{L^2(0,T;H^{-1}(\Om))}
 =\delta\|\mathcal R_\eps u_\delta|_{H_0^1(\Om)}\|_{L^2(0,T;H^{-1}(\Om))}
 \leq C_\eps\delta.
\]
Therefore $M_\delta\rightharpoonup\chi_{\Omeps}b^\eps u_\eps$ in
$H^1(0,T;H^{-1}(\Om))$, which gives \eqref{eq:initial-mass}, and passing to the
limit in \eqref{eq:regularized-interior-balance} gives
\eqref{eq:micro-interior-balance}.
\end{proof}

Lemma~\ref{lem:abstract-uniqueness}, applied with the non-negative and
possibly non-injective mass $B_\eps$, gives uniqueness of the microscopic
solution.

\begin{proposition}\label{prop:micro-uniqueness}
Under \((H_1)\)--\((H_8)\) and \((H_{10})\), for every fixed $\eps>0$ the
microscopic variational solution is unique in $X_\eps$.
\end{proposition}

\begin{proof}
Let $u^1,u^2$ be two solutions.  Since $v(t)\in K_\eps(t)$ if and only if
$v(t)-\Psi(t)\in C_{\eps,0}$, and since the initial quadratic term vanishes when
$v(0)=0$, \eqref{eq:microVI} with $\sigma=T$ is \eqref{eq:abstract-VI} with
$V=V_\eps$, $C=C_{\eps,0}$, $B=B_\eps$, $\mathfrak G=F_\eps$ and
$\mathcal A^j=\mathcal A_\eps(\cdot)u^j\in L^2(0,T;V_\eps')$.
Lemma~\ref{lem:abstract-uniqueness} and \eqref{ass:Adu} give
\[
 c_u\int_{\Qeps}\sum_{i=1}^Nw_i^\eps|D_i(u^1-u^2)|^2\dd x\dd t
 \leq\int_0^T\langle\mathcal A_\eps(t)u^1-\mathcal A_\eps(t)u^2,u^1-u^2\rangle\dd t
 \leq0 .
\]
Since $w_i>0$ a.e., $D(u^1-u^2)=0$; since $\Omeps$ is connected and
$\gamma_0(u^1-u^2)=0$ on $\Gamma_2$, we get $u^1=u^2$.
\end{proof}

\begin{remark}[Non-uniqueness under $(H_9)$]\label{rem:non-uniqueness}
The uniqueness proof uses only that $B_\eps$ is symmetric and non-negative;
strict monotonicity in \((H_{10})\) is what controls the gradient difference.
Under \((H_9)\) alone, consider $w_i\equiv1$, $b=\chi_E$ with $E\subset\Ystar$ of positive measure such
that $\Ystar\setminus\overline E$ contains a non-empty open set,
$\beta(s)=s-\operatorname{proj}_{[-1,1]}s$, $A_0(y,\tau,r)=\beta(r)$,
$A_i(y,\tau,r,\xi)=\beta(\xi_i)$, and $f=0$, $g=0$, $\Phi=0$.  Since $\beta$ is
continuous, nondecreasing and $1$-Lipschitz, $\beta(r)r\geq0$ and
$\beta(s)s\geq\frac12s^2-\frac12$, hypotheses \((H_1)\)--\((H_9)\) hold, with
$c_0=\frac12$ and $k_2=\frac N2$ in \eqref{ass:A2}.  For small $\eps$ choose an open set
$D_\eps\Subset\Omeps$ on which $b^\eps=0$ and
$0\neq\eta\in C_c^\infty(D_\eps)$.  For every $a\in C^\infty([0,T])$ with
$a(0)=0$ and every $\lambda>0$ such that $|\lambda a\eta|\leq1$ and
$|\lambda aD_i\eta|\leq1$, $i=1,\ldots,N$, the function
$u=\lambda a(t)\eta(x)$ is a microscopic variational solution: indeed
$b^\eps u=0$, $A_0^\eps[u]=0$ and $A^\eps[u]=0$, so that both \eqref{eq:initial-mass} and
\eqref{eq:micro-interior-balance} hold with zero terms.  Moreover,
$\int_0^\sigma\langle B_\eps v',u\rangle\dd t=0$, so \eqref{eq:microVI}
reduces to $\frac12\int_{\Omeps}b^\eps|v(\sigma)|^2\dd x\geq0$.  Hence
\((H_9)\) does not imply uniqueness when $b$ vanishes on an open set.
\end{remark}

\section{The surface lifting and uniform estimates}\label{sec:surface}

The zero-mean condition in \eqref{ass:surface} permits a weighted Neumann
lifting on the cell.  The lifting space is not periodic; this imposes zero flux
on the outer faces and makes the cellwise surface-to-volume identity exact.  The
lifting gives an exact bulk representation of the surface functional, so that
volume two-scale convergence suffices (compare
\cite{AllaireDamlamianHornung1996,NeussRadu1996}).

\begin{lemma}
For a.e. $\tau\in Z$ there is a unique
$\vartheta(\tau)\in X_\omega(\Ystar)/\R$ such that
\begin{equation}
 \sum_{i=1}^N\int_{\Ystar}w_iD_i\vartheta(\tau)D_i\zeta\dd y
 =\langle g(\tau),\gamma_0^y \zeta\rangle_S
 \quad\forall\,\zeta\in X_\omega(\Ystar)/\R.                   \label{eq:lifting}
\end{equation}
If $F_i=w_iD_i\vartheta$, then
\begin{equation}
 \sum_{i=1}^N\int_Z\int_{\Ystar}\frac{|F_i|^2}{w_i}\dd y\dd\tau<\infty,
 \quad \Div_yF=0,\quad F\cdot n=g\text{ on }S,\quad
 F\cdot n=0\text{ on }\partial Y.                        \label{eq:Fproperties}
\end{equation}
The boundary equalities are weak equalities.
\end{lemma}

\begin{proof}
On $\dot X_\omega=X_\omega(\Ystar)/\mathbb R$ the weighted gradient
seminorm is a norm (Lemma~\ref{lem:cell-poincare}), and
$\|z-z_{\Ystar}\|_{W^{1,q}(\Ystar)}\leq C\|[z]\|_{\dot X_\omega}$.  Since
$\langle g(\tau),1\rangle_S=0$, the trace theorem gives a well-defined
functional on the quotient with
\begin{equation}
 |\langle g(\tau),\gamma_0^yz\rangle_S|
 \leq C\|g(\tau)\|_{(W^{1-1/q,q}(S))'}
          \|[z]\|_{\dot X_\omega}.                    \label{eq:lifting-functional}
\end{equation}
Lax--Milgram gives a unique $\vartheta(\tau)$ with
\begin{equation}
 \|[\vartheta(\tau)]\|_{\dot X_\omega}
 \leq C\|g(\tau)\|_{(W^{1-1/q,q}(S))'},             \label{eq:lifting-energy-pointwise}
\end{equation}
depending linearly and boundedly on $g(\tau)$, hence measurably on $\tau$;
this gives the energy bound in \eqref{eq:Fproperties}.  Testing
\eqref{eq:lifting} with $\zeta\in C_c^\infty(\Ystar)$, with functions vanishing
near $S$, and with functions vanishing near $\partial Y$ gives the three weak
relations in \eqref{eq:Fproperties}.
\end{proof}

\begin{lemma}\label{lem:exact-surface}
For every $z\in V_\eps$ and a.e. $t$,
\begin{equation}
 \langle\ell_\eps(t),z\rangle
 =\int_{G_\eps}F(x/\eps,t/\eps^2)\cdot Dz(x)\dd x.        \label{eq:exact-surface}
\end{equation}
Moreover,
\begin{equation}
 |\langle\ell_\eps(t),z\rangle|
 \leq C\|g(t/\eps^2)\|_{(W^{1-1/q,q}(S))'}
 \left(\sum_{i=1}^N\|\sqrt{w_i^\eps}D_iz\|_{L^2(\Omeps)}^2\right)^{1/2},
                                                               \label{eq:surface-bound}
\end{equation}
with an $L^2(0,T)$ bound independent of $\eps$.
\end{lemma}

\begin{proof}
Put $z_k(y)=z(\eps(k+y))$.  Equation \eqref{eq:lifting} and scaling give
\begin{align*}
 \eps^{N-1}\langle g(t/\eps^2),\gamma_0^y z_k\rangle_S
 &=\eps^{N-1}\int_{\Ystar}F(y,t/\eps^2)\cdot D_yz_k\dd y\\
 &=\int_{\eps(k+\Ystar)}F(x/\eps,t/\eps^2)\cdot D_xz\dd x.
\end{align*}
Summing over $k\in I_\eps$ proves
\eqref{eq:exact-surface}.  Write $F_i^\eps(x,t)=F_i(x/\eps,t/\eps^2)$.  Weighted Cauchy--Schwarz gives
\begin{equation}
 |\langle\ell_\eps(t),z\rangle|\leq
 \left(\sum_{i=1}^N\int_{G_\eps}\frac{|F_i^\eps|^2}{w_i^\eps}\dd x\right)^{1/2}
 \left(\sum_{i=1}^N\int_{G_\eps}w_i^\eps|D_i z|^2\dd x\right)^{1/2}. \label{eq:surface-CS}
\end{equation}
The first squared factor can be calculated cell by cell:
\begin{align*}
 \sum_{i=1}^N\int_{G_\eps}\frac{|F_i(x/\eps,t/\eps^2)|^2}{w_i^\eps(x)}\dd x
 &=\eps^N\# I_\eps\sum_{i=1}^N\int_{\Ystar}
   \frac{|F_i(y,t/\eps^2)|^2}{w_i(y)}\dd y\\
 &\leq C\|g(t/\eps^2)\|_{(W^{1-1/q,q}(S))'}^2.
\end{align*}
Here $\eps^N\#I_\eps\leq|\Om|$ and
\eqref{eq:lifting-energy-pointwise} were used.  Since $G_\eps\subset\Omeps$,
substitution in \eqref{eq:surface-CS} proves \eqref{eq:surface-bound}.
Moreover, periodicity in $\tau$ gives
\begin{equation}
 \int_0^T\|g(t/\eps^2)\|_{(W^{1-1/q,q}(S))'}^2\dd t
 \leq (T+1)\|g\|_{L^2_\#(Z;(W^{1-1/q,q}(S))')}^2
 \qquad\forall\,\eps\in(0,1).                                      \label{eq:phi-periodic-bound}
\end{equation}
Thus the coefficient in \eqref{eq:surface-bound} is bounded in $L^2(0,T)$
uniformly in $\eps$.
\end{proof}

\begin{theorem}\label{thm:uniform}
Under \((H_1)\)--\((H_8)\), every microscopic variational solution satisfies
\begin{align}
 &\int_{\Qeps}\left(|u_\eps|^2+
        \sum_{i=1}^Nw_i^\eps|D_iu_\eps|^2\right)\dd x\dd t\leq C,
                                                               \label{eq:uniform-energy}\\[1.2ex]
 &\left\|\partial_t(\chi_{\Omeps}b^\eps u_\eps)\right\|_
 {L^2(0,T;H^{-1}(\Om))}\leq C,                            \label{eq:uniform-time}
\end{align}
where $C$ is independent of $\eps$.
\end{theorem}

\begin{proof}
The comparison $\Psi$ belongs to $\mathcal K_\eps^1(T)$, since
$\gamma_0\Psi=\varphi(t)-\varphi(0)\leq\varphi(t)$ on $\Gamma_1$, and
$\Psi(0)=0$.  Using it in \eqref{eq:microVI} with $z_\eps=u_\eps-\Psi$ gives
\begin{align}
 \int_{\Qeps}\bigl(A_0^\eps[u_\eps]u_\eps
              +A^\eps[u_\eps]\cdot Du_\eps\bigr)
 &\leq\int_{\Qeps}\bigl(A_0^\eps[u_\eps]\Psi
              +A^\eps[u_\eps]\cdot D\Psi+fz_\eps\bigr)\dd x\dd t\nonumber\\
 &\quad+\int_0^T\langle\ell_\eps,z_\eps\rangle\dd t
 -\int_0^T\langle B_\eps\partial_t\Psi,z_\eps\rangle\dd t. \label{eq:integrated-energy-start}
\end{align}
By periodicity,
\begin{equation}
 \int_Q|k_1^\eps|^2\dd x\dd t
 +\int_Q|k_2^\eps|\dd x\dd t\leq C.                 \label{eq:periodic-data-bound}
\end{equation}
Let $\mathcal E_\eps$ denote the left-hand side of \eqref{eq:uniform-energy}
and $\mathcal G_\eps=\sum_{i=1}^N\int_{\Qeps}w_i^\eps|D_iu_\eps|^2$ its gradient part.
By Lemma~\ref{lem:uniform-poincare}, $\mathcal E_\eps\leq(1+C_P)\mathcal G_\eps$.
Put $\delta=c_0/(8(1+C_P))$.  In what follows $C$ depends on $\delta$, $c_1$,
$\|w\|_{L^\infty}$, $\|b\|_{L^\infty}$, \eqref{eq:periodic-data-bound},
$\|f\|_{L^2(Q)}$, $\|\Psi\|_{H^1(0,T;H^1(\Om))}$ and
$\|g\|_{L^2_\#(Z;(W^{1-1/q,q}(S))')}$, but not on $\eps$; all integrals are over
$\Qeps$ unless indicated.
\begin{itemize}
\item[(a)] \emph{Left-hand side.}  By \eqref{ass:A2} and
$c_1k_1^\eps|u_\eps|\leq\delta|u_\eps|^2+C|k_1^\eps|^2$, it is at least
$c_0\mathcal G_\eps-\delta\mathcal E_\eps-C$.
\item[(b)] \emph{Zero-order term.}  By \eqref{ass:A1a},
$|A_0^\eps[u_\eps]\Psi|\leq c_1(k_1^\eps+|u_\eps|)|\Psi|
\leq\frac\delta2|u_\eps|^2+C(|k_1^\eps|^2+|\Psi|^2)$, so its integral is at
most $\frac\delta2\mathcal E_\eps+C$.
\item[(c)] \emph{Flux term.}  By \eqref{ass:A1b}, with
$X_\eps=k_1^\eps+|u_\eps|+\sum_{j=1}^N\sqrt{w_j^\eps}|D_ju_\eps|$ and
$X_\eps^2\leq(N+2)(|k_1^\eps|^2+|u_\eps|^2+\sum_{j=1}^Nw_j^\eps|D_ju_\eps|^2)$,
\[
 |A^\eps[u_\eps]\cdot D\Psi|
 \leq c_1\|w\|_{L^\infty}^{1/2}\sqrt N\,|D\Psi|\,X_\eps
 \leq\frac{\delta}{2(N+2)}X_\eps^2+C|D\Psi|^2,
\]
so its integral is at most $\frac\delta2\mathcal E_\eps+C$.
\item[(d)] \emph{Source term.}  $|fz_\eps|\leq|f|(|u_\eps|+|\Psi|)
\leq\frac\delta4|u_\eps|^2+C(|f|^2+|\Psi|^2)$, so its integral is at most
$\frac\delta4\mathcal E_\eps+C$.
\item[(e)] \emph{Mass term.}  $|\langle B_\eps\partial_t\Psi,z_\eps\rangle|
\leq\|b\|_{L^\infty}\int_{\Omeps}|\partial_t\Psi|(|u_\eps|+|\Psi|)\dd x$, so,
as in (d), its time integral is at most $\frac\delta4\mathcal E_\eps+C$.
\item[(f)] \emph{Surface term.}  By \eqref{eq:surface-bound} applied to
$z_\eps(t)$, Young's inequality and \eqref{eq:phi-periodic-bound},
\[
 \int_0^T|\langle\ell_\eps,z_\eps\rangle|\dd t
 \leq\frac\delta2\mathcal G_\eps
 +\frac\delta2\sum_{i=1}^N\|\sqrt{w_i^\eps}D_i\Psi\|_{L^2(\Qeps)}^2
 +C\int_0^T\|g(t/\eps^2)\|_{(W^{1-1/q,q}(S))'}^2\dd t
 \leq\delta\mathcal E_\eps+C.
\]
\end{itemize}
Inserting (a)--(f) in \eqref{eq:integrated-energy-start} gives
$c_0\mathcal G_\eps\leq4\delta\mathcal E_\eps+C\leq\frac{c_0}2\mathcal G_\eps+C$,
by $\mathcal E_\eps\leq(1+C_P)\mathcal G_\eps$, which bounds $\mathcal G_\eps$
and then $\mathcal E_\eps$; this is \eqref{eq:uniform-energy}.

For \eqref{eq:uniform-time}, \eqref{eq:micro-interior-balance} and
\eqref{eq:exact-surface} give, for $\zeta\in H_0^1(\Om)$,
\[
 \langle\partial_t(\chi_{\Omeps}b^\eps u_\eps),\zeta\rangle
 =\int_{\Omeps}f\zeta
 -\int_{\Omeps}(A_0^\eps[u_\eps]\zeta+A^\eps[u_\eps]\cdot D\zeta)
 +\int_{G_\eps}F^\eps\cdot D\zeta .
\]
By weighted Cauchy--Schwarz, \eqref{ass:A1a}--\eqref{ass:A1b},
\eqref{eq:lifting-energy-pointwise} and the boundedness of the weights, the
right-hand side is bounded by $M_\eps(t)\|\zeta\|_{H_0^1(\Om)}$, where
$\|M_\eps\|_{L^2(0,T)}$ is controlled by $\|f\|_{L^2(Q)}$,
\eqref{eq:periodic-data-bound}, \eqref{eq:uniform-energy} and
\eqref{eq:phi-periodic-bound}.
\end{proof}

\section{Weighted space--time two-scale compactness}\label{sec:compactness}

Strong compactness of the states follows from the detection inequality below
and time-translation estimates.

\begin{lemma}\label{lem:detection}
Let $(\rho_\eps)$ be bounded in $L^\infty(\Om)$ with
$\rho_\eps\weakto^*\rho$ in $L^\infty(\Om)$ as $\eps\to0$, where $\rho>0$
a.e. in $\Om$.  For every $\eta>0$ there are $C_\eta>0$ and $\eps_\eta>0$
such that
\[
 \|v\|_{L^2(\Om)}
 \leq \eta\|v\|_{W^{1,q}(\Om)}
      +C_\eta\|\rho_\eps v\|_{H^{-1}(\Om)}
\]
for $v\in W^{1,q}(\Om)$ and $0<\eps<\eps_\eta$.
\end{lemma}

\begin{proof}
If the assertion failed for some $\eta>0$, there would be $\eps_n\downarrow0$
and $v_n\in W^{1,q}(\Om)$ with $\|v_n\|_{L^2(\Om)}=1$ and
$\eta\|v_n\|_{W^{1,q}(\Om)}+n\|\rho_{\eps_n}v_n\|_{H^{-1}(\Om)}<1$.  By the
compact embedding $W^{1,q}(\Om)\hookrightarrow\hookrightarrow L^2(\Om)$, after extraction
$v_n\to v$ in $L^2(\Om)$ with $\|v\|_{L^2(\Om)}=1$.  For
$\zeta\in C_c^\infty(\Om)$, $\int_\Om\rho_{\eps_n}(v_n-v)\zeta\to0$ since
$(\rho_{\eps})$ is bounded in $L^\infty(\Om)$, and
$\int_\Om\rho_{\eps_n}v\zeta\to\int_\Om\rho v\zeta$ by the weak-star
convergence.  Hence $\rho_{\eps_n}v_n\to\rho v$ in distributions, while
$\|\rho_{\eps_n}v_n\|_{H^{-1}(\Om)}<1/n$.  Thus $\rho v=0$, which contradicts
$\rho>0$ a.e. and $\|v\|_{L^2(\Om)}=1$.
\end{proof}

We apply Lemma~\ref{lem:detection} to the perforated periodic density
\[
 \rho_\eps(x)=\chi_{\Omeps}(x)b(x/\eps),
 \qquad \rho_\eps\weakto^*\overline b>0
       \quad\hbox{in }L^\infty(\Om),
\]
the weak-star limit following from periodic averaging and \eqref{eq:strip}.

\begin{theorem}\label{thm:weighted-compactness}
Let $u_\eps\in L^2(0,T;V_\eps)$ satisfy \eqref{eq:uniform-energy},
\eqref{eq:uniform-time} and $(\chi_{\Omeps}b^\eps u_\eps)(0)=0$ in
$H^{-1}(\Om)$; this holds in particular for microscopic variational solutions,
by \eqref{eq:initial-mass} and Theorem~\ref{thm:uniform}.  There are a subsequence and
functions
\begin{align*}
 u_0&\in L^2(0,T;H^1(\Om)),\qquad \gamma_0 u_0=0\text{ on }\Gamma _2,\\
 u_1&\in L^2\bigl(Q;L^2_\#(Z;V_\omega)\bigr)
\end{align*}
such that, with $U=(u_0,u_1)$ and
\begin{equation*}
 \mathbb D U=D_xu_0+D_yu_1,
 \qquad (\mathbb D U)_i=D_{x_i}u_0+D_{y_i}u_1,
 \quad i=1,\ldots,N,
\end{equation*}
\begin{equation}
 P_\eps u_\eps\longrightarrow u_0\quad\hbox{strongly in }L^2(Q),
 \qquad P_\eps u_\eps\tto u_0.                          \label{eq:compact-u}
\end{equation}
Moreover,
\begin{align}
 \chi_{\Omeps}D_i u_\eps&\twto
     \chi_{\Ystar}(\mathbb D U)_i &&\text{weakly in }L^q(Q),
                                                               \label{eq:compact-gradq}\\
 \chi_{\Omeps}\sqrt{w_i^\eps}D_i u_\eps&\twto
     \chi_{\Ystar}\sqrt{w_i(y)}(\mathbb D U)_i
     &&\text{weakly in }L^2(Q).                                 \label{eq:compact-weighted}
\end{align}
The limit $u_0$ is independent of both fast variables, and
\begin{equation}
 u_0\in H^1(0,T;H^{-1}(\Om)),\qquad
 \overline b\,\partial_tu_0\in L^2(0,T;H^{-1}(\Om)),\qquad
 u_0(0)=0\quad\hbox{in }H^{-1}(\Om).                     \label{eq:macro-time}
\end{equation}
\end{theorem}

\begin{proof}
Lemmas~\ref{lem:weighted-control} and \ref{lem:extension} give
\begin{equation}
 \|P_\eps u_\eps\|_{L^2(0,T;W^{1,q}(\Om))}\leq C.          \label{eq:Pbound}
\end{equation}
Put $\bar u_\eps=P_\eps u_\eps$.  Since $\bar u_\eps=u_\eps$ in $\Omeps$,
\begin{equation}
 \rho_\eps \bar u_\eps=\chi_{\Omeps}b^\eps u_\eps,
 \qquad
 \|\partial_t(\rho_\eps \bar u_\eps)\|_{L^2(0,T;H^{-1}(\Om))}\leq C.
                                                               \label{eq:rho-time}
\end{equation}
For $0<h<T$, writing
$\rho_\eps(\bar u_\eps(t+h)-\bar u_\eps(t))=\int_t^{t+h}\partial_s(\rho_\eps\bar u_\eps)\dd s$
in $H^{-1}(\Om)$ and using Young's convolution inequality, we get
\begin{equation}
 \|\rho_\eps(\bar u_\eps(\cdot+h)-\bar u_\eps)\|_
 {L^2(0,T-h;H^{-1}(\Om))}\leq Ch.                       \label{eq:rho-translate}
\end{equation}
Apply Lemma~\ref{lem:detection} for a.e. $t$ to
$\bar u_\eps(t+h)-\bar u_\eps(t)$, and then take the $L^2(0,T-h)$ norm.  The triangle
inequality, \eqref{eq:Pbound} and \eqref{eq:rho-translate} give
$\|\bar u_\eps(\cdot+h)-\bar u_\eps\|_{L^2(0,T-h;L^2(\Om))}\leq 2C\eta+C_\eta Ch$.
Choosing first $\eta$ and then $h$ small, the time translations tend to zero
uniformly in $\eps$ (the finitely many terms of the sequence with
$\eps\geq\eps_\eta$ are handled individually), and \cite[Theorem~5]{Simon1986}, applied with the compact embedding
$W^{1,q}(\Om)\hookrightarrow\hookrightarrow L^2(\Om)$, gives, after extraction, $\bar u_\eps\to u_0$ strongly in $L^2(Q)$,
hence also strongly two-scale.  Proposition~\ref{prop:classical-gradient-decomposition}
with $p=q$ gives a periodic corrector $\widetilde u_1$ such that
\begin{equation}
 D_i\bar u_\eps\twto D_{x_i}u_0+D_{y_i}\widetilde u_1
 \quad\hbox{weakly in }L^q(Q),\qquad i=1,\ldots,N.          \label{eq:extended-gradq}
\end{equation}
By \eqref{eq:strip}, $\chi_{\Omeps}\tto\chi_{\Ystar}$ in every $L^p(Q)$,
$p<\infty$; since $\bar u_\eps=u_\eps$ on $\Omeps$, this yields
\begin{equation}
 \chi_{\Omeps}D_i u_\eps\twto
 \chi_{\Ystar}\bigl(D_{x_i}u_0+D_{y_i}\widetilde u_1\bigr)
 \quad\hbox{weakly in }L^q(Q),
 \qquad i=1,\ldots,N.
 \label{eq:material-gradq-tilde}
\end{equation}

By \eqref{eq:uniform-energy}, after extraction,
\begin{equation}
 \chi_{\Omeps}\sqrt{w_i^\eps}D_i u_\eps
 \twto \Upsilon_i
 \quad\hbox{weakly in }L^2(Q),
 \qquad i=1,\ldots,N.
 \label{eq:g-preliminary}
\end{equation}
Approximating $\sqrt{w_i}$ by $h_m\in C_\#^\infty(Y)$ in $L^{q'}(Y)$,
$q'=q/(q-1)$, testing with $(h_m\psi)^\eps$ for smooth periodic $\psi$ in
\eqref{eq:material-gradq-tilde}, and noting that the error is bounded by
$C\|\sqrt{w_i}-h_m\|_{L^{q'}(Y)}\|\psi\|_\infty$ uniformly in $\eps$, we obtain
\begin{equation}
 \Upsilon_i=\chi_{\Ystar}\sqrt{w_i(y)}
 \bigl(D_{x_i}u_0+D_{y_i}\widetilde u_1\bigr)
 \quad\hbox{in }\mathcal D'(Q\times Y\times Z),
 \qquad i=1,\ldots,N.
 \label{eq:g-distribution}
\end{equation}

We now identify the $L^2$ structure of $\Upsilon$, using the affine cell inequality
\eqref{eq:cell-affine}.  Set
$\mathcal H_\omega^{\rm aff}=\prod_{i=1}^NL^2(\Ystar)$ and define
\begin{equation}
 \mathcal T_\omega:\mathbb R^N\times V_\omega\longrightarrow\mathcal H_\omega^{\rm aff},
 \qquad
 \mathcal T_\omega(\xi,z)=
 \left(\sqrt{w_i}\,(\xi_i+D_i z)\right)_{i=1}^N.
\end{equation}
By \eqref{eq:cell-affine},
\begin{align}
 |\xi|^2+\|z\|_{L^2(\Ystar)}^2
 +\sum_{i=1}^N\int_{\Ystar}w_i|D_i z|^2\dd y
 \leq C\|\mathcal T_\omega(\xi,z)\|_{\mathcal H_\omega^{\rm aff}}^2.
 \label{eq:affine-range-lower}
\end{align}
Consequently $\mathcal T_\omega$ is injective and its range
$\mathcal R_\omega:=\mathcal T_\omega(\mathbb R^N\times V_\omega)
\subset\mathcal H_\omega^{\rm aff}$ is a closed subspace.

We claim that
\begin{equation}
 \Upsilon(x,t,\tau,\cdot)\in\mathcal R_\omega
 \quad\hbox{for a.e. }(x,t,\tau)\in Q\times Z.
 \label{eq:g-in-affine-range}
\end{equation}
Let $\Xi\in\mathcal R_\omega^\perp$, put $H_i=\sqrt{w_i}\,\Xi_i$ on $\Ystar$,
and extend $\Xi$ and $H$ by zero to $T_0$ and then periodically.
Orthogonality to $\mathcal T_\omega(\xi,0)$, $\xi\in\R^N$, and to
$\mathcal T_\omega(0,z)$, $z\in V_\omega$, gives
\begin{equation}
 \int_{\Ystar}H_i\dd y=0\quad(i=1,\ldots,N),           \label{eq:H-zero-mean}
 \qquad
 \sum_{i=1}^N\int_{\Ystar}H_iD_i z\dd y=0\quad\forall\,z\in V_\omega.
\end{equation}
For any $\zeta\in C_\#^\infty(Y)$, the function
$\zeta-\overline b^{-1}\int_{\Ystar}b\zeta$ belongs to $V_\omega$
and has the same gradient as $\zeta$.  Thus the second identity, after
extension of $H$ by zero to $T_0$, is exactly
\begin{equation}
 \int_Y H\cdot D_y\zeta\dd y=0\quad\forall\,\zeta\in C_\#^\infty(Y),
 \qquad \Div_yH=0\quad\hbox{in }\mathcal D'_\#(Y).
 \label{eq:H-divfree}
\end{equation}
Here $H\in L^2_\#(Y)^N$ since $w_i\in L^\infty(Y)$, and the periodic
distributional identity includes both the zero normal flux on $S$ and
cancellation of the traces across opposite faces of $Y$.

Let $\alpha\in C_c^\infty(Q)$ and $\beta\in C_\#^\infty(Z)$.  The weighted
Green formula gives
\begin{align}
 &\sum_{i=1}^N\int_Q\alpha(x,t)\beta(t/\eps^2)
 \chi_{\Omeps}(x)\sqrt{w_i^\eps(x)}D_i u_\eps(x,t)
 \Xi_i(x/\eps)\dd x\dd t                                      \nonumber\\
 &\qquad=-\int_Q\bar u_\eps(x,t)\beta(t/\eps^2)
 H(x/\eps)\cdot D_x\alpha(x,t)\dd x\dd t .
\end{align}
Indeed, for smooth $u_\eps$ the periodic distributional identity
\eqref{eq:H-divfree} applies to $\alpha u_\eps$ after extension across the
holes; there are no interior boundary terms.  The identity extends by
density to $H^1_\omega(\Omeps)$: the left-hand side is bounded by
$\|\sqrt{w_i^\eps}D_i u_\eps\|_{L^2}$ times
$\|\Xi_i\|_{L^2_\#(Y)}$, and the right-hand side by
$\|u_\eps\|_{L^2}$ times $\|H\|_{L^2_\#(Y)}$ (with constants for
fixed $\eps$).  Since $H^\eps=0$ inside the holes, $u_\eps$ can then
be replaced by $\bar u_\eps$ on the right.

Passing to the limit, on the left by \eqref{eq:g-preliminary} (the sequence $\alpha\beta^\eps\Xi(x/\eps)$, with $\Xi\in L^2_\#(Y)$,
converges strongly two-scale, because the mean-value property of products of
continuous and $L^1_\#$ functions \cite[Section~5]{Allaire1992} gives
convergence of its $L^2(Q)$ norm; hence the strong--weak product rule
\cite[Theorem~1.8]{Allaire1992} applies) and on the right by the strong convergence of
$\bar u_\eps$, we obtain
\begin{align*}
 &\sum_{i=1}^N\int_Q\int_Z\int_{\Ystar}
 \alpha(x,t)\beta(\tau)\Upsilon_i(x,t,y,\tau)\Xi_i(y)
       \dd y\dd\tau\dd x\dd t\\
 &\qquad=-\int_Q\int_Zu_0(x,t)\beta(\tau)D_x\alpha(x,t)\cdot
 \left(\int_{\Ystar}H(y)\dd y\right)\dd\tau\dd x\dd t=0,
\end{align*}
where the last equality follows from \eqref{eq:H-zero-mean}.
The identities hold first for a countable dense family of
$\alpha\in C_c^\infty(Q)$ and $\beta\in C_\#^\infty(Z)$; choosing also a
countable dense subset of the separable Hilbert space
$\mathcal R_\omega^\perp$ and then using density gives
\[
 \Upsilon(x,t,\tau,\cdot)\perp\mathcal R_\omega^\perp
 \quad\hbox{for a.e. }(x,t,\tau).
\]
Since $\mathcal R_\omega$ is closed, this proves
\eqref{eq:g-in-affine-range}.

Therefore there exist unique measurable
$\xi:Q\times Z\to\R^N$ and $u_1:Q\times Z\to V_\omega$ such that
\begin{equation}
 \Upsilon_i(x,t,y,\tau)=\chi_{\Ystar}(y)\sqrt{w_i(y)}
 \bigl(\xi_i(x,t,\tau)+D_{y_i}u_1(x,t,y,\tau)\bigr),
 \qquad i=1,\ldots,N.
 \label{eq:g-affine-representation}
\end{equation}
Indeed, the inverse of
$\mathcal T_\omega:\mathbb R^N\times V_\omega\to\mathcal R_\omega$ is bounded by
\eqref{eq:affine-range-lower}; hence measurability follows from measurability
of $\Upsilon$, and
\begin{equation}
 \|\xi\|_{L^2(Q\times Z)^N}^2+\|u_1\|_{L^2(Q\times Z;V_\omega)}^2
 \leq C\sum_{i=1}^N\|\Upsilon_i\|_{L^2(Q\times Y\times Z)}^2<\infty. \label{eq:affine-representation-estimate}
\end{equation}

It remains to identify $\xi$ and to compare $u_1$ with the corrector
$\widetilde u_1$ obtained from the $W^{1,q}$ two-scale decomposition.
Comparing \eqref{eq:g-distribution} and \eqref{eq:g-affine-representation} on
$\Ystar$ and using $w_i>0$ a.e., we get
\begin{equation}
 D_{y_i}(u_1-\widetilde u_1)=D_{x_i}u_0-\xi_i,
 \qquad i=1,\ldots,N.
 \label{eq:constant-gradient-difference}
\end{equation}
The right-hand side is independent of $y$, and a periodic function on the
connected cell cannot have a non-zero constant gradient (see the proof of
Lemma~\ref{lem:cell-poincare}).  Hence
\begin{equation}
 \xi(x,t,\tau)=D_xu_0(x,t)\quad\hbox{for a.e. }(x,t,\tau).
 \label{eq:xi-equals-macro-gradient}
\end{equation}
In particular, \eqref{eq:affine-representation-estimate} yields
$D_xu_0\in L^2(Q)^N$, that is, $u_0\in L^2(0,T;H^1(\Om))$.

After choosing the representative of $\widetilde u_1$ with zero $b$-mean on
$\Ystar$, equation \eqref{eq:constant-gradient-difference} and
\eqref{eq:xi-equals-macro-gradient} give $u_1=\widetilde u_1$ on $\Ystar$.
Together with \eqref{eq:material-gradq-tilde} and
\eqref{eq:g-affine-representation}, this proves
\eqref{eq:compact-gradq}--\eqref{eq:compact-weighted}.
Since
$\bar u_\eps\rightharpoonup u_0$ in $L^2(0,T;W^{1,q}(\Om))$ and
$\gamma_0 \bar u_\eps=0$ on $\Gamma _2$, continuity of the trace and weak closedness of
the zero-trace subspace give $\gamma_0 u_0=0$ on $\Gamma _2$.

Finally, $\rho_\eps \bar u_\eps\rightharpoonup\overline b\,u_0$ in $L^2(Q)$, since
$\bar u_\eps\to u_0$ strongly and $\rho_\eps\weakto^*\overline b$.  By
\eqref{eq:rho-time}, $(\rho_\eps \bar u_\eps)$ is bounded in $H^1(0,T;H^{-1}(\Om))$
and vanishes at $t=0$; hence $\overline b\,\partial_tu_0\in L^2(0,T;H^{-1}(\Om))$
and $\overline b\,u_0(0)=0$ in $H^{-1}(\Om)$, which gives \eqref{eq:macro-time}
since $\overline b>0$.
\end{proof}

The next lemma handles the factor $1/\eps$ that appears when the
microscopic time derivative is tested with an $O(\eps)$ corrector.  It is a
perforated, weighted version of the ``very weak'' corrector convergence used
in \cite[Lemma~2.9]{DouanlaSteklov2013}; see also the very weak multiscale
convergence of \cite{FlodenHolmbomOlssonPersson2010}.

\begin{lemma}\label{lem:very-weak}
Let $u_\eps,u_0,u_1$ be as in Theorem~\ref{thm:weighted-compactness}.  If
\begin{equation}
 \psi\in C_c^\infty(Q;C^\infty_\#(Z;C^\infty(\overline{\Ystar}))),
 \qquad \int_{\Ystar}b(y)\psi(x,t,y,\tau)\dd y=0,
\end{equation}
then, with the cellwise convention for $\psi^\eps$ specified below,
\begin{equation}
 \frac1\eps\int_{\Qeps}b^\eps u_\eps
       \psi^\eps\dd x\dd t
 \longrightarrow
 \int_Q\int_Z\int_{\Ystar}bu_1\psi\dd y\dd\tau\dd x\dd t.       \label{eq:very-weak}
\end{equation}
\end{lemma}

\begin{proof}
Extend $\psi(x,t,\cdot,\tau)$ by zero to $T_0$ and periodically; since the
$(x,t)$-support of $\psi$ is compact in $Q$, only complete cells contribute for
small $\eps$.  Let $q'=q/(q-1)$.  The zero $b$-mean of $\psi$ allows us to
apply a fixed linear Bogovskii operator on $\Ystar$
\cite[Chapter~III, Theorem~3.1]{Galdi2011}: there is
\begin{equation}
 \boldsymbol{\Theta}=\mathcal B_{\Ystar}(b\psi)\in W^{1,q'}_0(\Ystar)^N,
 \qquad \Div_y\boldsymbol{\Theta}=b\psi,                        \label{eq:bogovskii}
\end{equation}
with $\|\boldsymbol{\Theta}\|_{W^{1,q'}(\Ystar)}\leq C\|b\psi\|_{L^{q'}(\Ystar)}$,
and similarly for its derivatives in $(x,t,\tau)$.  Extended by zero to $T_0$
and periodically, $\boldsymbol\Theta$ has no jump across $S$ or $\partial Y$.

The chain rule gives
\begin{equation*}
 \Div_x\bigl[\boldsymbol{\Theta}(x,t,x/\eps,t/\eps^2)\bigr]
 =(\Div_x\boldsymbol{\Theta})^\eps+\frac1\eps(\Div_y\boldsymbol{\Theta})^\eps.
\end{equation*}
Because $\Div_y\boldsymbol{\Theta}=b\psi$ in $\Ystar$ and both sides are zero in $T_0$,
integration by parts on $\Om$ gives the exact identity
\begin{align}
 \frac1\eps\int_{\Qeps}b^\eps u_\eps\psi^\eps
 &=-\sum_{i=1}^N\int_QD_i(P_\eps u_\eps)\boldsymbol{\Theta}_i^\eps\dd x\dd t
   -\int_QP_\eps u_\eps(\Div_x\boldsymbol{\Theta})^\eps\dd x\dd t.   \label{eq:very-weak-identity}
\end{align}
The compositions of $\boldsymbol{\Theta}$ and $\Div_x\boldsymbol{\Theta}$ are admissible
tests in $L^{q'}$ by approximation and the uniform bound \eqref{eq:Pbound}.
Therefore
\eqref{eq:extended-gradq}, the strong $L^2$ convergence in
\eqref{eq:compact-u}, and \eqref{eq:very-weak-identity} give the limit
\begin{equation}
 -\int_Q\int_Z\int_{\Ystar}\bigl[(D_xu_0+D_yu_1)\cdot\boldsymbol{\Theta}+u_0\Div_x\boldsymbol{\Theta}\bigr]
 \dd y\dd\tau\dd x\dd t. \label{eq:very-weak-limit-expanded}
\end{equation}
For each fixed $(y,\tau)$, compact support in $x$ gives
\begin{equation*}
 \sum_{i=1}^N\int_\Om D_{x_i}u_0\boldsymbol{\Theta}_i\dd x
 +\int_\Om u_0\Div_x\boldsymbol{\Theta}\dd x=0.
\end{equation*}
These two parts cancel in \eqref{eq:very-weak-limit-expanded}.  Finally,
$\boldsymbol{\Theta}=0$ on $\partial\Ystar$ and \eqref{eq:bogovskii} imply
\begin{align*}
 -\sum_{i=1}^N\int_{\Ystar}D_{y_i}u_1\boldsymbol{\Theta}_i\dd y
 &=\int_{\Ystar}u_1\Div_y\boldsymbol{\Theta}\dd y
 =\int_{\Ystar}bu_1\psi\dd y.
\end{align*}
Integration over $(x,t,\tau)$ proves \eqref{eq:very-weak}.
\end{proof}

\begin{proposition}
\label{prop:micro-equation}
Let $u_\eps$ be microscopic variational solutions, in particular satisfying
the interior balance \eqref{eq:micro-interior-balance}, and let $u_0,u_1$ be
given by Theorem~\ref{thm:weighted-compactness}.  Let
\begin{align}
 a_0^\eps&=\chi_{\Omeps}A_0^\eps(x,t,u_\eps),\nonumber\\
 a_i^\eps&=\chi_{\Omeps}\frac{A_i^\eps(x,t,u_\eps,Du_\eps)}
                              {\sqrt{w_i^\eps}},\quad i=1,\ldots,N.
\end{align}
Then
\begin{equation}
 a_0^\eps\tto\chi_{\Ystar}A_0(y,\tau,u_0)
 \quad\hbox{strongly in }L^2(Q),                           \label{eq:A0-solution-limit}
\end{equation}
and, along a subsequence, $a_i^\eps\twto a_i$ weakly in $L^2(Q)$ for
$i=1,\ldots,N$.  Moreover,
$u_1\in L^2(Q;\mathcal W_\omega)$ and, for a.e. $(x,t)$,
\begin{equation}
 \int_Z\langle\partial_\tau(bu_1),\psi\rangle\dd\tau
 +\int_Z\int_{\Ystar}\sum_{i=1}^Na_i\sqrt{w_i}D_{y_i}\psi
 =\int_Z\langle g(\tau),\gamma_0^y \psi(\tau)\rangle_S\dd\tau
 \qquad\forall\,\psi\in\mathcal S_b,                       \label{eq:micro-unidentified}
\end{equation}
where $\mathcal S_b$ is the set of restrictions to $\Ystar\times Z$ of
functions $\psi\in C^\infty_\#(Y\times Z)$ such that
$\int_{\Ystar}b(y)\psi(y,\tau)\dd y=0$ for every $\tau\in Z$.
\end{proposition}

\begin{proof}
Put $\bar u_\eps=P_\eps u_\eps$.  Since $\bar u_\eps=u_\eps$ on $\Omeps$,
\begin{equation*}
 a_0^\eps=\chi_{\Omeps}A_0(x/\eps,t/\eps^2,\bar u_\eps).
\end{equation*}
Equation \eqref{eq:compact-u} and Proposition~\ref{prop:corrected31} show
that the second factor converges strongly two-scale to
$A_0(y,\tau,u_0)$.  Since
$\chi_{\Omeps}\tto\chi_{\Ystar}$ strongly in every finite $L^p$, the product
rule proves \eqref{eq:A0-solution-limit}.  Moreover,
\eqref{ass:A1b} gives, for a.e. $t$,
\begin{align*}
 \sum_{i=1}^N\|a_i^\eps(t)\|_{L^2(\Om)}^2
 \leq C\int_{\Omeps}\left(|k_1^\eps|^2+|u_\eps|^2
       +\sum_{j=1}^Nw_j^\eps|D_ju_\eps|^2\right)\dd x.
\end{align*}
Thus \eqref{eq:uniform-energy} and \eqref{eq:periodic-data-bound} make each
$a_i^\eps$ bounded in $L^2(Q)$, and two-scale compactness gives, after one
common extraction, $a_i^\eps\twto a_i$ for $i=1,\ldots,N$.

Choose a smooth
test $\psi(x,t,y,\tau)$ compactly supported in $(x,t)$, periodic in
$(y,\tau)$, and with zero $b$-mean.  Since its spatial support is compactly
contained in $\Om$, the function $\eps\psi^\eps$ is an admissible test in
the interior balance \eqref{eq:micro-interior-balance}.  We obtain
\begin{equation}
\begin{aligned}
 &\eps\int_0^T\langle
   \partial_t(\chi_{\Omeps}b^\eps u_\eps),\psi^\eps\rangle\dd t
 +\eps\int_{\Qeps}A_0^\eps[u_\eps]\psi^\eps\dd x\dd t
 +\int_{\Qeps}A^\eps[u_\eps]\cdot(D_y\psi+\eps D_x\psi)^\eps\dd x\dd t\\
 &\qquad=\eps\int_{\Qeps}f\psi^\eps\dd x\dd t
 +\int_0^T\langle\ell_\eps(t),\eps\psi^\eps(t)\rangle\dd t.
\end{aligned}
\end{equation}
Integration by parts in time transforms the first term into
\begin{equation*}
 -\eps\int_{\Qeps}b^\eps u_\eps(\partial_t\psi)^\eps\dd x\dd t
 -\frac1\eps\int_{\Qeps}b^\eps u_\eps(\partial_\tau\psi)^\eps
       \dd x\dd t.
\end{equation*}
The first integral is $O(\eps)$ by \eqref{eq:uniform-energy}.  Since
$\int_{\Ystar}b\,\partial_\tau\psi\dd y=0$, Lemma~\ref{lem:very-weak}
applied to $\partial_\tau\psi$ sends the second one to
\begin{equation*}
 -\int_Q\int_Z\int_{\Ystar}bu_1\partial_\tau\psi
       \dd y\dd\tau\dd x\dd t
 =\int_Q\int_Z\langle\partial_\tau(bu_1),\psi\rangle
       \dd\tau\dd x\dd t.
\end{equation*}
The terms containing $A_0^\eps$, $f$, and
$\eps D_{x_i}\psi$, $i=1,\ldots,N$, tend to zero by their $L^2$ bounds and
the explicit factor $\eps$.  For the remaining principal term, write
\begin{equation*}
 \sum_{i=1}^N\int_Qa_i^\eps
    \chi_{\Omeps}\sqrt{w_i^\eps}(D_{y_i}\psi)^\eps\dd x\dd t;
\end{equation*}
weak two-scale convergence of $a_i^\eps$ and admissibility of the fixed
periodic factor give
\begin{equation*}
 \int_Q\int_Z\int_{\Ystar}
       \sum_{i=1}^Na_i\sqrt{w_i}D_{y_i}\psi
       \dd y\dd\tau\dd x\dd t.
\end{equation*}
For the surface term, \eqref{eq:exact-surface} gives
\begin{align*}
 \int_0^T\langle\ell_\eps(t),\eps\psi^\eps(t)\rangle\dd t
 &=\int_{G_\eps\times(0,T)}F^\eps\cdot
       (D_y\psi+\eps D_x\psi)^\eps\dd x\dd t.
\end{align*}
The slow part is $O(\eps)$ by \eqref{eq:Fproperties}; periodic averaging
sends the other part to
\begin{align*}
 \int_Q\int_Z\int_{\Ystar}F\cdot D_y\psi
 &=\int_Q\int_Z\langle g(\tau),\gamma_0^y \psi\rangle_S
       \dd\tau\dd x\dd t,
\end{align*}
where the equality is \eqref{eq:lifting}.  We have therefore proved the
integrated form of \eqref{eq:micro-unidentified} for every smooth macro test.

For $z\in V_\omega$, \eqref{eq:lifting-functional} and Cauchy--Schwarz give
\begin{align*}
 &\left|\langle g(\tau),\gamma_0^y z\rangle_S
 -\int_{\Ystar}\sum_{i=1}^Na_i(x,t,y,\tau)
                     \sqrt{w_i(y)}D_{y_i}z(y)\dd y\right|\\
 &\quad\leq C\left(\|g(\tau)\|_{(W^{1-1/q,q}(S))'}
       +\left(\sum_{i=1}^N\|a_i(x,t,\cdot,\tau)\|_{L^2(\Ystar)}^2
        \right)^{1/2}\right)\|z\|_{V_\omega}.
\end{align*}
The coefficient belongs to $L^2(Q\times Z)$.  Hence the right-hand side of
\eqref{eq:micro-unidentified} minus its flux term defines an element of
$L^2(Q\times Z;V_\omega')$.  Fubini's theorem and density of the smooth
tests identify it with $\partial_\tau(bu_1)$, proving
$u_1\in L^2(Q;\mathcal W_\omega)$ and the pointwise-in-$(x,t)$ assertion.
\end{proof}

\begin{theorem}\label{thm:surface-limit}
Let $z_\eps\in L^2(0,T;V_\eps)$ satisfy
\begin{equation*}
 \sup_{0<\eps<1}\int_{\Qeps}\left(
 |z_\eps|^2+\sum_{i=1}^Nw_i^\eps|D_i z_\eps|^2\right)\dd x\dd t<\infty.
\end{equation*}
Assume that, along a subsequence, its leading two-scale limit
$z_0=z_0(x,t)$ is independent of $(y,\tau)$ and that its weighted gradients
converge as in \eqref{eq:compact-weighted}, with corrector $z_1$.  Then
\begin{equation}
 \int_0^T\langle\ell_\eps(t),z_\eps(t)\rangle\dd t\longrightarrow
 \int_Q\int_Z\langle g(\tau),\gamma_0^y z_1\rangle_S\dd\tau\dd x\dd t
 +\int_Qm_g\cdot D_xz_0\dd x\dd t. \label{eq:surface-limit}
\end{equation}
\end{theorem}

\begin{proof}
By Lemma~\ref{lem:exact-surface}, the left side is
\begin{equation}
 \sum_{i=1}^N\int_Q
 \chi_{G_\eps}(x)\frac{F_i(x/\eps,t/\eps^2)}{\sqrt{w_i^\eps(x)}}
       \sqrt{w_i^\eps(x)}D_iz_\eps\dd x\dd t.           \label{eq:surface-product}
\end{equation}
For each $i=1,\ldots,N$, put
\[
 \widehat F_i(y,\tau)=\frac{F_i(y,\tau)}{\sqrt{w_i(y)}}
             =\sqrt{w_i(y)}D_i\vartheta(y,\tau)
 \quad\hbox{in }\Ystar\times Z,
\]
and extend $\widehat F_i$ by zero to $T_0\times Z$ and periodically.  Then
$\widehat F_i\in L^2_\#(Y\times Z)$, and since the incomplete cells lie in a strip of
measure $O(\eps)$, periodic averaging gives
$\chi_{G_\eps}F_i^\eps/\sqrt{w_i^\eps}\tto\chi_{\Ystar}F_i/\sqrt{w_i}$.  By
\eqref{eq:strip}, $\chi_{G_\eps}\sqrt{w_i^\eps}D_iz_\eps$ has the same weak
two-scale limit as in \eqref{eq:compact-weighted}.  The strong--weak product rule sends \eqref{eq:surface-product} to
\begin{equation}
 \int_Q\int_Z\int_{\Ystar}F(y,\tau)\cdot
              (D_xz_0+D_yz_1)\dd y\dd\tau\dd x\dd t.     \label{eq:surface-intermediate}
\end{equation}
For a.e. $(x,t,\tau)$, approximate
$z_1(x,t,\cdot,\tau)\in V_\omega$ in its weighted norm by smooth periodic
functions and pass to the limit in \eqref{eq:lifting}; this gives
\begin{equation*}
 \int_{\Ystar}F\cdot D_yz_1\dd y
 =\langle g(\tau),\gamma_0^y z_1\rangle_S.
\end{equation*}
For the macroscopic part, the non-periodic function $y_i$ belongs to
$X_\omega(\Ystar)$, and \eqref{eq:lifting} gives, for every
$i=1,\ldots,N$,
\begin{equation*}
 \int_{\Ystar}F_i(y,\tau)\dd y
 =\langle g(\tau),y_i|_S\rangle_S.
\end{equation*}
Therefore
\begin{equation*}
 \int_Z\int_{\Ystar}F\cdot(D_xz_0+D_yz_1)\dd y\dd\tau
 =m_g\cdot D_xz_0+\int_Z\langle g(\tau),\gamma_0^y z_1\rangle_S\dd\tau.
\end{equation*}
Integrating this equality over $Q$ in
\eqref{eq:surface-intermediate} proves \eqref{eq:surface-limit}.
\end{proof}

\section{The homogenized problem}\label{sec:limit}

\subsection{The coupled two-scale problem}\label{sec:coupled}

Put $\theta=|\Ystar|$ and
\[
 H^1_{\Gamma _2}(\Om)=\{v\in H^1(\Om):\gamma_0 v=0
 \text{ on }\Gamma _2\}.
\]
The natural energy space and the regular comparison space are
\begin{align}
 \mathbb X_\omega={}&L^2(0,T;H^1_{\Gamma _2}(\Om))
 \times L^2\bigl(Q;L^2_\#(Z;V_\omega)\bigr),             \\
 \mathbb E_\omega={}&\{(v_0,v_1)\in\mathbb X_\omega:
 \partial_tv_0\in L^2\bigl(0,T;(H^1_{\Gamma _2}(\Om))'\bigr),
 \ v_1\in L^2(Q;\mathcal W_\omega)\}.
\end{align}
For $V=(v_0,v_1)$ write
\begin{equation*}
 \mathbb D V=D_xv_0+D_yv_1,
 \qquad (\mathbb D V)_i=D_{x_i}v_0+D_{y_i}v_1,
 \quad i=1,\ldots,N,
\end{equation*}
and set
\begin{equation}
 \mathbb K_\omega=\{V\in\mathbb X_\omega:
       \gamma_0 v_0\leq\varphi\text{ on }\Gamma _1\times(0,T)\}.
\end{equation}
For $V=(v_0,v_1)\in\mathbb E_\omega$ and
$\widehat V=(\hat v_0,\hat v_1)\in\mathbb X_\omega$, define
\begin{align}
 \langle\mathbb T V,\widehat V\rangle
 &:=\int_0^T\langle\overline b\,\partial_tv_0,\hat v_0\rangle\dd t
 +\int_Q\int_Z\langle\partial_\tau(bv_1),\hat v_1\rangle\dd\tau\dd x\dd t,
                                                               \\
 \langle\mathbb A(U),\widehat V\rangle
 &:=\int_Q\int_Z\int_{\Ystar}\bigl[A_0(y,\tau,u_0)\hat v_0
 +A(y,\tau,u_0,\mathbb D U)\cdot\mathbb D \widehat V\bigr]\dd y\dd\tau\dd x\dd t,
 \\
\langle\mathbb L,\widehat V\rangle
 &:=\theta\int_Qf\hat v_0\dd x\dd t+\int_{Q\times Z}\langle g(\tau),\gamma_0^y \hat v_1\rangle_S\dd\tau\dd x\dd t
 +\int_Qm_g\cdot D_x\hat v_0\dd x\dd t.
\end{align}
The operator $\mathbb T$ is defined on the comparison space
$\mathbb E_\omega$; the limit $U$ need not belong to it, since only
$\partial_tu_0\in L^2(0,T;H^{-1}(\Om))$ is known.

We use the fixed closed convex set obtained by translating the obstacle with
$\Psi(t)=\Phi(t)-\Phi(0)$, as in Section~\ref{sec:micro}:
\begin{align*}
 X_\Gamma&=L^2(0,T;H^1_{\Gamma_2}(\Om)),\qquad
 C_0=\{z\in H^1_{\Gamma_2}(\Om):\gamma_0z\leq\varphi(0)\text{ on }\Gamma_1\},\\
 \mathcal K_0&=\{z\in X_\Gamma:z(t)\in C_0\text{ for a.e. }t\},\qquad
 \mathcal D_0=\{z\in\mathcal K_0:\partial_tz\in X_\Gamma',\ z(0)=0\},
\end{align*}
where $X_\Gamma'=L^2(0,T;(H^1_{\Gamma_2}(\Om))')$, time derivatives are taken
in the Gelfand triple
$H^1_{\Gamma_2}(\Om)\subset L^2(\Om)\subset(H^1_{\Gamma_2}(\Om))'$, and
$\langle G,z\rangle_{X_\Gamma}=\int_0^T\langle G(t),z(t)\rangle\dd t$.  The set
$C_0$ is non-empty, closed and convex, $0\in C_0$ because $\varphi(0)\geq0$, and
$\Psi\in H^1(0,T;H^1_{\Gamma_2}(\Om))$.  Since
$\gamma_0\Psi(t)=\varphi(t)-\varphi(0)$ on $\Gamma_1$,
\[
 V=(v_0,v_1)\in\mathbb K_\omega\quad\Longleftrightarrow\quad v_0-\Psi\in\mathcal K_0 .
\]

\begin{lemma}\label{lem:limit-density}
If $V=(v_0,v_1)\in\mathbb K_\omega\cap\mathbb E_\omega$ and $v_0(0)=0$, then $V$ can
be approximated in the graph norm
\begin{equation*}
 \|V\|_{\mathrm{gr}}:=\|V\|_{\mathbb X_\omega}+\|\partial_tv_0\|_{L^2(0,T;(H^1_{\Gamma _2}(\Om))')}+\|\partial_\tau(bv_1)\|_{L^2(Q\times Z;V_\omega')}
\end{equation*}
by comparisons $V^n=(v_0^n,v_1^n)$ with
$v_0^n\in H^1(0,T;H^1_{\Gamma _2}(\Om))$, $v_0^n(0)=0$,
$v_0^n-\Psi\in\mathcal K_0$, and
\begin{equation}
 v_1^n(x,t,y,\tau)
 =\sum_{\ell=1}^{L_n}\alpha_{\ell n}(x,t)z_{\ell n}(y,\tau),
\end{equation}
where $\alpha_{\ell n}\in C_c^\infty(Q)$ and
$z_{\ell n}$ is the restriction to $\Ystar\times Z$ of a function in
$C_\#^\infty(Y\times Z)$ with zero $b$-mean, for
$\ell=1,\ldots,L_n$.
\end{lemma}

\begin{proof}
The function $z_0=v_0-\Psi$ belongs to $\mathcal K_0$ and satisfies $z_0(0)=0$,
so its extension by zero to $t<0$ keeps a time derivative in
$L^2(-1,T;(H^1_{\Gamma_2}(\Om))')$.  With a non-negative one-sided mollifier $\mu_n$
supported in $(0,1/n)$, the approximants $v_0^n=\Psi+\mu_n*z_0$ satisfy the
obstacle, $v_0^n(0)=0$, and converge in the graph norm.  For $v_1$ apply
Lemma~\ref{lem:graph-core}.
\end{proof}

\begin{theorem}\label{thm:main}
Under \((H_1)\)--\((H_8)\), let $\eps_n\downarrow0$ and assume that, for each $n$,
$u_{\eps_n}$ is a microscopic variational solution in the sense of
Section~\ref{sec:micro}.  Along the subsequence of
Theorem~\ref{thm:weighted-compactness}, the limit $U=(u_0,u_1)$ belongs to
$\mathbb K_\omega$, satisfies $u_1\in L^2(Q;\mathcal W_\omega)$
(Proposition~\ref{prop:micro-equation}), and solves the coupled two-scale Minty comparison problem
\begin{align}
 &\langle\mathbb T V-\mathbb L,V-U\rangle
 +\int_Q\int_Z\int_{\Ystar}\bigl[A_0(y,\tau,u_0)(v_0-u_0)\nonumber\\
 &\hspace{36mm}+A(y,\tau,u_0,\mathbb D V)\cdot(\mathbb D V-\mathbb D U)\bigr]
 \dd y\dd\tau\dd x\dd t\geq0\nonumber\\
 &\hspace{36mm}\forall\,V=(v_0,v_1)\in\mathbb K_\omega\cap\mathbb E_\omega\text{ with }v_0(0)=0.
                                                     \label{eq:two-scale-Minty}
\end{align}
\end{theorem}

\begin{proof}
For every $n$, Theorem~\ref{thm:uniform} gives constants independent of
$\eps_n$.  Theorem~\ref{thm:weighted-compactness} therefore supplies a
subsequence, still indexed by $n$, and functions $(u_0,u_1)$ satisfying
\eqref{eq:compact-u}--\eqref{eq:compact-weighted} and the time regularity and
initial condition in \eqref{eq:macro-time}.  Proposition~\ref{prop:micro-equation}
adds $u_1\in L^2(Q;\mathcal W_\omega)$.

Since $\gamma_0P_{\eps_n}u_{\eps_n}\leq\varphi$ on $\Gamma_1\times(0,T)$ and
this convex constraint is closed in $L^2(0,T;W^{1,q}(\Om))$, the weak limit
satisfies $\gamma_0u_0\leq\varphi$.  Therefore
$U\in\mathbb K_\omega$.

It remains to prove \eqref{eq:two-scale-Minty}.  Start with a comparison from the
core in Lemma~\ref{lem:limit-density}; in particular $v_0(0)=0$ and
$v_1$ is compactly supported in $(x,t)\in Q$.  Set
\begin{equation}
 v_\eps=v_0+\eps v_1^\eps.
\end{equation}
Since $v_1$ vanishes near $\partial\Om$ and near $t=0$, $v_\eps=v_0$ on the
exterior boundary, $v_\eps\in K_\eps(t)$, and $v_\eps(0)=0$.  For every
$i=1,\ldots,N$,
\begin{equation}
 D_iv_\eps=D_{x_i}v_0+(D_{y_i}v_1)^\eps
                  +\eps(D_{x_i}v_1)^\eps,               \label{eq:comparison-gradient}
\end{equation}
so Corollary~\ref{cor:minty-recovery} applies to the comparison state and its
normalized gradient.

The comparison $v_\eps$ belongs to $\mathcal K_\eps^1(T)$ for small
$\eps$, and $v_\eps(0)=0$.  Gradient monotonicity \eqref{ass:A3}, with the state frozen at
$u_\eps$, gives
\begin{equation*}
 \int_{\Omeps}A^\eps[u_\eps]\cdot D(v_\eps-u_\eps)
 \leq\int_{\Omeps}A^\eps(x,t,u_\eps,Dv_\eps)\cdot D(v_\eps-u_\eps).
\end{equation*}
Consequently
\begin{align}
 &\int_0^T\langle B_\eps\partial_tv_\eps,v_\eps-u_\eps\rangle\dd t
 +\int_{\Qeps}\bigl[A_0^\eps[u_\eps](v_\eps-u_\eps)
 +A^\eps(x,t,u_\eps,Dv_\eps)\cdot D(v_\eps-u_\eps)\bigr]\dd x\dd t\nonumber\\
 &\qquad\geq\int_{\Qeps}f(v_\eps-u_\eps)\dd x\dd t
 +\int_0^T\langle\ell_\eps(t),v_\eps-u_\eps\rangle\dd t.
\end{align}
We calculate the limit of the comparison time term.  The chain rule gives
\begin{equation*}
 \partial_tv_\eps=\partial_tv_0+\eps(\partial_tv_1)^\eps
                    +\eps^{-1}(\partial_\tau v_1)^\eps.
\end{equation*}
The contribution of $\partial_tv_0$ is
\begin{equation*}
 \int_{\Qeps}b^\eps\partial_tv_0(v_\eps-u_\eps)\dd x\dd t
 \longrightarrow
 \int_0^T\langle\overline b\,\partial_tv_0,v_0-u_0\rangle\dd t;
\end{equation*}
by strong convergence of $v_\eps-P_\eps u_\eps$ in $L^2(Q)$ and periodic averaging.  The term containing
$\eps\partial_tv_1$ is $O(\eps)$.  For the fast term, expand
\begin{align}
 \frac1\eps\int_{\Qeps}b^\eps(\partial_\tau v_1)^\eps(v_0+\eps v_1^\eps-u_\eps)\dd x\dd t
 &=\frac1\eps\int_{\Qeps}b^\eps v_0(\partial_\tau v_1)^\eps\dd x\dd t
 +\int_{\Qeps}b^\eps v_1^\eps(\partial_\tau v_1)^\eps\dd x\dd t\nonumber\\
 &\quad-\frac1\eps\int_{\Qeps}b^\eps u_\eps(\partial_\tau v_1)^\eps\dd x\dd t.
 \label{eq:minty-fast-three}
\end{align}
Put $\psi=\partial_\tau v_1$, which has zero $b$-mean, and
$\boldsymbol{\Theta}=\mathcal B_{\Ystar}(b\psi)$ as in Lemma~\ref{lem:very-weak}.
The identity \eqref{eq:very-weak-identity} with $v_0$ in place of $P_\eps u_\eps$
shows that the first term tends to
$-\int_Q\int_Z\int_{\Ystar}(D_xv_0\cdot\boldsymbol{\Theta}+v_0\Div_x\boldsymbol{\Theta})=0$,
by integration by parts in $x$.  The second tends, by periodic
averaging, to
\begin{equation*}
 \int_Q\int_Z\int_{\Ystar}bv_1\partial_\tau v_1
 =\frac12\int_Q\int_Z\partial_\tau
       \left(\int_{\Ystar}b|v_1|^2\dd y\right)\dd\tau\dd x\dd t=0.
\end{equation*}
The last term in \eqref{eq:minty-fast-three} tends by
Lemma~\ref{lem:very-weak} to
$-\int_Q\int_Z\int_{\Ystar}bu_1\partial_\tau v_1$.
Since $v_1$ is smooth, $\int_{\Ystar}bu_1\partial_\tau v_1\dd y
=\langle\partial_\tau(bv_1),u_1\rangle$ and
$\int_Z\langle\partial_\tau(bv_1),v_1\rangle\dd\tau=0$; this limit is therefore
\begin{equation*}
 \int_Q\int_Z
       \langle\partial_\tau(bv_1),v_1-u_1\rangle
       \dd\tau\dd x\dd t.
\end{equation*}
The macro and micro contributions together are
$\langle\mathbb T V,V-U\rangle$.

For the zero-order term, \eqref{eq:A0-solution-limit} and
$v_\eps-P_\eps u_\eps\to v_0-u_0$ strongly in $L^2(Q)$ give
\begin{align*}
 &\int_{\Qeps}A_0^\eps[u_\eps](v_\eps-u_\eps)\dd x\dd t\\
 &\quad\longrightarrow\int_Q\int_Z\int_{\Ystar}
       A_0(y,\tau,u_0)(v_0-u_0)\dd y\dd\tau\dd x\dd t.
\end{align*}
For every $i=1,\ldots,N$, Corollary~\ref{cor:minty-recovery} and
\eqref{eq:comparison-gradient} give
\begin{equation*}
 \chi_{\Omeps}\frac{A_i^\eps(x,t,u_\eps,Dv_\eps)}
                         {\sqrt{w_i^\eps}}
 \tto\chi_{\Ystar}\frac{A_i(y,\tau,u_0,\mathbb D V)}
                         {\sqrt{w_i(y)}}
\end{equation*}
strongly in $L^2(Q)$.  On the other hand,
\begin{equation*}
 \chi_{\Omeps}\sqrt{w_i^\eps}D_i(v_\eps-u_\eps)
 \twto\chi_{\Ystar}\sqrt{w_i}
          \bigl((\mathbb D V)_i-(\mathbb D U)_i\bigr)
\end{equation*}
weakly in $L^2(Q)$, by \eqref{eq:compact-weighted} and
\eqref{eq:comparison-gradient}.  Summing the strong--weak product limits over
$i=1,\ldots,N$ gives the principal term in \eqref{eq:two-scale-Minty}.

Since $v_\eps-P_\eps u_\eps\to v_0-u_0$ strongly in $L^2(Q)$ and the material
characteristic has mean $\theta$, the volume source converges to
$\theta\int_Qf(v_0-u_0)$.  Theorem~\ref{thm:surface-limit}, applied to
$z_\eps=v_\eps-u_\eps$, gives
\begin{equation*}
 \int_0^T\langle\ell_\eps(t),v_\eps-u_\eps\rangle\dd t
 \longrightarrow \int_Q\int_Z\langle g(\tau),\gamma_0^y(v_1-u_1)\rangle_S
 \dd\tau\dd x\dd t+\int_Qm_g\cdot D_x(v_0-u_0)\dd x\dd t.
\end{equation*}
This proves \eqref{eq:two-scale-Minty} for the core comparisons; the graph-norm
density of Lemma~\ref{lem:limit-density} and the continuity of the Nemytskii
operators extend it to all admissible $V$.
\end{proof}

For $z\in\mathcal K_0$, Lemma~\ref{lem:exp-regularization}, with
$V=H^1_{\Gamma_2}(\Om)$, $C=C_0$, $\sigma=T$ and $B$ the pivot embedding, gives
$z_h\in H^1(0,T;H^1_{\Gamma_2}(\Om))\cap\mathcal D_0$ with
\begin{equation}
 z_h\longrightarrow z\ \text{ in }X_\Gamma\ (h\to0),\qquad
 \langle\partial_tz_h,z_h-z\rangle_{X_\Gamma}\leq0.        \label{eq:exp-reg-props}
\end{equation}

The next proposition identifies the fast-time cell equation at the limit
state under the base hypotheses alone.

\begin{proposition}\label{prop:frozen-cell}
Under the hypotheses of Theorem~\ref{thm:main}, for a.e. $(x,t)\in Q$,
\begin{equation}
 \int_Z\langle\partial_\tau(bu_1),\psi\rangle\dd\tau
 +\int_Z\int_{\Ystar}A(y,\tau,u_0,\mathbb D U)\cdot D_y\psi\dd y\dd\tau
 =\int_Z\langle g(\tau),\gamma_0^y\psi\rangle_S\dd\tau
 \qquad\forall\,\psi\in L^2_\#(Z;V_\omega).                          \label{eq:cell-at-solution}
\end{equation}
\end{proposition}

\begin{proof}
Put $\tilde u=u_0-\Psi$.  Since $\gamma_0\Psi(t)=\varphi(t)-\varphi(0)$ on
$\Gamma_1$ and $\gamma_0\Psi=0$ on $\Gamma_2$, the obstacle for $u_0$ gives
$\tilde u\in\mathcal K_0$.  Fix $\psi\in L^2(Q;\mathcal W_\omega)$, $\delta>0$ and
$h>0$, and use in \eqref{eq:two-scale-Minty} the comparison
\[
 v_0=\Psi+\tilde u_h,\qquad v_1=u_1+\delta\psi,
\]
with $\tilde u_h$ from Lemma~\ref{lem:exp-regularization}.  Then
$v_0\in H^1(0,T;H^1_{\Gamma_2}(\Om))$, $v_0(0)=0$ and
$\gamma_0v_0\leq\varphi(t)-\varphi(0)+\varphi(0)=\varphi(t)$ on $\Gamma_1$;
moreover $v_1\in L^2(Q;\mathcal W_\omega)$.  Hence
$V=(v_0,v_1)\in\mathbb K_\omega\cap\mathbb E_\omega$, and $v_0-u_0=\tilde u_h-\tilde u$.

We examine every term of \eqref{eq:two-scale-Minty}.  The macroscopic time
term equals
\[
 \overline b\,\langle\partial_t\Psi,\tilde u_h-\tilde u\rangle_{X_\Gamma}
 +\overline b\,\langle\partial_t\tilde u_h,\tilde u_h-\tilde u\rangle_{X_\Gamma}.
\]
The first part tends to zero by \eqref{eq:exp-reg-props}, since
$\partial_t\Psi\in L^2(Q)$; the second is non-positive by
\eqref{eq:exp-reg-props}, so it may be discarded.  The fast-time term
equals
$\delta\int_Q\int_Z\langle\partial_\tau(bu_1),\psi\rangle
+\delta^2\int_Q\int_Z\langle\partial_\tau(b\psi),\psi\rangle$, and the second
part vanishes by \eqref{eq:fast-ibp}.  In $\langle\mathbb L,V-U\rangle$, the
source term and the term containing $m_g$ tend to zero by
\eqref{eq:exp-reg-props}, while the microscopic surface term equals
$\delta\int_Q\int_Z\langle g,\gamma_0^y\psi\rangle_S$.  The zero-order term
tends to zero.  Finally,
$\mathbb DV-\mathbb DU=D_x(\tilde u_h-\tilde u)+\delta D_y\psi$, and the normalized
fluxes $A_i(y,\tau,u_0,\mathbb DV)/\sqrt{w_i}$ converge strongly in
$L^2(Q\times\Ystar\times Z)$ to
$A_i(y,\tau,u_0,\mathbb DU+\delta D_y\psi)/\sqrt{w_i}$ by the Carath\'eodory
property and \eqref{ass:A1b}, through the truncation argument of
Proposition~\ref{prop:corrected31}.  Letting $h\to0$ and dividing by
$\delta$ yields
\begin{align*}
 &\int_Q\int_Z\langle\partial_\tau(bu_1),\psi\rangle\dd\tau\dd x\dd t
 +\int_Q\int_Z\int_{\Ystar}A(y,\tau,u_0,\mathbb DU+\delta D_y\psi)\cdot
       D_y\psi\dd y\dd\tau\dd x\dd t\\
 &\qquad-\int_Q\int_Z\langle g(\tau),\gamma_0^y\psi\rangle_S\dd\tau\dd x\dd t
 \geq0.
\end{align*}
Letting $\delta\downarrow0$ by the same continuity, and replacing $\psi$ by
$-\psi$, gives the integrated equality.  Its left-hand side is continuous in
$\psi$ for the norm of $L^2(Q\times Z;V_\omega)$, because
$\partial_\tau(bu_1)\in L^2(Q\times Z;V_\omega')$; by Lemma~\ref{lem:graph-core},
$L^2(Q;\mathcal W_\omega)$ is dense in that space.  Choosing
$\psi=\alpha(x,t)\zeta$, with $\alpha\in L^\infty(Q)$ and $\zeta$ in a countable
dense subset of $L^2_\#(Z;V_\omega)$, localizes the identity and proves
\eqref{eq:cell-at-solution}.
\end{proof}

\subsection{The cell problem under \texorpdfstring{$(H_{10})$}{(H10)}}\label{sec:decoupling}

In the rest of this section we assume \eqref{ass:Adu}.  Its joint
state--gradient form yields uniqueness of the cell corrector, continuity of
the cell law and strong monotonicity of the effective operator.

For $(r,\xi)\in\R\times\R^N$, the cell problem consists in finding
$\pi(r,\xi)\in\mathcal W_\omega$ such that
\begin{equation}
 \int_Z\langle\partial_\tau(b\pi),\psi\rangle\dd\tau
 +\int_Z\int_{\Ystar}A(y,\tau,r,\xi+D_y\pi)\cdot D_y\psi\dd y\dd\tau
 =\int_Z\langle g(\tau),\gamma_0^y\psi\rangle_S\dd\tau
 \qquad\forall\,\psi\in L^2_\#(Z;V_\omega),                  \label{eq:cell-parametric}
\end{equation}
and the effective zero-order term and flux are
\begin{align}
 A_0^{\rm hom}(r)&=\int_Z\int_{\Ystar}
 A_0(y,\tau,r)\dd y\dd\tau,                              \label{eq:A0hom}\\
 A^{\rm hom}(r,\xi)&=\int_Z\int_{\Ystar}
 A(y,\tau,r,\xi+D_y\pi(r,\xi))\dd y\dd\tau.             \label{eq:Ahom}
\end{align}
Note that $A^{\rm hom}$ depends on $g$ through $\pi$, and that the
effective flux is $A^{\rm hom}(r,\xi)-m_g$.

\begin{proposition}\label{prop:cell}
Assume \eqref{ass:Adu}.
\begin{enumerate}
\item[(a)] For every $(r,\xi)$, problem \eqref{eq:cell-parametric} has exactly
one solution $\pi(r,\xi)\in\mathcal W_\omega$, and there is $C>0$, independent
of $(r,\xi)$, such that
\begin{align}
 \|\pi(r,\xi)\|_{L^2_\#(Z;V_\omega)}
 +\Bigl(\sum_{i=1}^N\int_Z\int_{\Ystar}w_i|\xi_i+D_{y_i}\pi(r,\xi)|^2\dd y\dd\tau\Bigr)^{1/2}
 &\leq C(1+|r|+|\xi|), \label{eq:cell-solution-growth}\\
 \|\partial_\tau(b\pi(r,\xi))\|_{L^2_\#(Z;V_\omega')}
 &\leq C(1+|r|+|\xi|).                                      \label{eq:cell-time-growth}
\end{align}
\item[(b)] The map $(r,\xi)\mapsto\pi(r,\xi)$ is continuous for the graph norm
of $\mathcal W_\omega$, and $A_0^{\rm hom}$ and $A^{\rm hom}$ are continuous
and satisfy
\begin{equation}
 |A_0^{\rm hom}(r)|\leq C(1+|r|),\qquad
 |A^{\rm hom}(r,\xi)|\leq C(1+|r|+|\xi|).
 \label{eq:effective-growth}
\end{equation}
\item[(c)] There is $c_e>0$ such that, for all
$(r_j,\xi_j)\in\R\times\R^N$, $j=1,2$,
\begin{equation}
 [A_0^{\rm hom}(r_1)-A_0^{\rm hom}(r_2)](r_1-r_2)
 +[A^{\rm hom}(r_1,\xi_1)-A^{\rm hom}(r_2,\xi_2)]\cdot(\xi_1-\xi_2)
 \geq c_e|\xi_1-\xi_2|^2.                       \label{eq:effective-monotonicity}
\end{equation}
\end{enumerate}
\end{proposition}

\begin{proof}
(a) Let $\mathcal H=L^2_\#(Z;V_\omega)$ and
$\langle\mathcal C_{r,\xi}z,\psi\rangle
=\int_Z\int_{\Ystar}A(y,\tau,r,\xi+D_yz)\cdot D_y\psi$.  By \eqref{ass:A1b},
$\mathcal C_{r,\xi}:\mathcal H\to\mathcal H'$ is bounded, and it is
hemicontinuous.  Taking $r=\bar r$ in \eqref{ass:Adu},
\begin{equation}
 \langle\mathcal C_{r,\xi}z_1-\mathcal C_{r,\xi}z_2,z_1-z_2\rangle
 \geq c_u\sum_{i=1}^N\int_Z\int_{\Ystar}
       w_i|D_{y_i}(z_1-z_2)|^2\dd y\dd\tau.             \label{eq:cell-strong-monotone}
\end{equation}
With $\eta=\xi+D_yz$, \eqref{ass:A1a}--\eqref{ass:A2} and weighted
Cauchy--Schwarz give
\begin{equation}
 \langle\mathcal C_{r,\xi}z,z\rangle
 \geq\frac{c_0}{2}\sum_{i=1}^N\int_Z\int_{\Ystar}w_i|\xi_i+D_{y_i}z|^2\dd y\dd\tau
 -C\bigl(1+|r|^2+|\xi|^2\bigr), \label{eq:cell-coercive-expanded}
\end{equation}
and by \eqref{eq:cell-affine} the first term controls $\|z\|_{\mathcal H}^2$.
The surface functional
$\langle\mathcal F,\psi\rangle=\int_Z\langle g(\tau),\gamma_0^y\psi(\tau)\rangle_S$
belongs to $\mathcal H'$ by \eqref{eq:lifting-functional}.  By
Lemma~\ref{lem:periodic-degenerate-derivative}, $\mathcal Lz=\partial_\tau(bz)$
is maximal monotone with $\langle\mathcal Lz,z\rangle=0$.  The range theorem
\cite[Theorem~2.2]{Le2011}, applied with the maximal monotone operator
$\mathcal L$ (whose domain contains $0$), the bounded pseudomonotone operator
$\mathcal C_{r,\xi}$, $u_0=0$ and $R$ large, condition (2.1) there following
from \eqref{eq:cell-coercive-expanded} and \eqref{eq:cell-affine}, gives
$z\in\mathcal W_\omega$ with $\mathcal Lz+\mathcal C_{r,\xi}z=\mathcal F$.
Uniqueness follows from \eqref{eq:fast-ibp}, \eqref{eq:cell-strong-monotone}
and \eqref{eq:cell-poincare}.  Testing the equation by $z$ and using
\eqref{eq:cell-coercive-expanded}, \eqref{eq:lifting-functional} and
\eqref{eq:cell-affine} gives \eqref{eq:cell-solution-growth}; the equation
itself and \eqref{ass:A1b} then give \eqref{eq:cell-time-growth}.

(b) Let $(r_n,\xi_n)\to(r,\xi)$, $\pi_n=\pi(r_n,\xi_n)$, $\pi=\pi(r,\xi)$,
$\eta_n=\xi_n+D_y\pi_n$, $\eta=\xi+D_y\pi$.  Subtracting the cell equations,
testing by $\pi_n-\pi$ and using \eqref{eq:fast-ibp} and \eqref{ass:Adu} give
\[
 c_u\sum_{i=1}^N\int_Z\int_{\Ystar}w_i|\eta_{n,i}-\eta_i|^2
 \leq\int_Z\int_{\Ystar}\bigl([A_0(r_n)-A_0(r)](r_n-r)
 +[A(r_n,\eta_n)-A(r,\eta)]\cdot(\xi_n-\xi)\bigr)\longrightarrow0,
\]
by \eqref{ass:A1a}--\eqref{ass:A1b} and the bound
\eqref{eq:cell-solution-growth}.  By \eqref{eq:cell-affine},
$\pi_n\to\pi$ in $L^2_\#(Z;V_\omega)$.  To justify convergence of the fluxes,
put $E=\Ystar\times Z$ and $z_{n,j}=\sqrt{w_j}\eta_{n,j}$,
$z_j=\sqrt{w_j}\eta_j$.  Then $z_n\to z$ in $L^2(E)^N$.  For a.e.
$(y,\tau)$ define
\[
 \mathscr B_i(y,\tau,r,z)=\frac{1}{\sqrt{w_i(y)}}
 A_i\Bigl(y,\tau,r,(z_j/\sqrt{w_j(y)})_{j=1}^N\Bigr).
\]
It is Carath\'eodory and, by \eqref{ass:A1b},
$|\mathscr B_i(y,\tau,r,z)|\leq c_1(k_1+|r|+\sum_{j=1}^N|z_j|)$.
Every subsequence of $z_n$ has a further subsequence converging a.e.;
on that subsequence $\mathscr B_i(\cdot,r_n,z_n)$ converges a.e., while its
squares are uniformly integrable because $z_n\to z$ in $L^2(E)^N$ and
$k_1\in L^2(E)$.  Vitali's theorem thus gives
\[
 \frac{A_i(\cdot,r_n,\eta_n)}{\sqrt{w_i}}
 \longrightarrow\frac{A_i(\cdot,r,\eta)}{\sqrt{w_i}}
 \quad\text{in }L^2(E),\qquad i=1,\ldots,N.
\]
The same argument applies to $A_0(\cdot,r_n)$.  Consequently
$A_0^{\rm hom}$ and $A^{\rm hom}$ are continuous, and the cell equation
gives convergence of $\partial_\tau(b\pi_n)$ in
$L^2_\#(Z;V_\omega')$.  The growth bounds
\eqref{eq:effective-growth} follow from \eqref{ass:A1a}--\eqref{ass:A1b} and
\eqref{eq:cell-solution-growth}.

(c) Put $\pi_j=\pi(r_j,\xi_j)$ and
$\eta_j=\xi_j+D_y\pi_j$, $j=1,2$.  Subtracting the two cell equations and
testing by $\pi_1-\pi_2$, the periodic time term vanishes by
\eqref{eq:fast-ibp} and the common surface functional cancels, so
\begin{equation}
 \int_Z\int_{\Ystar}[A(y,\tau,r_1,\eta_1)-A(y,\tau,r_2,\eta_2)]
       \cdot D_y(\pi_1-\pi_2)\dd y\dd\tau=0.             \label{eq:cell-difference}
\end{equation}
Using \eqref{eq:cell-difference}, the definitions
\eqref{eq:A0hom}--\eqref{eq:Ahom}, and then \eqref{ass:Adu}, we obtain
\[
 \begin{aligned}
 &[A_0^{\rm hom}(r_1)-A_0^{\rm hom}(r_2)](r_1-r_2)
 +[A^{\rm hom}(r_1,\xi_1)-A^{\rm hom}(r_2,\xi_2)]\cdot(\xi_1-\xi_2)\\
 &\qquad\geq c_u\int_Z\int_{\Ystar}\sum_{i=1}^Nw_i
       |\eta_{1i}-\eta_{2i}|^2\dd y\dd\tau .
 \end{aligned}
\]
For each $\tau$, apply the affine cell inequality \eqref{eq:cell-affine} to
$(\xi_1-\xi_2,\pi_1(\tau)-\pi_2(\tau))$.  Integration over $Z$ controls
$|\xi_1-\xi_2|^2$ by the positive integral above and proves
\eqref{eq:effective-monotonicity}.
\end{proof}

\subsection{The macroscopic problem, uniqueness and full convergence}\label{sec:macro}

The decoupled macroscopic problem is identified in a Lions-type integrated
comparison formulation \cite[Chapter~2]{Lions1969}, in which the time
derivative falls only on the comparison.  The comparisons vanish at $t=0$,
and the initial condition $u_0(0)=0$ is carried separately in $H^{-1}(\Om)$;
for a weak formulation with general comparisons and a quadratic initial term,
see \cite[Section~4]{MigorskiDudek2025}.  Recall the
notation $X_\Gamma$, $C_0$, $\mathcal K_0$, $\mathcal D_0$ and
$\langle\cdot,\cdot\rangle_{X_\Gamma}$ introduced before
Lemma~\ref{lem:limit-density}, and denote by $\mathcal J$ the pivot map,
$\langle\mathcal J\eta,\zeta\rangle=\int_\Om\eta\zeta\dd x$ for $\eta\in L^2(\Om)$,
together with its pointwise realization in time.  For $r\in L^2(Q)$ and a.e.
$t\in(0,T)$ define $\widetilde Q_r(t):H^1_{\Gamma_2}(\Om)\to(H^1_{\Gamma_2}(\Om))'$ by
\begin{equation}
 \langle\widetilde Q_r(t)z,\zeta\rangle
 =\int_\Om A_0^{\rm hom}(r(t))\zeta\dd x
 +\int_\Om\bigl[A^{\rm hom}(r(t),D_xz+D_x\Psi(t))-m_g\bigr]
        \cdot D_x\zeta\dd x,
\end{equation}
let $\widetilde{\mathcal Q}_r:X_\Gamma\to X_\Gamma'$ be its Nemytskii
realization, and put
\[
 \widetilde{\mathcal Q}(z)=\widetilde{\mathcal Q}_{z+\Psi}(z),\qquad
 \widetilde F=\theta f-\overline b\,\mathcal J\partial_t\Psi\in X_\Gamma'.
\]
Thus $\widetilde{\mathcal Q}_r$ freezes the state at $r$, while
$\widetilde{\mathcal Q}$ evaluates it at $z+\Psi$.

\begin{lemma}\label{lem:effective-minty}
Assume \eqref{ass:Adu}, let $U=(u_0,u_1)$ be any subsequential limit furnished
by Theorem~\ref{thm:main}, and put $\tilde u=u_0-\Psi$.  Then
$\tilde u\in\mathcal K_0$ and
\begin{equation}
 \langle\overline b\,\partial_tz+\widetilde{\mathcal Q}_{u_0}(z)-\widetilde F,
       z-\tilde u\rangle_{X_\Gamma}\geq0
 \qquad\forall\,z\in\mathcal D_0.                                  \label{eq:translated-Minty}
\end{equation}
Moreover, $\widetilde{\mathcal Q}_{u_0}$ and $\widetilde{\mathcal Q}$ are
continuous from $X_\Gamma$ to $X_\Gamma'$ and bounded on bounded sets,
$\widetilde{\mathcal Q}_{u_0}(\tilde u)=\widetilde{\mathcal Q}(\tilde u)$, and
there is $m_e>0$ such that
\begin{equation}
 \langle\widetilde{\mathcal Q}(z_1)-\widetilde{\mathcal Q}(z_2),
       z_1-z_2\rangle_{X_\Gamma}
 \geq m_e\|z_1-z_2\|_{X_\Gamma}^2
 \qquad\forall\,z_1,z_2\in X_\Gamma.                               \label{eq:Q-monotone}
\end{equation}
\end{lemma}

\begin{proof}
Since $\gamma_0\Psi(t)=\varphi(t)-\varphi(0)$ on $\Gamma_1$ and
$\gamma_0\Psi=0$ on $\Gamma_2$, the obstacle for $u_0$ gives
$\gamma_0\tilde u\leq\varphi(0)$ on $\Gamma_1$; hence $\tilde u\in\mathcal K_0$.

Let $z\in\mathcal D_0$, put $v=z+\Psi$ and insert in \eqref{eq:two-scale-Minty}
the cell solution at the frozen state,
\[
 v_1(x,t)=\pi\bigl(u_0(x,t),D_xv(x,t)\bigr).
\]
Since $(x,t)\mapsto(u_0(x,t),D_xv(x,t))$ is strongly measurable,
Proposition~\ref{prop:cell}(b) and the separability of $L^2_\#(Z;V_\omega)$
and $L^2_\#(Z;V_\omega')$ show that $v_1$ and $\partial_\tau(bv_1)$ are
strongly measurable, and Proposition~\ref{prop:cell}(a) gives
\[
 \|v_1\|_{L^2(Q;L^2_\#(Z;V_\omega))}
 +\|\partial_\tau(bv_1)\|_{L^2(Q\times Z;V_\omega')}
 \leq C\bigl(1+\|u_0\|_{L^2(Q)}+\|v\|_{L^2(0,T;H^1(\Om))}\bigr).
\]
Since $\partial_tv\in X_\Gamma'$, $v(0)=0$ and $\gamma_0v\leq\varphi$ on
$\Gamma_1\times(0,T)$, the pair $V=(v,v_1)$ belongs to
$\mathbb K_\omega\cap\mathbb E_\omega$ and is admissible in
\eqref{eq:two-scale-Minty}.  By construction
$A(y,\tau,u_0,\mathbb DV)=A(y,\tau,u_0,D_xv+D_yv_1)$ is the flux of the cell
problem \eqref{eq:cell-parametric} with $(r,\xi)=(u_0,D_xv)$.  Testing this
cell problem by $v_1-u_1\in L^2(Q;L^2_\#(Z;V_\omega))$ and integrating over
$Q$ cancels, in \eqref{eq:two-scale-Minty}, the fast-time term, the
microscopic surface term and the $D_y(v_1-u_1)$ part of the flux.  By
\eqref{eq:A0hom}--\eqref{eq:Ahom}, what remains is
\begin{align*}
 &\int_0^T\langle\overline b\,\partial_tv,v-u_0\rangle\dd t
 +\int_QA_0^{\rm hom}(u_0)(v-u_0)\dd x\dd t
 +\int_Q\bigl(A^{\rm hom}(u_0,D_xv)-m_g\bigr)\cdot D_x(v-u_0)\dd x\dd t\\
 &\qquad\geq\theta\int_Qf(v-u_0)\dd x\dd t.
\end{align*}
Since $v-u_0=z-\tilde u$ and $\partial_tv=\partial_tz+\partial_t\Psi$, this is
exactly \eqref{eq:translated-Minty}.

By Proposition~\ref{prop:cell}(b), $A_0^{\rm hom}$ and $A^{\rm hom}$ are
continuous with linear growth.  If $z_n\to z$ in $X_\Gamma$, then, along a
subsequence, $(z_n,D_xz_n)\to(z,D_xz)$ a.e. in $Q$ with an $L^2(Q)$ majorant;
continuity, \eqref{eq:effective-growth} and dominated convergence give
convergence in $L^2(Q)$ of the compositions defining
$\widetilde{\mathcal Q}_{u_0}(z_n)$ and $\widetilde{\mathcal Q}(z_n)$, and the
subsequence principle yields continuity of both operators.  Boundedness on
bounded sets follows from \eqref{eq:effective-growth} and
$\Psi\in H^1(0,T;H^1(\Om))$, and
$\widetilde{\mathcal Q}_{u_0}(\tilde u)=\widetilde{\mathcal Q}(\tilde u)$
because $\tilde u+\Psi=u_0$.  Finally, \eqref{eq:effective-monotonicity}
applied at the points $(z_j+\Psi,D_xz_j+D_x\Psi)$, $j=1,2$, gives
\[
 \langle\widetilde{\mathcal Q}(z_1)-\widetilde{\mathcal Q}(z_2),
       z_1-z_2\rangle_{X_\Gamma}
 \geq c_e\|D_x(z_1-z_2)\|_{L^2(Q)^N}^2,
\]
the constant vector $m_g$ cancelling.  Since $|\Gamma_2|_{N-1}>0$, the
Poincar\'e inequality on $H^1_{\Gamma_2}(\Om)$ yields \eqref{eq:Q-monotone}.
\end{proof}

The comparison form \eqref{eq:translated-Minty} evaluates the effective flux at
the gradient of the test function.  The next lemma converts this comparison
inequality into the integrated weak formulation, without additional time
regularity of $\tilde u$.

\begin{lemma}
\label{lem:minty-to-lions}
Assume \eqref{ass:Adu} and let $\tilde u$ be as in Lemma~\ref{lem:effective-minty}.
Then
\begin{equation}
 \langle\overline b\,\partial_tw+\widetilde{\mathcal Q}(\tilde u)-\widetilde F,
       w-\tilde u\rangle_{X_\Gamma}\geq0
 \qquad\forall\,w\in\mathcal D_0.                                  \label{eq:lions-form}
\end{equation}
\end{lemma}

\begin{proof}
Fix $w\in\mathcal D_0$, $\lambda\in(0,1)$ and $h>0$, and let $\tilde u_h$
be the regularization of Lemma~\ref{lem:exp-regularization}.  The convex
combination
\[
 z=(1-\lambda)\tilde u_h+\lambda w
\]
belongs to $\mathcal D_0$, and
$z-\tilde u=(1-\lambda)(\tilde u_h-\tilde u)+\lambda(w-\tilde u)$.  Expanding the time term,
\begin{align*}
 \langle\partial_tz,z-\tilde u\rangle_{X_\Gamma}
 &=(1-\lambda)^2\langle\partial_t\tilde u_h,\tilde u_h-\tilde u\rangle_{X_\Gamma}
 +(1-\lambda)\lambda\langle\partial_t\tilde u_h,w-\tilde u\rangle_{X_\Gamma}\\
 &\quad+\lambda(1-\lambda)\langle\partial_tw,\tilde u_h-\tilde u\rangle_{X_\Gamma}
 +\lambda^2\langle\partial_tw,w-\tilde u\rangle_{X_\Gamma}.
\end{align*}
The first term is non-positive by \eqref{eq:exp-reg-props}.  The third tends
to zero as $h\to0$ by \eqref{eq:exp-reg-props}.  For the second, use
\eqref{eq:exp-reg-props} once more and then the product rule in the Gelfand triple,
which applies because $\tilde u_h\in H^1(0,T;H^1_{\Gamma_2}(\Om))$, $w\in\mathcal D_0$ and
$\tilde u_h(0)=w(0)=0$:
\begin{align*}
 \langle\partial_t\tilde u_h,w-\tilde u\rangle_{X_\Gamma}
 &\leq\langle\partial_t\tilde u_h,w-\tilde u_h\rangle_{X_\Gamma}
 =(\tilde u_h(T),w(T))_{L^2(\Om)}-\tfrac12\|\tilde u_h(T)\|_{L^2(\Om)}^2
 -\langle\partial_tw,\tilde u_h\rangle_{X_\Gamma}\\
 &\leq\tfrac12\|w(T)\|_{L^2(\Om)}^2-\langle\partial_tw,\tilde u_h\rangle_{X_\Gamma}.
\end{align*}
Since
$\tfrac12\|w(T)\|_{L^2(\Om)}^2=\langle\partial_tw,w\rangle_{X_\Gamma}$, letting
$h\to0$ gives
\[
 \limsup_{h\downarrow0}\langle\overline b\,\partial_tz,z-\tilde u\rangle_{X_\Gamma}
 \leq\overline b\bigl[(1-\lambda)\lambda+\lambda^2\bigr]
       \langle\partial_tw,w-\tilde u\rangle_{X_\Gamma}
 =\lambda\langle\overline b\,\partial_tw,w-\tilde u\rangle_{X_\Gamma}.
\]
By continuity of $\widetilde{\mathcal Q}_{u_0}$ (Lemma~\ref{lem:effective-minty}),
\[
 \langle\widetilde{\mathcal Q}_{u_0}(z)-\widetilde F,z-\tilde u\rangle_{X_\Gamma}
 \longrightarrow
 \lambda\langle\widetilde{\mathcal Q}_{u_0}(\tilde u+\lambda(w-\tilde u))-\widetilde F,
       w-\tilde u\rangle_{X_\Gamma}
 \qquad(h\to0).
\]
Combining these limits with \eqref{eq:translated-Minty} gives
\[
 0\leq\lambda\langle\overline b\,\partial_tw,w-\tilde u\rangle_{X_\Gamma}
 +\lambda\langle\widetilde{\mathcal Q}_{u_0}(\tilde u+\lambda(w-\tilde u))-\widetilde F,
       w-\tilde u\rangle_{X_\Gamma}.
\]
Divide by $\lambda$ and let $\lambda\downarrow0$, using once more the
continuity of $\widetilde{\mathcal Q}_{u_0}$.  Since
$\widetilde{\mathcal Q}_{u_0}(\tilde u)=\widetilde{\mathcal Q}(\tilde u)$, this
is \eqref{eq:lions-form}.
\end{proof}

\begin{theorem}\label{thm:decoupled}
Assume \((H_1)\)--\((H_8)\) and \((H_{10})\), and let $U=(u_0,u_1)$ be any
subsequential limit furnished by Theorem~\ref{thm:main}.  Then:
\begin{enumerate}
\item[(i)] the corrector is given by the cell law,
\begin{equation}
 u_1(x,t,\cdot,\cdot)=\pi(u_0(x,t),D_xu_0(x,t))
 \quad\text{for a.e. }(x,t)\in Q;                              \label{eq:corrector-decoupled}
\end{equation}
\item[(ii)] $u_0\in L^2(0,T;H^1_{\Gamma_2}(\Om))\cap H^1(0,T;H^{-1}(\Om))$,
$u_0(0)=0$ in $H^{-1}(\Om)$, $\gamma_0u_0\leq\varphi$ on $\Gamma_1\times(0,T)$,
and
\begin{align}
 &\int_0^T\langle\overline b\,\partial_tv,v-u_0\rangle\dd t
 +\int_QA_0^{\rm hom}(u_0)(v-u_0)\dd x\dd t\nonumber\\
 &\quad+\int_Q\bigl(A^{\rm hom}(u_0,D_xu_0)-m_g\bigr)\cdot D_x(v-u_0)\dd x\dd t
 \geq\theta\int_Qf(v-u_0)\dd x\dd t\nonumber\\
 &\hspace{60mm}\forall\,v\in\mathcal V_{\rm ad},                \label{eq:macroVI}
\end{align}
where $\mathcal V_{\rm ad}$ is the set of $v\in L^2(0,T;H^1_{\Gamma_2}(\Om))$ such
that $\partial_tv\in L^2(0,T;(H^1_{\Gamma_2}(\Om))')$, $v(0)=0$ and
$\gamma_0v\leq\varphi$ on $\Gamma_1\times(0,T)$;
\item[(iii)] problem \eqref{eq:macroVI} has exactly one solution in the class
$\{u\in L^2(0,T;H^1_{\Gamma_2}(\Om)):\gamma_0u\leq\varphi
\text{ on }\Gamma_1\times(0,T)\}$.
\end{enumerate}
Here the time derivative of the comparison is understood in the Gelfand triple
$H^1_{\Gamma_2}(\Om)\subset L^2(\Om)\subset(H^1_{\Gamma_2}(\Om))'$.
\end{theorem}

\begin{proof}
(i) By Proposition~\ref{prop:frozen-cell}, $u_1(x,t)$ solves, for a.e.
$(x,t)$, the cell problem \eqref{eq:cell-parametric} with parameters
$(r,\xi)=(u_0(x,t),D_xu_0(x,t))$, and $u_1(x,t)\in\mathcal W_\omega$ because
$u_1\in L^2(Q;\mathcal W_\omega)$.  Proposition~\ref{prop:cell}(a) gives
\eqref{eq:corrector-decoupled}.

(ii) The spatial and time regularity, the initial condition and the obstacle
are part of Theorems~\ref{thm:weighted-compactness} and \ref{thm:main}.  By
Lemmas~\ref{lem:effective-minty} and \ref{lem:minty-to-lions}, $\tilde u=u_0-\Psi$
satisfies \eqref{eq:lions-form}.  Given an admissible $v$, put $w=v-\Psi$.
Since $\Psi(0)=0$ and $\gamma_0\Psi=\varphi-\varphi(0)$ on $\Gamma_1$,
$w\in\mathcal D_0$; conversely, every $w\in\mathcal D_0$ arises in this way.
Since $w-\tilde u=v-u_0$, $\partial_tw=\partial_tv-\partial_t\Psi$,
$\widetilde F=\theta f-\overline b\,\mathcal J\partial_t\Psi$ and
$\widetilde{\mathcal Q}(\tilde u)$ is the operator in \eqref{eq:macroVI} evaluated at
$u_0$, the terms containing $\partial_t\Psi$ cancel and \eqref{eq:lions-form}
is exactly \eqref{eq:macroVI}.

(iii) Existence is (ii), applied to the microscopic solutions supplied by
Theorem~\ref{thm:existence-eps}, since \((H_{10})\) implies \((H_9)\).  If $u^1$ and $u^2$ are two solutions, the
translation of (ii) shows that $\tilde u^j=u^j-\Psi\in\mathcal K_0$ solve
\eqref{eq:lions-form}.  Since every $w\in H^1(0,T;H^1_{\Gamma_2}(\Om))$ with
$w(t)\in C_0$ and $w(0)=0$ belongs to $\mathcal D_0$, they satisfy
\eqref{eq:abstract-VI} with $V=H^1_{\Gamma_2}(\Om)$, $C=C_0$,
$B=\overline b\,\mathcal J$, $\Psi=0$, $\mathfrak G=\widetilde F$ and
$\mathcal A^j=\widetilde{\mathcal Q}(\tilde u^j)$.
Lemma~\ref{lem:abstract-uniqueness} and \eqref{eq:Q-monotone} give
$m_e\|\tilde u^1-\tilde u^2\|_{X_\Gamma}^2\leq0$.
\end{proof}

\begin{remark}[Strong form and the role of $m_g$]\label{rem:strong-form}
Since $m_g$ is constant,
$\int_\Om m_g\cdot D_xw\dd x=\int_{\Gamma_1}(m_g\cdot\nu)\gamma_0w\dd\sigma$
for $w\in H^1_{\Gamma_2}(\Om)$.  Hence, formally, \eqref{eq:macroVI} is the
weak form of
\[
 \overline b\,\partial_tu-\Div A^{\rm hom}(u,Du)+A_0^{\rm hom}(u)=\theta f
 \ \text{in }Q,\qquad u=0\ \text{on }\Gamma_2\times(0,T),\qquad u(0)=0,
\]
\[
 u\leq\varphi,\quad \nu\cdot(A^{\rm hom}(u,Du)-m_g)\leq0,\quad
 (u-\varphi)\,\nu\cdot(A^{\rm hom}(u,Du)-m_g)=0
 \ \text{on }\Gamma_1\times(0,T).
\]
Thus the direct contribution of the first moment of $g$ is the boundary
term $m_g\cdot\nu$ in the unilateral condition on $\Gamma_1$; this term need
not be constant.

In the linear case $A(y,\tau,r,\xi)=W(y)\xi$, $W=\mathrm{diag}(w_1,\ldots,w_N)$,
with $A_0$ nondecreasing and satisfying \eqref{ass:A1a} (so that
\eqref{ass:A2} holds, and \((H_{10})\) holds with $c_u=1$), the constant can be
computed.  Let $\kappa_i\in V_\omega$
be the stationary correctors,
$\int_{\Ystar}W(e_i+D_y\kappa_i)\cdot D_y\psi\dd y=0$ $\forall\,\psi\in V_\omega$,
and $A^{\rm lin}\xi=\int_{\Ystar}W(\xi+\sum_{j=1}^N\xi_jD_y\kappa_j)\dd y$.  Since
$\partial_\tau(b\kappa_j)=0$, uniqueness in Proposition~\ref{prop:cell}(a) gives
$\pi(r,\xi)=\sum_{j=1}^N\xi_j\kappa_j+\pi_g$, where $\pi_g$ solves
\eqref{eq:cell-parametric} with $\xi=0$, so that
$A^{\rm hom}(\xi)=A^{\rm lin}\xi+c_g$ with
$c_g=\int_Z\int_{\Ystar}WD_y\pi_g$.  Since $W$ is symmetric,
\[
 (c_g)_i=\int_Z\int_{\Ystar}D_y\pi_g\cdot W(e_i+D_y\kappa_i)
 -\int_Z\int_{\Ystar}WD_y\pi_g\cdot D_y\kappa_i .
\]
The first integral vanishes by the corrector equation tested with
$\pi_g(\tau)$.  The second is the equation of $\pi_g$ tested with
$\kappa_i$; its time term vanishes by \eqref{eq:fast-ibp} since $\kappa_i$ does not
depend on $\tau$.  Hence $(c_g)_i=-\int_Z\langle g,\kappa_i\rangle_S$ and
\begin{equation}
 A^{\rm hom}(\xi)-m_g=A^{\rm lin}\xi-d_g,\qquad
 (d_g)_i=\int_Z\bigl\langle g(\tau),(y_i+\kappa_i)|_S\bigr\rangle_S\dd\tau .
 \label{eq:dphi}
\end{equation}
The vector $d_g$ does not depend on the choice of origin or on the
normalization of $\kappa_i$, since $\langle g,1\rangle_S=0$.  In this case $g$
enters the limit only through its $\tau$-average and only through the
condition on $\Gamma_1$; in the nonlinear case it also enters $A^{\rm hom}$
through the cell problem \eqref{eq:cell-parametric}.
\end{remark}

\begin{remark}[Time regularity]\label{rem:time-regularity}
If $\partial_tu_0\in L^2(0,T;(H^1_{\Gamma_2}(\Om))')$ were known, the comparisons
$u_0+\lambda(w-u_0)$, $\lambda\downarrow0$, and the regularization of
Lemma~\ref{lem:exp-regularization} would give \eqref{eq:macroVI} with
$\partial_tu_0$ in place of $\partial_tv$, for every
$v\in L^2(0,T;H^1_{\Gamma_2}(\Om))$ with $\gamma_0v\leq\varphi$.  This
regularity does not follow from \eqref{eq:uniform-time}.
\end{remark}

\begin{corollary}\label{cor:macro}
Assume \((H_1)\)--\((H_8)\) and \((H_{10})\), let $u_\eps$ be the microscopic
solution, and let $u_0$ be the solution of \eqref{eq:macroVI}.  Then
$P_\eps u_\eps\to u_0$ strongly in $L^2(Q)$ as $\eps\to0$, without extracting
a subsequence, and \eqref{eq:compact-gradq}--\eqref{eq:compact-weighted} hold
for the whole family with $u_1=\pi(u_0,D_xu_0)$.
\end{corollary}

\begin{proof}
By Theorem~\ref{thm:decoupled}, every cluster point solves \eqref{eq:macroVI},
which has exactly one solution $u_0$.  If full convergence failed, there would
be $d>0$ and a sequence $\eps_n\downarrow0$ with
$\|P_{\eps_n}u_{\eps_n}-u_0\|_{L^2(Q)}\geq d$.  Weighted compactness gives a
strongly convergent subsequence; Theorem~\ref{thm:decoupled} identifies its
limit as the solution of \eqref{eq:macroVI}, namely $u_0$, a contradiction.
The same cluster-point argument in the weak two-scale topologies, together
with \eqref{eq:corrector-decoupled} and uniqueness of $\pi(u_0,D_xu_0)$,
proves the full-family gradient convergences.
\end{proof}

\section*{Declarations}

\noindent\textbf{Data availability.} Data sharing is not applicable to this
article, as no datasets were generated or analysed during the current study.

\smallskip
\noindent\textbf{Conflict of interest.} The authors declare that they have no
conflict of interest.

\smallskip
\noindent\textbf{Author contributions.} All authors contributed to the
analysis and to the writing of the manuscript, and read and approved the final
version.


\begin{thebibliography}{99}

\bibitem{AcerbiChiadoDalMasoPercivale1992}
E.~Acerbi, V.~Chiad\`o Piat, G.~Dal Maso and D.~Percivale,
\emph{An extension theorem from connected sets, and homogenization in general
periodic domains},
Nonlinear Anal. \textbf{18} (1992), no.~5, 481--496.

\bibitem{Allaire1992}
G.~Allaire,
\emph{Homogenization and two-scale convergence},
SIAM J. Math. Anal. \textbf{23} (1992), no.~6, 1482--1518,
\href{https://doi.org/10.1137/0523084}{doi:10.1137/0523084}.

\bibitem{AllaireDamlamianHornung1996}
G.~Allaire, A.~Damlamian and U.~Hornung,
\emph{Two-scale convergence on periodic surfaces and applications},
in \emph{Proceedings of the International Conference on Mathematical Modelling
of Flow through Porous Media}, A.~Bourgeat et al., eds.,
World Scientific, Singapore, 1996, 15--25.

\bibitem{AltLuckhaus1983}
H.~W. Alt and S.~Luckhaus,
\emph{Quasilinear elliptic--parabolic differential equations},
Math. Z. \textbf{183} (1983), 311--341,
\href{https://doi.org/10.1007/BF01176474}{doi:10.1007/BF01176474}.

\bibitem{Avila2024}
J.~Avila,
\emph{Homogenization of quasilinear problems with semilinear terms and
Signorini boundary conditions in perforated domains},
NoDEA Nonlinear Differential Equations Appl. \textbf{31} (2024), no.~4,
Paper No.~64, 19~pp.,
\href{https://doi.org/10.1007/s00030-024-00957-0}{doi:10.1007/s00030-024-00957-0}.

\bibitem{Avila2025}
J.~Avila,
\emph{Homogenization and corrector results of elliptic problems with Signorini
boundary conditions in perforated domains},
Ann. Mat. Pura Appl. (4) \textbf{204} (2025), no.~3, 879--901,
\href{https://doi.org/10.1007/s10231-024-01510-7}{doi:10.1007/s10231-024-01510-7}.

\bibitem{Beliaev2004}
A.~Beliaev,
\emph{Homogenization of a parabolic operator with Signorini boundary
conditions in perforated domains},
Asymptot. Anal. \textbf{40} (2004), no.~3--4, 255--268,
\href{https://doi.org/10.3233/ASY-2004-654}{doi:10.3233/ASY-2004-654}.

\bibitem{BensoussanLionsPapanicolaou1978}
A.~Bensoussan, J.-L. Lions and G.~Papanicolaou,
\emph{Asymptotic Analysis for Periodic Structures},
North-Holland, Amsterdam, 1978.

\bibitem{Bonafede1993}
S.~Bonafede,
\emph{Quasilinear degenerate elliptic variational inequalities with
discontinuous coefficients},
Comment. Math. Univ. Carolin. \textbf{34} (1993), no.~1, 55--61.

\bibitem{Brezis1972}
H.~Br\'ezis,
\emph{Probl\`emes unilat\'eraux},
J. Math. Pures Appl. (9) \textbf{51} (1972), 1--168.

\bibitem{Brezis1973}
H.~Br\'ezis,
\emph{Op\'erateurs maximaux monotones et semi-groupes de contractions dans les
espaces de Hilbert},
North-Holland Mathematics Studies 5, North-Holland, Amsterdam, 1973.

\bibitem{Carrillo1999}
J.~Carrillo,
\emph{Entropy solutions for nonlinear degenerate problems},
Arch. Ration. Mech. Anal. \textbf{147} (1999), 269--361,
\href{https://doi.org/10.1007/s002050050152}{doi:10.1007/s002050050152}.

\bibitem{CioranescuGrisoDamlamian2012}
D.~Cioranescu, A.~Damlamian, P.~Donato, G.~Griso and R.~Zaki,
\emph{The periodic unfolding method in domains with holes},
SIAM J. Math. Anal. \textbf{44} (2012), no.~2, 718--760,
\href{https://doi.org/10.1137/100817942}{doi:10.1137/100817942}.

\bibitem{CioranescuDonato1988}
D.~Cioranescu and P.~Donato,
\emph{Homog\'en\'eisation du probl\`eme de Neumann non homog\`ene dans des
ouverts perfor\'es},
Asymptot. Anal. \textbf{1} (1988), no.~2, 115--138,
\href{https://doi.org/10.3233/ASY-1988-1203}{doi:10.3233/ASY-1988-1203}.

\bibitem{CioranescuSaintJeanPaulin1979}
D.~Cioranescu and J.~Saint Jean Paulin,
\emph{Homogenization in open sets with holes},
J. Math. Anal. Appl. \textbf{71} (1979), 590--607.

\bibitem{ClarkPacker1999}
G.~W. Clark and L.~A. Packer,
\emph{Two-scale homogenization of non-linear degenerate evolution equations},
J. Math. Anal. Appl. \textbf{238} (1999), no.~1, 316--328,
\href{https://doi.org/10.1006/jmaa.1999.6547}{doi:10.1006/jmaa.1999.6547}.

\bibitem{ConcaMuratTimofte2003}
C.~Conca, F.~Murat and C.~Timofte,
\emph{A generalized strange term in Signorini's type problems},
ESAIM Math. Model. Numer. Anal. \textbf{37} (2003), no.~5, 773--805,
\href{https://doi.org/10.1051/m2an:2003055}{doi:10.1051/m2an:2003055}.

\bibitem{DiazPodolskiyShaposhnikova2026}
J.~I. D\'iaz, A.~V. Podolskyi and T.~A. Shaposhnikova,
\emph{Homogenization of dynamic Signorini-type problems on critically
oscillating boundaries under time-periodic forcing},
arXiv:2607.03141v2 [math.AP] (2026),
\href{https://doi.org/10.48550/arXiv.2607.03141}{doi:10.48550/arXiv.2607.03141}.

\bibitem{DiazShaposhnikovaPodolskiy2024}
J.~I.~D\'iaz, T.~A.~Shaposhnikova and A.~V.~Podolskiy,
\emph{Strange non-local operators homogenizing the Poisson equation with
 dynamical unilateral boundary conditions: asymmetric particles of critical size},
Electron. J. Differential Equations \textbf{2024} (2024), no.~03, 1--29,
\href{https://doi.org/10.58997/ejde.2024.03}{doi:10.58997/ejde.2024.03}.

\bibitem{DiBenedettoShowalter1982}
E.~DiBenedetto and R.~E. Showalter,
\emph{A pseudo-parabolic variational inequality and Stefan problem},
Nonlinear Anal. \textbf{6} (1982), no.~3, 279--291,
\href{https://doi.org/10.1016/0362-546X(82)90095-5}{doi:10.1016/0362-546X(82)90095-5}.

\bibitem{DouanlaSteklov2013}
H.~Douanla,
\emph{Homogenization of Steklov spectral problems with indefinite density
function in perforated domains},
Acta Appl. Math. \textbf{123} (2013), 261--284,
\href{https://doi.org/10.1007/s10440-012-9765-4}{doi:10.1007/s10440-012-9765-4}.

\bibitem{DouanlaTetsadjio2017}
H.~Douanla and E.~Tetsadjio,
\emph{Reiterated homogenization of hyperbolic--parabolic equations in domains
with tiny holes},
Electron. J. Differential Equations \textbf{2017} (2017), no.~59, 1--22.

\bibitem{DouanlaWoukeng2015}
H.~Douanla and J.~L.~Woukeng,
\emph{Homogenization of reaction--diffusion equations in fractured porous media},
Electron. J. Differential Equations \textbf{2015} (2015), no.~253, 1--23.

\bibitem{DuvautLions1972}
G.~Duvaut and J.-L. Lions,
\emph{Les in\'equations en m\'ecanique et en physique},
Dunod, Paris, 1972.

\bibitem{FisherLiu2002}
B.~Fisher and Z.~Liu,
\emph{On the solvability of double degenerate quasilinear parabolic equations},
Acta Math. Hungar. \textbf{96} (2002), no.~1--2, 117--124,
\href{https://doi.org/10.1023/A:1015669813154}{doi:10.1023/A:1015669813154}.

\bibitem{FlodenHolmbomOlssonPersson2010}
L.~Flod\'en, A.~Holmbom, M.~Olsson and J.~Persson,
\emph{Very weak multiscale convergence},
Appl. Math. Lett. \textbf{23} (2010), 1170--1173,
\href{https://doi.org/10.1016/j.aml.2010.05.005}{doi:10.1016/j.aml.2010.05.005}.

\bibitem{Galdi2011}
G.~P. Galdi,
\emph{An Introduction to the Mathematical Theory of the Navier--Stokes
Equations: Steady-State Problems}, 2nd ed.,
Springer Monographs in Mathematics, Springer, New York, 2011,
\href{https://doi.org/10.1007/978-0-387-09620-9}{doi:10.1007/978-0-387-09620-9}.

\bibitem{Holmbom1997}
A.~Holmbom,
\emph{Homogenization of parabolic equations: an alternative approach and some
corrector-type results},
Appl. Math. \textbf{42} (1997), no.~5, 321--343,
\href{https://doi.org/10.1023/A:1023049608047}{doi:10.1023/A:1023049608047}.

\bibitem{HopkerBohm2014}
M.~H\"opker and M.~B\"ohm,
\emph{A note on the existence of extension operators for Sobolev spaces on
periodic domains},
C. R. Math. Acad. Sci. Paris \textbf{352} (2014), no.~10, 807--810,
\href{https://doi.org/10.1016/j.crma.2014.09.002}{doi:10.1016/j.crma.2014.09.002}.

\bibitem{Huang2006}
Y.~Huang, N.~Su and X.~Zhang,
\emph{Homogenization of degenerate quasilinear parabolic equations with
periodic structure},
Asymptot. Anal. \textbf{48} (2006), no.~1--2, 77--89,
\href{https://doi.org/10.3233/ASY-2006-757}{doi:10.3233/ASY-2006-757}.

\bibitem{Kuttler1984}
K.~L.~Kuttler, Jr.,
\emph{Degenerate variational inequalities of evolution},
Nonlinear Anal. \textbf{8} (1984), no.~8, 837--850,
\href{https://doi.org/10.1016/0362-546X(84)90106-8}{doi:10.1016/0362-546X(84)90106-8}.

\bibitem{KuttlerShillor1999}
K.~L.~Kuttler and M.~Shillor,
\emph{Set-valued pseudomonotone maps and degenerate evolution inclusions},
Commun. Contemp. Math. \textbf{1} (1999), no.~1, 87--123,
\href{https://doi.org/10.1142/S0219199799000067}{doi:10.1142/S0219199799000067}.

\bibitem{Le2011}
V.~K.~Le,
\emph{A range and existence theorem for pseudomonotone perturbations of maximal
monotone operators},
Proc. Amer. Math. Soc. \textbf{139} (2011), no.~5, 1645--1658,
\href{https://doi.org/10.1090/S0002-9939-2010-10594-4}{doi:10.1090/S0002-9939-2010-10594-4}.

\bibitem{Lions1969}
J.-L. Lions,
\emph{Quelques m\'ethodes de r\'esolution des probl\`emes aux limites non
lin\'eaires},
Dunod, Gauthier-Villars, Paris, 1969.


\bibitem{LiuLiu2013}
G.-F.~Liu and Y.-L.~Liu,
\emph{On double degenerate quasilinear parabolic variational inequalities},
Acta Math. Appl. Sin. Engl. Ser. \textbf{29} (2013), 861--868,
\href{https://doi.org/10.1007/s10255-013-0263-x}{doi:10.1007/s10255-013-0263-x}.

\bibitem{Liu2000}
Z.~Liu,
\emph{On the solvability of degenerate quasilinear parabolic equations of
second order},
Acta Math. Sin. (Engl. Ser.) \textbf{16} (2000), no.~2, 313--324,
\href{https://doi.org/10.1007/s101140000052}{doi:10.1007/s101140000052}.

\bibitem{Lobkova2019}
T.~Lobkova,
\emph{Homogenization of linear parabolic equations with a certain resonant
matching between rapid spatial and temporal oscillations in periodically
perforated domains},
Acta Math. Appl. Sin. Engl. Ser. \textbf{35} (2019), no.~2, 340--358,
\href{https://doi.org/10.1007/s10255-019-0810-1}{doi:10.1007/s10255-019-0810-1}.

\bibitem{Maitre2003}
E.~Maitre,
\emph{On a nonlinear compactness lemma in $L^p(0,T;B)$},
Int. J. Math. Math. Sci. \textbf{2003} (2003), no.~27, 1725--1730,
\href{https://doi.org/10.1155/S0161171203106175}{doi:10.1155/S0161171203106175}.

\bibitem{MelnykNakvasiuk2012}
T.~A. Mel'nyk and Iu.~A. Nakvasiuk,
\emph{Homogenization of a parabolic Signorini boundary value problem in a
thick plane junction},
J. Math. Sci. \textbf{181} (2012), no.~5, 613--631,
\href{https://doi.org/10.1007/s10958-012-0708-4}{doi:10.1007/s10958-012-0708-4}.

\bibitem{MielkeTimofte2007}
A.~Mielke and A.~M. Timofte,
\emph{Two-scale homogenization for evolutionary variational inequalities via
the energetic formulation},
SIAM J. Math. Anal. \textbf{39} (2007), no.~2, 642--668,
\href{https://doi.org/10.1137/060672790}{doi:10.1137/060672790}.

\bibitem{MigorskiDudek2025}
S.~Mig\'orski and S.~Dudek,
\emph{A class of evolutionary history-dependent variational-hemivariational
inequalities with strong and weak solutions},
SeMA J. \textbf{82} (2025), 179--203,
\href{https://doi.org/10.1007/s40324-024-00367-3}{doi:10.1007/s40324-024-00367-3}.

\bibitem{Moussa2016}
A.~Moussa,
\emph{Some variants of the classical Aubin--Lions lemma},
J. Evol. Equ. \textbf{16} (2016), 65--93,
\href{https://doi.org/10.1007/s00028-015-0293-3}{doi:10.1007/s00028-015-0293-3}.

\bibitem{NandakumaranRajesh2001}
A.~K. Nandakumaran and M.~Rajesh,
\emph{Homogenization of a nonlinear degenerate parabolic differential equation},
Electron. J. Differential Equations \textbf{2001} (2001), no.~17, 1--19.

\bibitem{NeussRadu1996}
M.~Neuss-Radu,
\emph{Some extensions of two-scale convergence},
C. R. Acad. Sci. Paris S\'er. I Math. \textbf{322} (1996), no.~9, 899--904.

\bibitem{Nguetseng1989}
G.~Nguetseng,
\emph{A general convergence result for a functional related to the theory of
homogenization},
SIAM J. Math. Anal. \textbf{20} (1989), no.~3, 608--623,
\href{https://doi.org/10.1137/0520043}{doi:10.1137/0520043}.

\bibitem{Otto1996}
F.~Otto,
\emph{$L^1$-contraction and uniqueness for quasilinear elliptic--parabolic
equations},
J. Differential Equations \textbf{131} (1996), 20--38,
\href{https://doi.org/10.1006/jdeq.1996.0155}{doi:10.1006/jdeq.1996.0155}.

\bibitem{Pastukhova2001}
S.~E.~Pastukhova,
\emph{Homogenization of a mixed problem with Signorini condition for an
elliptic operator in a perforated domain},
Sb. Math. \textbf{192} (2001), no.~2, 245--260,
\href{https://doi.org/10.1070/SM2001v192n02ABEH000544}{doi:10.1070/SM2001v192n02ABEH000544}.

\bibitem{PodolskiyShaposhnikova2022}
A.~V.~Podolskiy and T.~A.~Shaposhnikova,
\emph{Homogenization of a parabolic equation in a perforated domain with a
 unilateral dynamic boundary condition: critical case},
Contemporary Mathematics. Fundamental Directions \textbf{68} (2022), no.~4,
671--685, \href{https://doi.org/10.22363/2413-3639-2022-68-4-671-685}{doi:10.22363/2413-3639-2022-68-4-671-685}.

\bibitem{Ptashnyk2019}
M.~Ptashnyk,
\emph{Homogenization of some degenerate pseudoparabolic variational
inequalities},
J. Math. Anal. Appl. \textbf{469} (2019), 44--75,
\href{https://doi.org/10.1016/j.jmaa.2018.08.047}{doi:10.1016/j.jmaa.2018.08.047}.

\bibitem{Sandrakov2005}
G.~V.~Sandrakov,
\emph{Homogenization of variational inequalities for non-linear diffusion
problems in perforated domains},
Izv. Math. \textbf{69} (2005), no.~5, 1035--1059,
\href{https://doi.org/10.1070/IM2005v069n05ABEH002287}{doi:10.1070/IM2005v069n05ABEH002287}.

\bibitem{ScorzaDragoni1948}
G.~Scorza Dragoni,
\emph{Un teorema sulle funzioni continue rispetto ad una e misurabili
rispetto ad un'altra variabile},
Rend. Sem. Mat. Univ. Padova \textbf{17} (1948), 102--106.

\bibitem{Showalter1997}
R.~E. Showalter,
\emph{Monotone Operators in Banach Space and Nonlinear Partial Differential
Equations},
Mathematical Surveys and Monographs 49, American Mathematical Society,
Providence, RI, 1997,
\href{https://doi.org/10.1090/surv/049}{doi:10.1090/surv/049}.

\bibitem{Simon1986}
J.~Simon,
\emph{Compact sets in the space $L^p(0,T;B)$},
Ann. Mat. Pura Appl. (4) \textbf{146} (1986), 65--96,
\href{https://doi.org/10.1007/BF01762360}{doi:10.1007/BF01762360}.

\bibitem{Woukeng2010}
J.-L.~Woukeng,
\emph{Reiterated homogenization of nonlinear pseudo monotone degenerate
parabolic operators},
Commun. Math. Anal. \textbf{9} (2010), no.~2, 98--129.

\bibitem{Yosifian2001}
G.~A. Yosifian,
\emph{Homogenization of some contact problems for the system of elasticity in
perforated domains},
Rend. Sem. Mat. Univ. Padova \textbf{105} (2001), 37--64.

\bibitem{ZhangHuang2005}
X.-Y.~Zhang and Y.~Huang,
\emph{Homogenization for degenerate quasilinear parabolic equations of second
order},
Acta Math. Appl. Sin. Engl. Ser. \textbf{21} (2005), 93--100,
\href{https://doi.org/10.1007/s10255-005-0219-x}{doi:10.1007/s10255-005-0219-x}.

\bibitem{Zhikov1998}
V.~V. Zhikov,
\emph{Weighted Sobolev spaces},
Sb. Math. \textbf{189} (1998), no.~8, 1139--1170,
\href{https://doi.org/10.1070/SM1998v189n08ABEH000344}{doi:10.1070/SM1998v189n08ABEH000344}.
\end{thebibliography}
\end{document}